\documentclass{amsart}
\usepackage{amsmath, amsthm, amssymb, graphicx, tikz, zref} % Required for inserting images
\usetikzlibrary{arrows}
\usepackage{enumitem}
\newtheorem{theorem}{Theorem}[section]

\newtheorem{lemma}[theorem]{Lemma}

\newtheorem{definition}[theorem]{Definition}

\usepackage[margin=1in]{geometry} 
\usepackage{tikz-cd}
\usepackage{placeins}
\usepackage{etoolbox}
\pretocmd{\section}{\FloatBarrier}{}{}
\pretocmd{\subsection}{\FloatBarrier}{}{}

\usepackage{pdfpages}
\usepackage{hyperref}
\usepackage{faktor}

\usepackage{tikz-cd}
\usepackage{mathabx,epsfig}

\hypersetup{
    colorlinks,
    linkcolor={red!50!black},
    citecolor={blue!50!black},
    urlcolor={blue!80!black}
}
\newcommand{\R}{\mathbb{R}}

\newcommand{\PP}{\mathbb{P}}\newcommand{\Z}{\mathbb{Z}}\newcommand{\C}{\mathbb{C}}
\newcommand{\Ext}{\mathrm{Ext}}

\newcommand{\OO}{\mathcal{O}}

\usepackage{makecell}

\title{Exceptional Collections of Line Bundles for Smooth Toric Fano Surfaces and Threefolds}
\author{Reginald Anderson}
\email{reginala@uci.edu}
\date{May 2023}

\begin{document}

\maketitle
\begin{abstract} The cellular resolution of the diagonal given by Bayer-Popescu-Sturmfels for unimodular projective toric varieties yields a full, strong exceptional collection of line bundles on unimodular projective toric surfaces. The Hanlon-Hicks-Lazarev resolution of the diagonal yields a full, strong exceptional collection of line bundles for 16 of the 18 smooth toric Fano threefolds. 

\end{abstract}

\section{Introduction} We work over an algebraically closed field $k$ of characteristic $0$. Consider the question of whether two varieties $X$ and $Y$ are isomorphic. For affine varieties, this is equivalent to asking whether their affine coordinate rings $A(X)$ and $A(Y)$ are isomorphic. For projective varieties, the homogeneous coordinate rings $S(X)$ and $S(Y)$ are dependent upon an embedding into $\PP^n$, and thus are no longer isomorphism invariants. Though explicitly giving an isomorphism between varieties $X$ and $Y$ which are isomorphic can be difficult in general, Bondal-Orlov proved that if $X$ is a smooth projective variety which is Fano or general type, the bounded derived category of coherent sheaves $D^b_{Coh}(X)$ is an isomorphism invariant \cite{Bondal_2001}. 

Let $\mathcal{D}$ be a $k$-linear triangulated category. We recall the following definitions. 

\begin{definition} An object $E$ of $\mathcal{D}$ is \textbf{exceptional} if \[ Hom_\mathcal{D}(E,E[\ell]) = \begin{cases} k & \ell = 0, \\ 0 & \text{ else. } \end{cases}\]  \end{definition}
\begin{definition} An ordered collection $\mathcal{E}=(E_0,E_1,\dots,E_{N-1})$ of exceptional objects is \textbf{exceptional} if \[ Hom_\mathcal{D}(E_i, E_j[\ell]) = 0 \] for $0\leq j<i\leq N-1$ and all $\ell\in\mathbb{Z}$. \end{definition}

\begin{definition} An exceptional collection $\mathcal{E}$ is \textbf{strong} if \[ Hom_\mathcal{D}(E_i, E_j[\ell]) = 0 \] for $\ell\neq 0$ and all $0\leq i,j\leq N-1$. \end{definition}
\begin{definition} An exceptional collection $\mathcal{E}$ is \textbf{full} if the smallest triangulated subcategory generated by $\mathcal{E}$ is all of $\mathcal{D}$.  \end{definition}

Now let X be a smooth projective Fano variety over k, and write $D^b(X) = D^b_{Coh}(X)$. By the Bondal–Orlov reconstruction theorem \cite{Bondal_2001}, $D^b(X)$ determines $X$. In the presence of a strong, full exceptional collection $\mathcal{E}$, call $G = \bigoplus_i E_i$ a generator (or ``tilting object"), and let $A = End(G)$. Bondal showed in \cite{BondalRepnAssocAlg} that there is an exact equivalence \[ D^b(X) \stackrel{\sim}{\rightarrow} D^b(\text{mod}-A) \] from $D^b(X)$ to right modules over the endomorphism algebra $A$ induced from taking $RHom(G, -)$ as $k$-linear categories. That is, the presence of a strong, full exceptional collection allows for a nicest possible presentation of $D^b(X)$. King conjectured in unpublished notes that any smooth projective toric variety has a strong, full exceptional collection of line bundles. King's conjecture was proven false \cite{hille2006counterexample} \cite{Efimov2014}, and is currently checked on a case-by-case basis whether a given smooth projective toric variety has a strong, full exceptional collection of line bundles, and if so, how to construct such a collection. We now specialize to smooth projective toric Fano varieties of dimensions two and three. 

In this paper, we construct a strong, full exceptional collection of line bundles for unimodular projective toric surfaces. For smooth complete toric Fano threefolds --- which are automatically projective --- Bernardi–Tirabassi \cite{bernardi2010derivedcategoriestoricfano} constructed the relevant strong exceptional collections and obtained fullness conditional on Bondal’s conjecture; Uehara \cite{uehara2015exceptionalcollectionstoricfano} subsequently proved the required generation statement and thereby established unconditionally that every smooth toric Fano threefold over $\C$ admits a full strong exceptional collection of line bundles. We therefore examine explicitly whether particular resolutions of the diagonal yield full strong exceptional collections of line bundles. Since the Bayer–Popescu–Sturmfels resolution agrees with the Hanlon–Hicks–Lazarev resolution for smooth projective toric Fano surfaces, we concentrate on the latter and ask whether it yields such a collection for smooth projective toric Fano varieties of dimensions two and three.

Bayer-Popescu-Sturmfels \cite{bayer-popescu-sturmfels} gave a cellular resolution $(\mathcal{F}^\bullet_{\mathcal{H}_L/L}, \partial)$ of the diagonal for a proper subclass of smooth toric projective varieties whose lattice of principal divisors is unimodular which we recall in Section~\ref{sec:3.3}. A floor-function labeling extends this cellular construction to smooth projective toric varieties with centrally symmetric fans; without this additional hypothesis, the undeformed complex need not resolve the diagonal \cite{anderson2023resolutiondiagonalsmoothtoric}. A related construction for global-quotient toric Deligne--Mumford stacks arising as quotients of smooth toric varieties by finite abelian groups is given in \cite{anderson2023resolutiondiagonaltoricdelignemumford}; these stacks correspond to simplicial toric varieties. There is much recent work on virtual resolutions of the diagonal and of toric subvarieties of smooth projective toric varieties, including work of Hanlon-Hicks-Lazarev \cite{hanlon2023resolutions}, Bruce-Cranton Heller-Sayrafi \cite{bruce2022bounds}, and Brown-Erman \cite{brown2023short}, using the notion of a virtual resolution introduced by Berkesch-Erman-Smith \cite{BerkeschErmanSmith2020}. This paper focuses on the following theorems.

\begin{theorem}\label{thm: MainTheorem1} The Bayer-Popescu-Sturmfels cellular resolution $(\mathcal{F}^\bullet_{\mathcal{H}_L/L}, \partial) \rightarrow \mathcal{O}_\Delta$ yields a strong full, exceptional collection of line bundles for any unimodular projective toric surface. \end{theorem}

Here, we give a classification of smooth projective toric surfaces which are unimodular in the sense of Bayer-Popescu-Sturmfels, and show that this classification agrees with the 5 smooth toric Fano surfaces given in Craw \cite{Craw2007EXPLICITMF}. For smooth projective toric Fano threefolds, we give the following result.

\begin{theorem}\label{thm: MainTheorem2} For 16 of the 18 smooth, projective toric Fano 3-folds, the Hanlon-Hicks-Lazarev resolution of the diagonal yields a full, strong, exceptional collection of line bundles. \end{theorem} 

Theorem~\ref{thm: MainTheorem2} agrees with Bondal's observation \cite{bondalOberwolfach} that for all smooth toric Fano threefolds except for two, his construction yields a full strong exceptional collection of line bundles.

\section*{Acknowledgements} 
The author is grateful to Gabriel Kerr for many helpful discussions on derived categories in algebraic geometry. The author also thanks Jay Yang for help in creating Macaulay2 code for the Hanlon-Hicks-Lazarev resolution of a toric subvariety, as well as Jesse Huang for many helpful discussions on full, strong exceptional collections. Lastly, the author also thanks an anonymous reviewer for many helpful comments and suggestions for revision.

\section*{AI acknowledgements}
ChatGPT 5.6 Sol Ultra was used to re-format this paper, summarize previous results, and improve or re-create images in TeX.

\section{Preliminaries and prior results}
\label{sec:3.3}
Here, we follow the conventions of \cite{anderson2023resolutiondiagonalsmoothtoric} and \cite{anderson2023resolutiondiagonaltoricdelignemumford}, which follow the conventions of Bayer-Sturmfels \cite{bayer-sturmfels} and Bayer-Popescu-Sturmfels \cite{bayer-popescu-sturmfels} to resolve the diagonal for a unimodular toric variety. Bayer-Popescu-Sturmfels give a cellular resolution of the Lawrence ideal $J_L$ corresponding to a unimodular lattice $L$. The lattice $L$ given by the image of $B$ in the fundamental exact sequence \cite{C-L-S}
\begin{equation}\label{eqn: fundexactseq}
0 \rightarrow M \stackrel{B}{\rightarrow} \Z^{|\Sigma(1)|} \stackrel{\pi}{\rightarrow} \mathrm{Cl}(X_\Sigma) \rightarrow 0,   
\end{equation}  
describes the lattice of principal divisors of the projective toric variety $X_\Sigma$. Here $M$ gives the lattice of characters for the dense open torus of $X_\Sigma$, and we conflate $B$ as a matrix with the homomorphism it represents. Let $n$ denote $|\Sigma(1)|$. We adopt Bayer-Popescu-Sturmfels' notion of unimodularity for the lattice $L$ and toric variety $X_\Sigma$. 
\begin{definition} We call $L$ and $X_\Sigma$ \textbf{unimodular} when $B$ has linearly independent columns and every maximal minor of $B$ lies in the set $\{0,\pm1\}$. \end{definition} 
Throughout this paper, ``unimodular" refers to this restrictive condition of Bayer-Popescu-Sturmfels. Unimodular toric varieties form a proper subclass of smooth toric varieties. Let $X_\Sigma$ have homogeneous coordinate ring $R = \C[x_\rho \text{ }|\text{ } \rho \in \Sigma(1)]$ graded by $\mathrm{Cl}(X_\Sigma)$ via Equation~\ref{eqn: fundexactseq}, and $X_\Sigma \times X_\Sigma$ have homogeneous coordinate ring $S = R\otimes_\C R$. We use the lattice $L$ to construct a hyperplane arrangement $\mathcal{H}_L$ by intersecting $\R L$, the real span of the lattice $L = \text{Im }(B) \subset \Z^n$, with all lattice translates of the coordinate hyperplanes $\{x_i=j\}, 1\leq i \leq n, j\in \Z$. This gives an infinite cellular complex $\R L \cap \mathcal{H}_L$, which we quotient by the action of the lattice $L$ to yield the finite cellular complex $\faktor{\R L \cap \mathcal{H}_L}{L}$, which supports the minimal $\mathrm{Cl}(X_\Sigma) \times \mathrm{Cl}(X_\Sigma)$-graded resolution of 

\begin{align*} J_L = \left<\textbf{x}^\textbf{a}\textbf{y}^\textbf{b} - \textbf{x}^\textbf{b}\textbf{y}^\textbf{a} \text{ }|\text{ }\textbf{a}-\textbf{b} \in L \right> = \left< \textbf{x}^\textbf{u}\textbf{y}^{\textbf{v}} - \textbf{x}^\textbf{v} \textbf{y}^\textbf{u} \text{ }|\text{ } \pi(\textbf{u}) = \pi(\textbf{v}) \text{ in } \mathrm{Cl}(X_\Sigma) \right> \subset S
\end{align*} 

corresponding to $X_\Sigma \cong \Delta \subset X_\Sigma \times X_\Sigma$.

\subsection{The Bayer--Popescu--Sturmfels resolution}
\label{subsec:BPS-construction}

We now give the construction of the Bayer--Popescu--Sturmfels
resolution used below. Retain the notation of
Equation~\ref{eqn: fundexactseq}, write
\[
d=\dim X_\Sigma=\operatorname{rank}(L),
\qquad
n=\lvert\Sigma(1)\rvert,
\]
and identify the lattice of principal divisors with
\[
L=\operatorname{Im}(B)=\ker(\pi)\subseteq\Z^n.
\]
Enumerate the rays as
\(\Sigma(1)=\{\rho_1,\ldots,\rho_n\}\), and write
\(D_i=D_{\rho_i}\) for the corresponding torus-invariant prime
divisor. 
Write
\[
S=\C[x_1,\ldots,x_n,y_1,\ldots,y_n]
\]
for the Cox ring of \(X_\Sigma\times X_\Sigma\), with its
\(\mathrm{Cl}(X_\Sigma)\times\mathrm{Cl}(X_\Sigma)\)-grading given by
\[
\deg(x_i)=([D_i],0),
\qquad
\deg(y_i)=(0,[D_i]).
\]
For \(\mathbf{u}=(u_1,\ldots,u_n)\in\Z^n\), set
\[
\mathbf{x}^{\mathbf{u}}
=
x_1^{u_1}\cdots x_n^{u_n},
\qquad
\mathbf{y}^{\mathbf{u}}
=
y_1^{u_1}\cdots y_n^{u_n},
\]
as Laurent monomials whenever some \(u_i\) is negative. If
\(\mathbf{u}=\mathbf{u}^{+}-\mathbf{u}^{-}\), where
\(\mathbf{u}^{+},\mathbf{u}^{-}\in\Z_{\geq0}^n\) have disjoint support, the
Lawrence ideal of \(L\) is
\[
J_L
=
\left\langle
\mathbf{x}^{\mathbf{u}^{+}}\mathbf{y}^{\mathbf{u}^{-}}
-
\mathbf{x}^{\mathbf{u}^{-}}\mathbf{y}^{\mathbf{u}^{+}}
\ \middle|\
\mathbf{u}\in L
\right\rangle
\subseteq S.
\]
Equivalently,
\[
J_L
=
\left\langle
\mathbf{x}^{\mathbf{a}}\mathbf{y}^{\mathbf{b}}
-
\mathbf{x}^{\mathbf{b}}\mathbf{y}^{\mathbf{a}}
\ \middle|\
\mathbf{a},\mathbf{b}\in\Z_{\geq0}^n,\ 
\mathbf{a}-\mathbf{b}\in L
\right\rangle.
\]
Since \(L=\ker(\pi)\), the last condition is equivalent to
\(\pi(\mathbf{a})=\pi(\mathbf{b})\). Bayer--Popescu--Sturmfels identify
this Lawrence ideal with the homogeneous ideal representing the
diagonal in Cox coordinates
\cite[Section~6, Proposition~6.1]{bayer-popescu-sturmfels}. More
precisely, it is the kernel of
\[
\Phi\colon
S\longrightarrow
\C[\mathrm{Cl}(X_\Sigma)]\otimes_{\C}R,
\qquad
\mathbf{x}^{\mathbf{a}}\mathbf{y}^{\mathbf{b}}
\longmapsto
[\pi(\mathbf{a})]\otimes
\mathbf{x}^{\mathbf{a}+\mathbf{b}},
\]
where \([g]\) denotes the basis monomial of the group algebra
\(\C[\mathrm{Cl}(X_\Sigma)]\) corresponding to
\(g\in\mathrm{Cl}(X_\Sigma)\). Its Cox sheafification on
\(X_\Sigma\times X_\Sigma\) satisfies
\(\widetilde{S/J_L}\cong\mathcal{O}_\Delta\).

We next construct the cellular complex, making precise the notation
fixed in the preceding paragraph. Let
\[
\R L=L\otimes_{\Z}\R\subseteq\R^n
\]
and let \(\mathcal{H}_L\) denote the ambient periodic arrangement of
all integral translates of the coordinate hyperplanes in \(\R^n\):
\[
\mathcal{H}_L
=
\bigcup_{i=1}^{n}\ \bigcup_{j\in\Z}
\left\{
\mathbf{p}=(p_1,\ldots,p_n)\in\R^n
\ \middle|\
p_i=j
\right\}.
\]
The restrictions of these hyperplanes to \(\R L\) induce a periodic
polyhedral cell decomposition of \(\R L\), which we denote by
\(\mathcal{A}_L\). Its cells are the relatively open faces of every
dimension cut out by the restricted hyperplanes. To retain the
notation used throughout the paper, we write
\[
\mathcal{A}_L/L
=
\faktor{\R L\cap\mathcal{H}_L}{L}.
\]
Thus the expression on the right denotes the quotient of the induced
cell decomposition. We retain the established shorthands
\(\mathcal{F}_{\mathcal{H}_L,\bullet}\) and
\(\mathcal{F}_{\mathcal{H}_L/L,\bullet}\) for the associated
algebraic complexes. Equivalently, if the rows of the
\(n\times d\) matrix \(B\) are
\(\mathbf{b}_1,\ldots,\mathbf{b}_n\), then under the identification
\(\R L\cong\R^d\), the arrangement consists of the hyperplanes
\[
H_{ij}
=
\left\{
\mathbf{q}\in\R^d
\ \middle|\
\mathbf{b}_i\cdot\mathbf{q}=j
\right\},
\qquad
1\leq i\leq n,\quad j\in\Z.
\]
The role of unimodularity is fundamental: the vertices of
\(\mathcal{A}_L\) are precisely the lattice points of \(L\) if and
only if \(L\) is unimodular
\cite[Proposition~2.2]{bayer-popescu-sturmfels}.

Henceforth in this subsection, assume that \(L\) is unimodular.

Let
\[
T
=
\C[x_1^{\pm1},\ldots,x_n^{\pm1},
y_1^{\pm1},\ldots,y_n^{\pm1}]
\]
and define the monomial \(S\)-submodule
\[
M_L
=
S\cdot
\left\{
\mathbf{x}^{\mathbf{u}}\mathbf{y}^{-\mathbf{u}}
\ \middle|\
\mathbf{u}\in L
\right\}
\subseteq T.
\]
The vertex \(\mathbf{u}\in L\) of \(\mathcal{A}_L\) is labeled by
\[
m_{\mathbf{u}}
=
\mathbf{x}^{\mathbf{u}}\mathbf{y}^{-\mathbf{u}}.
\]
For every nonempty face \(F\) of \(\mathcal{A}_L\), define \(m_F\) to
be the least common multiple of the labels of the vertices of \(F\),
where the least common multiple of Laurent monomials is taken by
taking the coordinatewise maximum of their exponent vectors. There is
also an intrinsic description of this label. If
\(\mathbf{p}\in\operatorname{relint}(F)\), then
\[
m_F
=
\mathbf{x}^{\lceil\mathbf{p}\rceil}
\mathbf{y}^{-\lfloor\mathbf{p}\rfloor},
\]
where the floor and ceiling are taken coordinatewise. This expression
is constant on \(\operatorname{relint}(F)\).

Choose orientations on the faces and let
\(\epsilon(F,G)\in\{0,\pm1\}\) be the corresponding incidence
function. Thus \(\epsilon(F,G)\neq0\) precisely when \(G\) is a facet
of \(F\). Define
\[
\mathcal{F}_{\mathcal{H}_L,i}
=
\bigoplus_{\substack{F\in\mathcal{A}_L\\ \dim F=i}}
S\bigl(-\deg_{\Z^{2n}}(m_F)\bigr)e_F.
\]
The cellular differential is
\[
\partial_i(e_F)
=
\sum_{\substack{G\subset F\\ \dim G=i-1}}
\epsilon(F,G)\frac{m_F}{m_G}e_G.
\]
Since \(G\subset F\), the Laurent monomial \(m_G\) divides \(m_F\), so
\(m_F/m_G\) is an ordinary monomial in \(S\). The augmentation in
degree zero sends
\[
e_{\mathbf{u}}
\longmapsto
\mathbf{x}^{\mathbf{u}}\mathbf{y}^{-\mathbf{u}}
\in M_L.
\]
Bayer--Popescu--Sturmfels prove that
\[
\left(
\mathcal{F}_{\mathcal{H}_L,\bullet},\partial
\right)
\longrightarrow M_L
\]
is a minimal \(\Z^{2n}\)-graded free \(S\)-resolution
\cite[Theorem~3.1]{bayer-popescu-sturmfels}. The exactness can be seen
directly from the cellular criterion of
\cite[Proposition~1.2]{bayer-sturmfels}. Indeed, for
\((\mathbf{a},\mathbf{b})\in\Z^{2n}\), the subcomplex consisting of
faces whose labels divide
\(\mathbf{x}^{\mathbf{a}}\mathbf{y}^{\mathbf{b}}\) is the
polyhedral subdivision of
\[
\left\{
\mathbf{p}\in\R L
\ \middle|\
-\mathbf{b}\leq\mathbf{p}\leq\mathbf{a}
\right\}.
\]
Whenever it is nonempty, this set is convex and hence contractible.
Moreover, distinct faces have distinct labels, which gives
minimality.

Translation by \(L\) preserves \(\mathcal{A}_L\). For
\(\mathbf{u}\in L\),
\[
F\longmapsto F+\mathbf{u},
\qquad
m_{F+\mathbf{u}}
=
\mathbf{x}^{\mathbf{u}}\mathbf{y}^{-\mathbf{u}}m_F.
\]
Introduce the Lawrence lifting
\[
\Lambda(L)
=
\left\{
(\mathbf{u},-\mathbf{u})\in\Z^{2n}
\ \middle|\
\mathbf{u}\in L
\right\}.
\]
The relevant quotient is the finite quotient face poset
\[
\faktor{\R L\cap\mathcal{H}_L}{L}.
\]
Its finiteness is proved in
\cite[Proposition~3.3 and Corollary~3.4]
{bayer-popescu-sturmfels}. To describe the induced algebraic
quotient, let
\[
S[L]
=
\bigoplus_{\mathbf{u}\in L}S\,z^{\mathbf{u}},
\qquad
z^{\mathbf{u}}z^{\mathbf{v}}
=
z^{\mathbf{u}+\mathbf{v}},
\]
where \(z^{\mathbf{u}}\) has degree
\((\mathbf{u},-\mathbf{u})\). Thus \(S[L]\) is the group algebra
\(S[\Lambda(L)]\) under the canonical identification
\(L\xrightarrow{\sim}\Lambda(L)\),
\(\mathbf{u}\mapsto(\mathbf{u},-\mathbf{u})\). It acts on the cellular
complex by
\[
z^{\mathbf{u}}e_F=e_{F+\mathbf{u}}.
\]
The compatible action on the augmentation module is
\[
z^{\mathbf{u}}\cdot
\bigl(\mathbf{x}^{\mathbf{v}}\mathbf{y}^{-\mathbf{v}}\bigr)
=
\mathbf{x}^{\mathbf{u}+\mathbf{v}}
\mathbf{y}^{-(\mathbf{u}+\mathbf{v})}.
\]
With an \(L\)-invariant choice of incidence signs,
\(\mathcal{F}_{\mathcal{H}_L,\bullet}\) is a finite-rank free
\(S[L]\)-complex. The quotient functor is obtained from
\[
S[L]\longrightarrow S,
\qquad
z^{\mathbf{u}}\longmapsto1,
\]
so put \(\Pi(-):=(-)\otimes_{S[L]}S\). This functor coarsens the
grading from \(\Z^{2n}\) to
\(\Z^{2n}/\Lambda(L)\). Although \(S\) is not flat over \(S[L]\),
\(\Pi\) is exact on the relevant graded categories: it is the
equivalence between \(\Z^{2n}\)-graded \(S[L]\)-modules and
\(\Z^{2n}/\Lambda(L)\)-graded \(S\)-modules established in
\cite[Theorem~3.2]{bayer-sturmfels}. Moreover,
\(\Pi(M_L)\cong S/J_L\) by
\cite[Lemma~3.1]{bayer-sturmfels}. This is the lattice-module quotient
construction of \cite[Lemma~3.5 and Corollary~3.7]{bayer-sturmfels}.

Choose one representative of each \(L\)-orbit of faces. The resulting
finite complex is
\[
\mathcal{F}_{\mathcal{H}_L/L,\bullet}
=
\mathcal{F}_{\mathcal{H}_L,\bullet}
\otimes_{S[L]}S,
\]
with
\[
\mathcal{F}_{\mathcal{H}_L/L,i}
=
\bigoplus_{\substack{[F]\in\faktor{\R L\cap\mathcal{H}_L}{L}\\
\dim F=i}}
S\bigl(-[\deg_{\Z^{2n}}(m_F)]\bigr)e_{[F]}.
\]
Here the brackets denote the class in
\(\Z^{2n}/\Lambda(L)\). If \(F\) is a chosen representative and
\(G\) is an actual facet of \(F\), let \(G_0\) be the chosen
representative of the orbit of \(G\), and choose
\(\mathbf{u}_G\in L\) with \(G=G_0+\mathbf{u}_G\). Before applying the
quotient functor, the differential in the \(S[L]\)-basis of chosen
representatives is
\[
\partial_i(e_F)
=
\sum_{\substack{G\subset F\\ \dim G=i-1}}
\epsilon(F,G)\frac{m_F}{m_G}
z^{\mathbf{u}_G}e_{G_0}.
\]
After applying \(z^{\mathbf{u}}\mapsto1\), this becomes
\[
\partial_i(e_{[F]})
=
\sum_{\substack{G\subset F\\ \dim G=i-1}}
\epsilon(F,G)\frac{m_F}{m_G}e_{[G_0]}.
\]
In this formula \(G\) ranges over the actual facets of \(F\); if
several such facets lie in the same \(L\)-orbit, their contributions
are added in the same component. Thus, for example, the two endpoints
of an edge can become the same vertex in the quotient, in which case
their two monomial contributions form a binomial. A different choice
of orbit representatives gives an isomorphic graded complex.

The main result of Bayer--Popescu--Sturmfels is that, when \(L\) is
unimodular, the augmented quotient complex is the minimal free
resolution
\[
0
\longrightarrow
\mathcal{F}_{\mathcal{H}_L/L,d}
\longrightarrow\cdots\longrightarrow
\mathcal{F}_{\mathcal{H}_L/L,1}
\longrightarrow
\mathcal{F}_{\mathcal{H}_L/L,0}=S
\longrightarrow
S/J_L
\longrightarrow0
\]
\cite[Theorem~3.5]{bayer-popescu-sturmfels}.

It remains to pass from the Cox ring to
\(X_\Sigma\times X_\Sigma\). There is a well-defined coarsening map
\[
\Z^{2n}/\Lambda(L)
\longrightarrow
\mathrm{Cl}(X_\Sigma)\times\mathrm{Cl}(X_\Sigma),
\qquad
[(\mathbf{a},\mathbf{b})]
\longmapsto
(\pi(\mathbf{a}),\pi(\mathbf{b})).
\]
If
\[
m_F
=
\mathbf{x}^{\boldsymbol{\alpha}_F}
\mathbf{y}^{\boldsymbol{\beta}_F},
\]
put
\[
\lambda_F=\pi(\boldsymbol{\alpha}_F),
\qquad
\mu_F=\pi(\boldsymbol{\beta}_F).
\]
These classes do not depend on the representative of \([F]\), since
translation by \(\mathbf{u}\in L\) changes the two exponent vectors
by \(\mathbf{u}\) and \(-\mathbf{u}\), respectively, and
\(\pi(\mathbf{u})=0\). After coarsening the grading, the \(i\)-th term
is therefore
\[
\bigoplus_{\substack{[F]\in\faktor{\R L\cap\mathcal{H}_L}{L}\\
\dim F=i}}
S(-\lambda_F,-\mu_F).
\]
Because \(X_\Sigma\) is smooth, every divisor class determines a line
bundle. Cox sheafification is exact, since it is computed locally by
homogeneous localization followed by taking degree zero. It consequently
gives
\[
\widetilde{\mathcal{F}}_{\mathcal{H}_L/L,i}
=
\bigoplus_{\substack{[F]\in\faktor{\R L\cap\mathcal{H}_L}{L}\\
\dim F=i}}
\mathcal{O}_{X_\Sigma}(-\lambda_F)
\boxtimes
\mathcal{O}_{X_\Sigma}(-\mu_F).
\]
Using the identification
\(\widetilde{S/J_L}\cong\mathcal{O}_\Delta\) from
\cite[Section~6, Proposition~6.1]{bayer-popescu-sturmfels}, we obtain
the locally free resolution
\[
0
\longrightarrow
\widetilde{\mathcal{F}}_{\mathcal{H}_L/L,d}
\longrightarrow\cdots\longrightarrow
\widetilde{\mathcal{F}}_{\mathcal{H}_L/L,1}
\longrightarrow
\mathcal{O}_{X_\Sigma\times X_\Sigma}
\longrightarrow
\mathcal{O}_\Delta
\longrightarrow0.
\]
This is the Bayer--Popescu--Sturmfels resolution of the diagonal used
in the remainder of the paper. After reversing homological indices,
so that
\(\mathcal{F}^{-i}_{\mathcal{H}_L/L}
=\widetilde{\mathcal{F}}_{\mathcal{H}_L/L,i}\) for
\(0\leq i\leq d\), we denote this augmented locally free complex by
\[
\left(\mathcal{F}^{\bullet}_{\mathcal{H}_L/L},\partial\right)
\longrightarrow\mathcal{O}_\Delta.
\]
The line bundles on the first factor
that it produces are
\[
\mathcal{O}_{X_\Sigma}(-\lambda_F),
\qquad
[F]\in\faktor{\R L\cap\mathcal{H}_L}{L},
\]
and similarly the line bundles on the second factor are
\(\mathcal{O}_{X_\Sigma}(-\mu_F)\). When forming the collection used
below, repetitions are removed and an ordering is then specified.

\subsection{The Hanlon--Hicks--Lazarev resolution}
\label{subsec:hhl-resolution}

We next recall the resolution of a toric subvariety constructed by
Hanlon--Hicks--Lazarev. The collection of line bundles is defined in
\cite[Section~2.3, equations~(6)--(7)]{hanlon2023resolutions}, the
boundary morphisms are defined in
\cite[Proposition~2.14, Definition~2.15, and
Proposition~2.16]{hanlon2023resolutions}, and the complex itself is
constructed in
\cite[Sections~3.4--3.5 and Appendix~A]{hanlon2023resolutions}.
Exactness is \cite[Theorem~3.5]{hanlon2023resolutions}.
Hanlon--Hicks--Lazarev state their theorem more generally for a closed
toric substack of a toric stack covered by smooth stacky charts. We
record the scheme case, which already applies to every closed toric
subvariety of a smooth toric variety and is the form used below.

Let \(Z=X_{\Sigma_Z}\) be a smooth toric variety with cocharacter
lattice \(N_Z\), character lattice \(M_Z=N_Z^*\), and torus-invariant
prime divisors \(D_\rho\), indexed by \(\rho\in\Sigma_Z(1)\). Write
\(u_\rho\in N_Z\) for the primitive generator of \(\rho\). For
\(m\in (M_Z)_\R\), set
\[
D(m):=
\sum_{\rho\in\Sigma_Z(1)}
\left\lfloor-\langle m,u_\rho\rangle\right\rfloor D_\rho.
\]
If \(m_0\in M_Z\), then
\[
D(m+m_0)=D(m)-\operatorname{div}(\chi^{m_0}),
\]
so the isomorphism class of \(\mathcal{O}_Z(D(m))\) depends only on
\([m]\in (M_Z)_\R/M_Z\). The finite set
\[
\Theta_Z:=
\left\{
\mathcal{O}_Z(D(m)):\ [m]\in (M_Z)_\R/M_Z
\right\}
\subseteq\operatorname{Pic}(Z)
\]
is the Thomsen collection. Equivalently, if
\[
F_m(u_\rho):=\left\lceil\langle m,u_\rho\rangle\right\rceil,
\]
then \(D(m)=-\sum_\rho F_m(u_\rho)D_\rho\). Thus the collection is
computed directly from the primitive ray generators of the fan. It
also has an intrinsic description in terms of toric Frobenius:
\(\Theta_Z\) is the set of distinct line-bundle summands of
\((F_\ell)_*\mathcal{O}_Z\) for \(\ell\) sufficiently large and highly
divisible; see \cite[Theorem~5.1 and
Corollary~5.5]{hanlon2023resolutions}.

Now let
\[
\iota\colon Y=X_{\Sigma_Y}\hookrightarrow Z=X_{\Sigma_Z}
\]
be a closed toric immersion, and let \(N_Y\) and
\(M_Y=N_Y^*\) be the cocharacter and character lattices of \(Y\),
respectively. The immersion is induced by an injective map
\[
\varphi\colon N_Y\hookrightarrow N_Z
\]
with torsion-free cokernel, and the immersion condition on the fans is
\[
\Sigma_Y
=
\left\{
\varphi_\R^{-1}(\sigma):\ \sigma\in\Sigma_Z
\right\}.
\]
The dual map
\(\varphi^*\colon M_Z\to M_Y\) induces a map of real tori, and we put
\[
T^\iota
:=\ker\left(
(M_Z)_\R/M_Z\longrightarrow (M_Y)_\R/M_Y
\right).
\]
Since \(\operatorname{coker}(\varphi)\) is torsion-free,
\(\varphi^*\) is surjective. If \(K=\ker(\varphi^*)\), the inclusion
\(K_\R\hookrightarrow(M_Z)_\R\) therefore induces an isomorphism
\[
K_\R/K\stackrel{\sim}{\longrightarrow}T^\iota.
\]
In particular,
\[
\dim T^\iota=\dim Z-\dim Y=\operatorname{codim}_Z(Y).
\]
For every ray \(\rho\in\Sigma_Z(1)\) whose primitive generator does
not lie in \(\operatorname{im}(\varphi)\), define
\[
H_\rho
:=\left\{
[m]\in T^\iota:\ \langle m,u_\rho\rangle\in\Z
\right\}.
\]
Each \(H_\rho\) may be a disjoint union of translates of a
codimension-one subtorus. These sets form a toric hyperplane
arrangement on \(T^\iota\). Let \(S^\iota\) be the resulting
stratification. Its lift to \(K_\R\) is the periodic affine hyperplane
arrangement
\[
\langle m,u_\rho\rangle=j,
\qquad \rho\in\Sigma_Z(1),\quad
u_\rho\notin\operatorname{im}(\varphi),\quad j\in\Z.
\]

For a stratum \(\sigma\in S^\iota\), choose a lift
\(\widetilde\sigma\subset K_\R\) and a point
\(m_\sigma\in\widetilde\sigma\). The integers
\[
F_{\widetilde\sigma}(u_\rho)
:=\left\lceil\langle m_\sigma,u_\rho\rangle\right\rceil
\]
are constant on \(\widetilde\sigma\). Define
\[
\mathcal{L}_\sigma
:=\mathcal{O}_Z\!\left(
\sum_{\rho\in\Sigma_Z(1)}
\left\lfloor-\langle m_\sigma,u_\rho\rangle\right\rfloor D_\rho
\right).
\]
Changing the lift changes the displayed divisor by a principal
divisor, so \(\mathcal{L}_\sigma\) is well defined up to toric
isomorphism. In particular,
\(\mathcal{L}_\sigma\in\Theta_Z\).

Let \(\widehat S_Z\) denote the full periodic Bondal stratification
of \((M_Z)_\R\) induced by all affine hyperplanes
\[
\langle m,u_\rho\rangle=j,
\qquad \rho\in\Sigma_Z(1),\quad j\in\Z.
\]
For \(\widehat\sigma\in\widehat S_Z\), its ceiling function is
\[
F_{\widehat\sigma}(u_\rho)
:=\left\lceil\langle m,u_\rho\rangle\right\rceil,
\qquad m\in\widehat\sigma;
\]
this is independent of the chosen \(m\) in the stratum. We write
\(\tau<\sigma\) when \(\tau\) is contained in the closure of
\(\sigma\), and use the same convention for strata of
\(\widehat S_Z\).

For \(p\in T^\iota\), let \(\dim_S(p)\) denote the dimension of the
stratum containing \(p\). For strata \(\tau<\sigma\), let
\(\operatorname{EP}(\sigma,\tau)\) be the set of homotopy classes of
stratum-monotone paths \(\gamma\colon[0,1]\to T^\iota\) satisfying
\[
\gamma(0)\in\sigma,\qquad
\gamma(1)\in\tau,\qquad
t_1\leq t_2\Longrightarrow
\dim_S(\gamma(t_1))\geq\dim_S(\gamma(t_2)).
\]
With the strata as objects, constant paths as identities, and
concatenation as composition, these homotopy classes form the
exit-path category \(\operatorname{EP}(S^\iota)\).
These classes retain the different incidences which may be identified
after passing from the periodic arrangement in \(K_\R\) to the torus
\(K_\R/K\). For
\(\gamma\in\operatorname{EP}(\sigma,\tau)\), choose restricted lifts
\(\widetilde\sigma,\widetilde\tau\subset K_\R\) representing
\(\gamma\), together with strata
\(\widehat\sigma,\widehat\tau\in\widehat S_Z\), such that the
restricted lifts lie in
their intersections with \(K_\R\) and
\(\widehat\tau<\widehat\sigma\). Hanlon--Hicks--Lazarev prove that
\[
F_{\widehat\sigma}(u_\rho)
-F_{\widehat\tau}(u_\rho)\geq 0
\qquad\text{for every }\rho\in\Sigma_Z(1).
\]
Consequently, in Cox coordinates \(x_\rho\) on \(Z\), the monomial
\[
x_\gamma
:=\prod_{\rho\in\Sigma_Z(1)}
 x_\rho^{F_{\widehat\sigma}(u_\rho)
-F_{\widehat\tau}(u_\rho)}
\]
defines a morphism
\[
\delta_\gamma\colon
\mathcal{L}_\sigma\longrightarrow\mathcal{L}_\tau.
\]
This is the boundary morphism of
\cite[Definition~2.15]{hanlon2023resolutions}; equivalently, it is the
monomial section indexed by \(0\) of the line bundle associated to
\(F_{\widehat\tau}-F_{\widehat\sigma}\). It is independent of the
chosen lifts after the preceding toric identifications. For composable
exit paths \(\gamma_1\colon\sigma\to\tau\) and
\(\gamma_2\colon\tau\to\upsilon\), the boundary morphisms satisfy
\[
\delta_{\gamma_2}\circ\delta_{\gamma_1}
=\delta_{\gamma_2\circ\gamma_1};
\]
see \cite[Proposition~2.16 and Proposition~3.9]
{hanlon2023resolutions}. Thus the assignments
\[
\mathcal{O}^\iota(\sigma):=\mathcal{L}_\sigma,
\qquad
\mathcal{O}^\iota(\gamma):=\delta_\gamma
\]
define the coefficient functor
\(\mathcal{O}^\iota\colon\operatorname{EP}(S^\iota)\longrightarrow
\operatorname{Coh}(Z)\).

Choose an orientation of every stratum of \(S^\iota\). An exit path
\(\gamma\) from a stratum \(\sigma\) to a codimension-one boundary
stratum \(\tau\) receives the incidence sign
\(\operatorname{sgn}(\gamma)\in\{\pm1\}\): it is \(+1\) when the
boundary orientation induced from \(\sigma\), using the outward
direction determined by \(\gamma\), agrees with the chosen orientation
of \(\tau\), and it is \(-1\) otherwise. Define
\[
C_k(S^\iota,\mathcal{O}^\iota)
:=\bigoplus_{\substack{\sigma\in S^\iota\\
\dim\sigma=k}}\mathcal{L}_\sigma
\]
and define the component of the differential from the
\(\mathcal{L}_\sigma\)-summand to the
\(\mathcal{L}_\tau\)-summand by
\[
(d_k)^\tau_\sigma
:=\sum_{\gamma\in\operatorname{EP}(\sigma,\tau)}
\operatorname{sgn}(\gamma)\,\delta_\gamma.
\]
The length-two boundary paths are paired by the Morse-quiver
involution \cite[Definition A.1]{hanlon2023resolutions}, and the two paths in each pair induce the same monomial
morphism with opposite signs.
Therefore \(d_{k-1}d_k=0\); this is
\cite[Propositions~3.9--3.10 and Appendix~A,
equations~(18)--(19)]{hanlon2023resolutions}. The distinguished
zero-dimensional stratum \(\sigma_0\) containing the identity of
\(T^\iota\) is labeled by \(\mathcal{O}_Z\); the identity on this
summand defines the morphism
\[
\alpha^\iota\colon
\mathcal{O}_Z\longrightarrow C_0(S^\iota,\mathcal{O}^\iota)
\]
which is an epimorphism on homology. Theorem~3.5 of
\cite{hanlon2023resolutions} gives an augmentation
\[
0\longrightarrow C_c(S^\iota,\mathcal{O}^\iota)
\longrightarrow\cdots\longrightarrow
C_1(S^\iota,\mathcal{O}^\iota)
\longrightarrow C_0(S^\iota,\mathcal{O}^\iota)
\stackrel{\epsilon}{\longrightarrow}
\iota_*\mathcal{O}_Y\longrightarrow0,
\qquad c=\operatorname{codim}_Z(Y),
\]
where \(\epsilon\) is the composite
\[
C_0(S^\iota,\mathcal{O}^\iota)
\twoheadrightarrow
H_0\!\left(C_\bullet(S^\iota,\mathcal{O}^\iota)\right)
\stackrel{\sim}{\longrightarrow}
\iota_*\mathcal{O}_Y.
\]
The augmented complex is exact. This is the
Hanlon--Hicks--Lazarev resolution of the toric
subvariety \(Y\) in \(Z\). Notice that neither projectivity nor the
Fano condition is needed for this construction; smoothness of the
ambient toric variety is the relevant hypothesis.

We now specialize this construction to the diagonal in the setting
used in this paper. Let \(X=X_\Sigma\) be a smooth projective toric
variety of dimension \(d\), with cocharacter lattice \(N\) and
character lattice \(M=N^*\). For \(\rho\in\Sigma(1)\), write
\(u_\rho\in N\) for the primitive ray generator and \(D_\rho\) for
the corresponding torus-invariant prime divisor. Take
\[
Z=X\times X,
\qquad
Y=X,
\qquad
\iota=\Delta\colon X\hookrightarrow X\times X.
\]
Thus \(N_Y=N\), \(N_Z=N\oplus N\), \(M_Y=M\), and
\(M_Z=M\oplus M\). On cocharacter lattices the diagonal is
\(N\to N\oplus N\), \(v\mapsto(v,v)\), and its dual on character
lattices is
\[
M\oplus M\longrightarrow M,
\qquad (m_1,m_2)\longmapsto m_1+m_2.
\]
Consequently, the map
\[
M_\R/M\stackrel{\sim}{\longrightarrow}T^\Delta,
\qquad
[m]\longmapsto[(m,-m)]
\]
is an isomorphism. The two rays \((u_\rho,0)\) and
\((0,u_\rho)\) of the product fan give
the same underlying hyperplane
\[
\left\{[m]\in M_\R/M:\
\langle m,u_\rho\rangle\in\Z\right\}.
\]
Thus \(S^\Delta\) is the stratification induced on \(M_\R/M\) by the
\(M\)-periodic arrangement
\[
\langle m,u_\rho\rangle=j,
\qquad \rho\in\Sigma(1),\quad j\in\Z.
\]
The realification of the principal-divisor map in
Equation~\eqref{eqn: fundexactseq} is the linear isomorphism
\[
B_\R\colon M_\R\stackrel{\sim}{\longrightarrow}\R L,
\qquad
m\longmapsto
\bigl(\langle m,u_\rho\rangle\bigr)_{\rho\in\Sigma(1)}.
\]
Since \(B(M)=L\), it descends to an isomorphism
\[
\overline{B}_\R\colon
M_\R/M\stackrel{\sim}{\longrightarrow}\R L/L.
\]
The hyperplane \(\langle m,u_\rho\rangle=j\) is carried to
\(\R L\cap\{p_\rho=j\}\). Therefore \(\overline{B}_\R\) induces an
isomorphism of stratified spaces
\[
\bigl(M_\R/M,S^\Delta\bigr)
\stackrel{\sim}{\longrightarrow}
\faktor{\R L\cap\mathcal{H}_L}{L}.
\]
When \(L\) is unimodular, the right-hand side supports the
Bayer--Popescu--Sturmfels resolution by
\cite[Theorem~3.5]{bayer-popescu-sturmfels}; without unimodularity it
is the same finite quotient stratification whose exit-path data enter
the Hanlon--Hicks--Lazarev construction.

For \(m\in M_\R\), set
\[
D_X(m):=
\sum_{\rho\in\Sigma(1)}
\left\lfloor-\langle m,u_\rho\rangle\right\rfloor D_\rho.
\]
For a lift \(\widetilde\sigma\) of a stratum and
\(m_\sigma\in\widetilde\sigma\), put
\[
D^+_\sigma
:=\sum_{\rho\in\Sigma(1)}
\left\lfloor-\langle m_\sigma,u_\rho\rangle\right\rfloor D_\rho,
\qquad
D^-_\sigma
:=\sum_{\rho\in\Sigma(1)}
\left\lfloor\langle m_\sigma,u_\rho\rangle\right\rfloor D_\rho.
\]
When \(L\) is unimodular, put
\(\mathbf{p}=B_\R(m_\sigma)\), and let \(F\) be the face of
\(\mathcal A_L\) whose relative interior contains \(\mathbf p\).
Its \(L\)-orbit corresponds to \(\sigma\) under
\(\overline{B}_\R\). Since
\(\lceil-\mathbf{p}\rceil=-\lfloor\mathbf{p}\rfloor\), the monomial
\[
\mathbf{x}^{\lceil\mathbf{p}\rceil}
\mathbf{y}^{-\lfloor\mathbf{p}\rfloor}
=
\mathbf{x}^{\lceil\mathbf{p}\rceil}
\mathbf{y}^{\lceil-\mathbf{p}\rceil}
\]
is precisely the Bayer--Popescu--Sturmfels face label, while its
exponent vector
\(\bigl(\lceil\mathbf p\rceil,\lceil-\mathbf p\rceil\bigr)\)
is the ceiling support-function datum defining the
Hanlon--Hicks--Lazarev line-bundle label. Moreover,
\[
\lambda_F=\pi(\lceil\mathbf{p}\rceil),
\qquad
\mu_F=\pi(-\lfloor\mathbf{p}\rfloor),
\]
while
\[
[D^+_\sigma]=-\lambda_F,
\qquad
[D^-_\sigma]=-\mu_F.
\]
The summand attached to \(\sigma\) is the box product
\[
\mathcal{L}^\Delta_\sigma
=\mathcal{O}_X(D^+_\sigma)
\boxtimes\mathcal{O}_X(D^-_\sigma).
\]
Therefore
\[
C_k^\Delta
=\bigoplus_{\substack{\sigma\in S^\Delta\\
\dim\sigma=k}}
\mathcal{O}_X(D^+_\sigma)
\boxtimes\mathcal{O}_X(D^-_\sigma),
\qquad 0\leq k\leq d.
\]
If \(x_\rho\) and \(y_\rho\) denote the Cox variables on the first
and second factors, respectively, set
\[
a^{+}_{\rho,\widetilde\sigma}
=\left\lceil\langle m_\sigma,u_\rho\rangle\right\rceil,
\qquad
a^{-}_{\rho,\widetilde\sigma}
=\left\lceil-\langle m_\sigma,u_\rho\rangle\right\rceil.
\]
For an exit path, choose incident restricted lifts
\(\widetilde\tau<\widetilde\sigma\) whose containing ambient strata
are ordered as above. Its boundary morphism is
multiplication by
\[
\prod_{\rho\in\Sigma(1)}
 x_\rho^{a^+_{\rho,\widetilde\sigma}
-a^+_{\rho,\widetilde\tau}}
 y_\rho^{a^-_{\rho,\widetilde\sigma}
-a^-_{\rho,\widetilde\tau}},
\]
and the differential is the signed sum of these morphisms over all
exit paths. We obtain an exact length-\(d\) resolution
\[
0\longrightarrow C_d^\Delta\longrightarrow\cdots
\longrightarrow C_0^\Delta\longrightarrow
\mathcal{O}_\Delta\longrightarrow0.
\]
This is precisely the diagonal case of
\cite[Corollary~B and Remark~1.3]{hanlon2023resolutions}. In
particular, every term is a direct sum of box products of line bundles
in \(\Theta_X\), and after repetitions are removed, the line bundles
appearing on either factor are exactly the Thomsen collection computed
from the primitive ray generators by the floor formula above. Indeed,
the first-factor labels are
\[
\left\{
\mathcal{O}_X(D_X(m)):\ [m]\in M_\R/M
\right\}
=\Theta_X,
\]
because the strata partition \(M_\R/M\); the second-factor labels give
the same set because \([m]\mapsto[-m]\) is a bijection. In the
notation used below, these distinct fixed-factor line bundles form
\(\mathcal{L}_X\); an ordering
\(\mathcal{E}=(E_0,\ldots,E_{N-1})\) is imposed only when the
exceptional-collection conditions are considered.

\section{Listing all complete unimodular toric surfaces}

The condition that the matrix $B$ with row vectors given by primitive ray generators $\{\rho_i\} = \Sigma(1)$ has linearly independent columns and that all maximal minors lie in $\{0,\pm 1\}$ is rather stringent.

\begin{lemma} The list of complete unimodular toric surfaces is \[ \PP^1 \times \PP^1, \PP^2, Bl_p\PP^2, Bl_{p,q}\PP^2, Bl_{p,q,r}\PP^2 \] either $\PP^1 \times \PP^1$ or $\PP^2$ blown up at any of 0-3 torus invariant points.  
\end{lemma}
\begin{proof}
Let $\rho_1,\ldots,\rho_m$ be the primitive ray generators of $\Sigma$. Since
$B$ has rank $2$, it contains two linearly independent rows, say $\rho_i$
and $\rho_j$. Their determinant is a nonzero maximal minor of $B$, and hence
\[
\det(\rho_i,\rho_j)=\pm 1.
\]
Thus $\rho_i$ and $\rho_j$ form a basis of $N$. After applying an element of
$GL(2,\Z)$, we may assume that $\rho_i=e_1$ and $\rho_j=e_2$.

If $\rho=ae_1+be_2$ is any other primitive ray generator, then
\[
a=\det(\rho,e_2),
\qquad
b=\det(e_1,\rho).
\]
Both determinants are maximal minors of $B$, so
\[
a,b\in\{0,\pm1\}.
\]
Consequently, every primitive ray generator belongs to
\[
\{\pm e_1,\pm e_2,\pm(e_1+e_2),\pm(e_1-e_2)\}.
\]
Moreover, every pair of consecutive ray generators has determinant $\pm1$,
so $\Sigma$ is smooth.

A toric blowdown deletes a ray and therefore preserves the condition that all
maximal minors lie in $\{0,\pm1\}$. By the classification of smooth complete
toric surfaces \cite[Chapter~10]{C-L-S}, repeated toric blowdowns reduce
$X_\Sigma$ to either $\PP^2$ or a Hirzebruch surface $\mathbb{F}_a$. The fan
of $\mathbb{F}_a$ may be taken to have primitive ray generators
\[
e_1,\quad e_2,\quad -e_1+ae_2,\quad -e_2.
\]
Since
\[
\left|\det(e_1,-e_1+ae_2)\right|=a,
\]
unimodularity forces $a\leq1$. Here
$\mathbb{F}_0\cong\PP^1\times\PP^1$, while
$\mathbb{F}_1\cong Bl_p\PP^2$.

For $\PP^2$, take the primitive ray generators
\[
v_1=e_1,\qquad v_2=-e_1+e_2,\qquad v_3=-e_2.
\]
A toric blowup at the fixed point corresponding to
$\operatorname{Cone}(v_i,v_{i+1})$ inserts the ray generated by
$v_i+v_{i+1}$. The three possible new rays are therefore
\[
e_1-e_2,\qquad e_2,\qquad -e_1.
\]
They correspond to blowing up the three torus-fixed points of $\PP^2$.
A further blowup adjacent to one of these newly inserted rays would introduce
a vector of the form
\[
2v_i+v_{i+1}
\qquad\text{or}\qquad
v_i+2v_{i+1},
\]
whose determinant with $v_{i+1}$ or $v_i$, respectively, has absolute value
$2$. Such a blowup therefore violates unimodularity. Hence this branch gives
precisely
\[
\PP^2,\qquad Bl_p\PP^2,\qquad
Bl_{p,q}\PP^2,\qquad Bl_{p,q,r}\PP^2.
\]

Finally, the fan of
$\mathbb{F}_0\cong\PP^1\times\PP^1$ has ray generators
\[
e_1,\qquad e_2,\qquad -e_1,\qquad -e_2.
\]
Blowing up one of its four torus-fixed points inserts one of
\[
e_1+e_2,\qquad -e_1+e_2,\qquad
-e_1-e_2,\qquad e_1-e_2.
\]
Any two adjacent vectors in this list have determinant of absolute value $2$.
Thus at most two can be inserted, and if two are inserted they must be
opposite. These blowups produce no additional isomorphism classes: one
insertion gives $Bl_{p,q}\PP^2$, and two opposite insertions give
$Bl_{p,q,r}\PP^2$. This proves the result.
\end{proof}

Figure~\ref{fig: Fan3} illustrates the $\PP^2$ branch of the preceding proof.
Solid black arrows denote rays already present in the fan, while the dashed
blue arrow in each successive fan denotes the ray inserted by the next toric
blowup. The labels beneath the fans identify the resulting isomorphism
classes of complete unimodular toric surfaces.

\begin{figure}[!htbp]
\centering
\begin{tikzpicture}[
  x=1cm,
  y=1cm,
  fan ray/.style={->,>=stealth,thick},
  new ray/.style={fan ray,blue,densely dashed},
  ray label/.style={font=\scriptsize},
  variety label/.style={font=\small},
  transition/.style={font=\Large}
]
  % The fan of P^2.
  \begin{scope}[shift={(0,0)}]
    \fill (0,0) circle (1.2pt);
    \draw[fan ray] (0,0) -- (1.05,0)
      node[ray label,right] {$e_1$};
    \draw[fan ray] (0,0) -- (-1.05,1.05)
      node[ray label,above left] {$-e_1+e_2$};
    \draw[fan ray] (0,0) -- (0,-1.05)
      node[ray label,below] {$-e_2$};
    \node[variety label] at (0,-1.75) {$\PP^2$};
  \end{scope}

  % Passage to the fan after the first toric blowup.
  \node[font=\scriptsize] at (1.90,0.48) {toric blowup};
  \node[transition] at (1.90,0) {$\Longrightarrow$};

  % Blow up the cone generated by -e_2 and e_1.
  \begin{scope}[shift={(3.8,0)}]
    \fill (0,0) circle (1.2pt);
    \draw[fan ray] (0,0) -- (1.05,0)
      node[ray label,right] {$e_1$};
    \draw[fan ray] (0,0) -- (-1.05,1.05)
      node[ray label,above left] {$-e_1+e_2$};
    \draw[fan ray] (0,0) -- (0,-1.05)
      node[ray label,below left] {$-e_2$};
    \draw[new ray] (0,0) -- (1.05,-1.05)
      node[ray label,below right] {$e_1-e_2$};
    \node[variety label] at (0,-1.75) {$Bl_p\PP^2$};
  \end{scope}

  % Passage to the fan after the second toric blowup.
  \node[font=\scriptsize] at (5.70,0.48) {toric blowup};
  \node[transition] at (5.70,0) {$\Longrightarrow$};

  % Blow up the cone generated by e_1 and -e_1+e_2.
  \begin{scope}[shift={(7.6,0)}]
    \fill (0,0) circle (1.2pt);
    \draw[fan ray] (0,0) -- (1.05,0)
      node[ray label,right] {$e_1$};
    \draw[fan ray] (0,0) -- (-1.05,1.05)
      node[ray label,above left] {$-e_1+e_2$};
    \draw[fan ray] (0,0) -- (0,-1.05)
      node[ray label,below left] {$-e_2$};
    \draw[fan ray] (0,0) -- (1.05,-1.05)
      node[ray label,below right] {$e_1-e_2$};
    \draw[new ray] (0,0) -- (0,1.05)
      node[ray label,above right] {$e_2$};
    \node[variety label] at (0,-1.75) {$Bl_{p,q}\PP^2$};
  \end{scope}

  % Passage to the fan after the third toric blowup.
  \node[font=\scriptsize] at (9.50,0.48) {toric blowup};
  \node[transition] at (9.50,0) {$\Longrightarrow$};

  % Blow up the cone generated by -e_1+e_2 and -e_2.
  \begin{scope}[shift={(11.4,0)}]
    \fill (0,0) circle (1.2pt);
    \draw[fan ray] (0,0) -- (1.05,0)
      node[ray label,right] {$e_1$};
    \draw[fan ray] (0,0) -- (-1.05,1.05)
      node[ray label,above left] {$-e_1+e_2$};
    \draw[fan ray] (0,0) -- (0,-1.05)
      node[ray label,below left] {$-e_2$};
    \draw[fan ray] (0,0) -- (1.05,-1.05)
      node[ray label,below right] {$e_1-e_2$};
    \draw[fan ray] (0,0) -- (0,1.05)
      node[ray label,above right] {$e_2$};
    \draw[new ray] (0,0) --
      node[ray label,midway,below] {$-e_1$} (-1.05,0);
    \node[variety label] at (0,-1.75) {$Bl_{p,q,r}\PP^2$};
  \end{scope}
\end{tikzpicture}
\caption{The $\PP^2$ branch of the classification. From left to right,
the dashed rays inserted by toric blowups are $e_1-e_2$, $e_2$, and
$-e_1$.}
\label{fig: Fan3}
\end{figure}
Figure~\ref{fig: Fan32} illustrates the
$\mathbb{F}_0\cong\PP^1\times\PP^1$ branch of the preceding proof.
Solid black arrows denote rays already present in the fan, while the
dashed blue arrow in each successive fan denotes the ray inserted by
the next toric blowup.

\begin{figure}[!htbp]
\centering
\begin{tikzpicture}[
  x=1cm,
  y=1cm,
  fan ray/.style={->,>=stealth,thick},
  new ray/.style={fan ray,blue,densely dashed},
  ray label/.style={font=\scriptsize},
  variety label/.style={font=\small},
  transition/.style={font=\Large}
]
  % The fan of F_0 = P^1 x P^1.
  \begin{scope}[shift={(0,0)}]
    \fill (0,0) circle (1.2pt);
    \draw[fan ray] (0,0) -- (1.05,0)
      node[ray label,right] {$e_1$};
    \draw[fan ray] (0,0) -- (0,1.05)
      node[ray label,above] {$e_2$};
    \draw[fan ray] (0,0) -- (-1.05,0)
      node[ray label,left] {$-e_1$};
    \draw[fan ray] (0,0) -- (0,-1.05)
      node[ray label,below] {$-e_2$};
    \node[variety label] at (0,-1.75)
      {$\PP^1\times\PP^1$};
  \end{scope}

  % Passage to the fan after the first toric blowup.
  \node[font=\scriptsize] at (2.05,0.48) {toric blowup};
  \node[transition] at (2.05,0) {$\Longrightarrow$};

  % Blow up the cone generated by e_1 and e_2.
  \begin{scope}[shift={(4.1,0)}]
    \fill (0,0) circle (1.2pt);
    \draw[fan ray] (0,0) -- (1.05,0)
      node[ray label,right] {$e_1$};
    \draw[fan ray] (0,0) -- (0,1.05)
      node[ray label,above left] {$e_2$};
    \draw[fan ray] (0,0) -- (-1.05,0)
      node[ray label,left] {$-e_1$};
    \draw[fan ray] (0,0) -- (0,-1.05)
      node[ray label,below right] {$-e_2$};
    \draw[new ray] (0,0) -- (1.05,1.05)
      node[ray label,above right] {$e_1+e_2$};
    \node[variety label] at (0,-1.75)
      {$Bl_{p,q}\PP^2$};
  \end{scope}

  % Passage to the fan after the second toric blowup.
  \node[font=\scriptsize] at (6.15,0.48) {toric blowup};
  \node[transition] at (6.15,0) {$\Longrightarrow$};

  % The only allowed second insertion is the opposite ray.
  \begin{scope}[shift={(8.2,0)}]
    \fill (0,0) circle (1.2pt);
    \draw[fan ray] (0,0) -- (1.05,0)
      node[ray label,right] {$e_1$};
    \draw[fan ray] (0,0) -- (0,1.05)
      node[ray label,above left] {$e_2$};
    \draw[fan ray] (0,0) -- (-1.05,0)
      node[ray label,left] {$-e_1$};
    \draw[fan ray] (0,0) -- (0,-1.05)
      node[ray label,below right] {$-e_2$};
    \draw[fan ray] (0,0) -- (1.05,1.05)
      node[ray label,above right] {$e_1+e_2$};
    \draw[new ray] (0,0) -- (-1.05,-1.05)
      node[ray label,below left] {$-e_1-e_2$};
    \node[variety label] at (0,-1.75)
      {$Bl_{p,q,r}\PP^2$};
  \end{scope}
\end{tikzpicture}
\caption{The $\mathbb{F}_0\cong\PP^1\times\PP^1$ branch of the
classification. Up to symmetry, the first blowup inserts $e_1+e_2$.
The only unimodular second insertion is the opposite ray
$-e_1-e_2$, since either adjacent choice has determinant of absolute
value $2$ with $e_1+e_2$. The fans after one and two blowups represent
$Bl_{p,q}\PP^2$ and $Bl_{p,q,r}\PP^2$, respectively.}
\label{fig: Fan32}
\end{figure}

\section{Full strong exceptional collections for toric unimodular surfaces}
\label{sec: 4}
We now continue on a case-by-case basis to show that the Bayer-Popescu-Sturmfels resolution recalled in Section~\ref{subsec:BPS-construction} yields a strong, full exceptional collection of line bundles for each of the unimodular complete toric surfaces $\PP^1\times \PP^1, \PP^2, Bl_p\PP^2, Bl_{p,q}\PP^2$, and $Bl_{p,q,r}\PP^2.$ In particular, in each case the construction of $\faktor{\R L \cap \mathcal{H}_L}{L}$ is as in Section~\ref{subsec:BPS-construction}. Note that fullness of the resulting collections studied for each space follows from the fact that each collection appears on one side (say, the left-hand factor) of a locally-free resolution of the diagonal by box products of line bundles. Strongness holds for the collection $\mathcal{E}$ on each surface $X$ considered in this section, as $\mathrm{Ext}^k(E_i, E_j)$ is computed to be isomorphic to $0$ for all $k\neq 0$ and all $i,j$ in $\mathcal{E}$.

\subsection{$\PP^1\times \PP^1$}
The primitive ray generators $\Sigma(1)$ for $X_\Sigma=\PP^1 \times \PP^1$ are $\{\rho_1 = e_1, \rho_2 = e_2, \rho_3 = -e_1, \rho_4 = -e_2\}$. This gives $B = \left[ \begin{matrix} 1 & 0 \\ 0 & 1 \\ -1 & 0 \\ 0 & -1 \end{matrix}\right]$. From $D_1 \sim D_3, D_2 \sim D_4$ in $\mathrm{Cl}(\PP^1\times \PP^1)$, we choose the presentation $\mathrm{Cl}(\PP^1\times \PP^1) \cong \left<D_1, D_2\right>$ and have the finite cellular complex $\faktor{\R L \cap \mathcal{H}_L}{L}$ as in Section~\ref{subsec:BPS-construction} given in Figure \hyperref[fig: P1P1]{3}  with monomials and orientation listed. 

\begin{figure}[!htbp]
\centering
\begin{tikzpicture}[
  x=1cm,
  y=1cm,
  cell edge/.style={line width=0.8pt},
  orientation/.style={->,>=stealth,line width=0.8pt},
  double orientation/.style={->>,>=stealth,line width=0.8pt},
  vertex/.style={circle,fill=black,inner sep=2.2pt},
  monomial/.style={font=\small}
]
  \coordinate (v00) at (0,0);
  \coordinate (v10) at (5.8,0);
  \coordinate (v01) at (0,5.8);
  \coordinate (v11) at (5.8,5.8);

  % The four one-dimensional cells.
  \draw[cell edge] (v00) -- (v10);
  \draw[cell edge] (v01) -- (v11);
  \draw[cell edge] (v00) -- (v01);
  \draw[cell edge] (v10) -- (v11);

  % Horizontal pair: one right-pointing arrow on each edge.
  \draw[orientation] (2.45,0) -- (3.35,0);
  \draw[orientation] (2.45,5.8) -- (3.35,5.8);

  % Vertical pair: two upward-pointing arrows on each edge.
  \draw[double orientation] (0,2.45) -- (0,3.35);
  \draw[double orientation] (5.8,2.45) -- (5.8,3.35);

  % Vertex labels.
  \node[monomial,anchor=north east] at (-0.12,-0.15) {$1$};
  \node[monomial,anchor=north west] at (5.95,-0.20)
    {$\dfrac{x_1y_3}{x_3y_1}$};
  \node[monomial,anchor=south east] at (-0.15,5.98)
    {$\dfrac{x_2y_4}{x_4y_2}$};
  \node[monomial,anchor=south west] at (5.95,5.98)
    {$\dfrac{x_1x_2y_3y_4}{x_3x_4y_1y_2}$};

  % Edge labels.
  \node[monomial,anchor=north] at (2.9,-0.35)
    {$x_1y_3$};
  \node[monomial,anchor=south] at (2.9,6.12)
    {$\dfrac{x_1x_2y_3y_4}{x_4y_2}$};
  \node[monomial,anchor=east] at (-0.27,2.9)
    {$x_2y_4$};
  \node[monomial,anchor=west] at (6.07,2.9)
    {$\dfrac{x_1x_2y_3y_4}{x_3y_1}$};

  % Two-dimensional cell label.
  \node[monomial] at (2.9,2.9)
    {$x_1x_2y_3y_4$};

  % Draw the vertices last so that they lie above the edges.
  \node[vertex] at (v00) {};
  \node[vertex] at (v10) {};
  \node[vertex] at (v01) {};
  \node[vertex] at (v11) {};
\end{tikzpicture}
\caption{Finite cellular complex
$\faktor{\mathbb{R}L\cap\mathcal{H}_L}{L}$ for
$\mathbb{P}^1\times\mathbb{P}^1$. Each cell is labeled by the least
common multiple of the Laurent monomial labels of its vertices, and the
arrows record the chosen orientations.}
\label{fig: P1P1}
\end{figure}

The corresponding locally free resolution of $\mathcal{O}_\Delta$ is given by

\[
\begin{array}{ccccccc}
0 \rightarrow \mathcal{O}((-1,-1)\times(-1,-1)) & \rightarrow & \makecell{ \mathcal{O}((-1,0)\times(-1,0)) \\ \oplus \\  \mathcal{O}((0,-1)\times(0,-1)) } & \rightarrow & \mathcal{O} & \rightarrow & 0 \\
& & & & \downarrow & & \\
& & & & \mathcal{O}_\Delta & \rightarrow & 0. \end{array} 
\]

Reading off the collection of line bundles from either the right or left side here yields \[ \mathcal{E} = \{\mathcal{O}(-1,-1), \mathcal{O}(-1,0), \mathcal{O}(0,-1), \mathcal{O}(0,0) \} \] 

which is exceptional since \[ \{ \mathcal{O}_{\PP^1}(-1), \mathcal{O}_{\PP^1} \} \] forms a strong, full exceptional collection of line bundles on $\PP^1$. For instance, the required vanishing for the pair $(\mathcal{O}_{\PP^1\times \PP^1}, \mathcal{O}_{\PP^1\times \PP^1}(-1,-1))$ follows from  

\begin{align*}
        \Ext^k_{D^b_{Coh}(\PP^1\times \PP^1)}(\mathcal{O}_{\PP^1\times \PP^1}, \mathcal{O}(-1,-1))
    &\cong H^k(\PP^1\times \PP^1, \mathcal{O}(-1,-1))\\
    &\cong \bigoplus_{p+q=k} H^p(\PP^1,\mathcal{O}(-1)) \otimes H^q(\PP^1,\mathcal{O}(-1))\\
    &= 0
\end{align*}

for all $k$. 

\subsection{$\PP^2$}

 Labeling primitive ray generators $\Sigma(1) = \{\rho_1 = e_1, \rho_2 = e_2, \rho_3=-e_1-e_2\}$ gives $B = \left[\begin{matrix} 1 & 0 \\ 0 & 1 \\ -1 & -1 \end{matrix}\right]$ and $\mathrm{Pic}(\PP^2) \simeq \left<D_1\right>$ from $D_1 \sim D_3 \sim D_2 \in \mathrm{Cl}(\PP^2)$ gives the labeled finite quotient cellular complex in Figure \hyperref[fig: P2cellcx]{4}.

\begin{figure}[!htbp]
\centering
\begin{tikzpicture}[
  x=1cm,
  y=1cm,
  cell edge/.style={line width=0.8pt},
  vertex/.style={
    circle,
    fill=black,
    inner sep=0pt,
    minimum size=7pt
  },
  cell label/.style={
    font=\small,
    fill=white,
    inner sep=1.5pt
  },
  edge label/.style={
    font=\small,
    fill=white,
    inner sep=1.5pt
  },
  vertex label/.style={
    font=\small,
    fill=white,
    inner sep=1pt
  }
]
  \coordinate (vSW) at (0,0);
  \coordinate (vSE) at (6,0);
  \coordinate (vNE) at (6,6);
  \coordinate (vNW) at (0,6);

  % The two two-cells and their common diagonal edge.
  \draw[cell edge]
    (vSW) -- (vSE) -- (vNE) -- (vNW) -- cycle;
  \draw[cell edge] (vNW) -- (vSE);

  % One arrowhead marks the identified horizontal edges.
  \foreach \y in {0,6} {
    \draw[cell edge]
      (2.86,\y-0.14) -- (3.00,\y) -- (2.86,\y+0.14);
  }

  % Two arrowheads mark the diagonal edge, oriented northwest.
  \draw[cell edge] (2.84,3.16) -- (3.08,3.16);
  \draw[cell edge] (2.84,3.16) -- (2.84,2.92);
  \draw[cell edge] (3.00,3.00) -- (3.24,3.00);
  \draw[cell edge] (3.00,3.00) -- (3.00,2.76);

  % Three arrowheads mark the identified vertical edges.
  \foreach \x in {0,6} {
    \draw[cell edge]
      (\x-0.14,3.24) -- (\x,3.10) -- (\x+0.14,3.24);
    \draw[cell edge]
      (\x-0.14,3.06) -- (\x,2.92) -- (\x+0.14,3.06);
    \draw[cell edge]
      (\x-0.14,2.88) -- (\x,2.74) -- (\x+0.14,2.88);
  }

  % Vertices.
  \node[vertex] at (vSW) {};
  \node[vertex] at (vSE) {};
  \node[vertex] at (vNE) {};
  \node[vertex] at (vNW) {};

  % Labels of the two-cells.
  \node[cell label] at (1.65,1.65)
    {$x_1x_2y_3$};
  \node[cell label] at (3.85,4.40)
    {$\dfrac{x_1x_2y_3^2}{x_3}$};

  % Labels of the one-cells.
  \node[edge label] at (3.00,-0.48)
    {$x_1y_3$};
  \node[edge label] at (3.00,6.58)
    {$\dfrac{x_1x_2y_3^2}{x_3y_2}$};
  \node[edge label] at (-0.48,3.00)
    {$x_2y_3$};
  \node[edge label] at (6.92,3.00)
    {$\dfrac{x_1x_2y_3^2}{x_3y_1}$};
  \node[edge label] at (2.12,3.05)
    {$\dfrac{x_1x_2y_3}{x_3}$};

  % Labels of the vertices in the chosen fundamental domain.
  \node[vertex label,anchor=east] at (-0.18,6.00)
    {$\dfrac{x_2y_3}{x_3y_2}$};
  \node[vertex label,anchor=west] at (6.28,6.00)
    {$\dfrac{x_1x_2y_3^2}{x_3^2y_1y_2}$};
  \node[vertex label,anchor=north east] at (-0.08,-0.15)
    {$1$};
  \node[vertex label,anchor=north west] at (6.28,-0.15)
    {$\dfrac{x_1y_3}{x_3y_1}$};
\end{tikzpicture}
\caption{$\faktor{\R L\cap\mathcal{H}_L}{L}$ for $\PP^2$.}
\label{fig: P2cellcx}
\end{figure}

The corresponding locally-free resolution of $\mathcal{O}_\Delta$ in $D^b_{Coh}(\PP^2\times \PP^2)$ is given by 

\[
\begin{array}{ccccccc}
0 \rightarrow \makecell{\mathcal{O}(-1,-2) \\ \oplus \\ \mathcal{O}(-2,-1)} & \rightarrow & \makecell{ \mathcal{O}(-1,-1) \\ \oplus \\ \mathcal{O}(-1,-1) \\  \oplus \\ \mathcal{O}(-1,-1) } & \rightarrow&  \mathcal{O}&  \rightarrow & 0 \\
& & & & \downarrow \\
& & & & \mathcal{O}_\Delta & \rightarrow & 0. \end{array} \]

The line bundles appearing on either the left or right-hand side give \[ \mathcal{E} = \{ \mathcal{O}(-2), \mathcal{O}(-1), \mathcal{O} \}\] which is a well-known strong, full exceptional collection of line bundles on $\PP^2$ due to Beilinson \cite{Beilinson1978}. 

\subsection{$Bl_p \PP^2$}
For $Bl_p \PP^2$ having primitive ray generators $\Sigma(1) = \{ \rho_1 = e_1, \rho_2 = e_2, \rho_3 = -e_1-e_2, \rho_4 = e_1+e_2 \}$, we have $D_1+D_4 \sim D_3 \sim D_2 + D_4$ in $\mathrm{Cl}(Bl_p \PP^2)$, so we choose the presentation $\mathrm{Cl}(Bl_p \PP^2) \cong \left< D_1, D_4 \right>$. We have labeled finite quotient cellular complex given in Figure~\ref{fig: BlpP2}. 

\begin{figure}[!htbp]
\centering
\begin{tikzpicture}[
  x=1cm,
  y=1cm,
  cell edge/.style={line width=0.8pt},
  vertex/.style={
    circle,
    fill=black,
    inner sep=0pt,
    minimum size=7pt
  },
  cell label/.style={
    font=\small,
    fill=white,
    inner sep=1.5pt
  },
  edge label/.style={
    font=\small,
    fill=white,
    inner sep=1.5pt
  },
  vertex label/.style={
    font=\small,
    fill=white,
    inner sep=1pt
  }
]
  \coordinate (vSW) at (0,0);
  \coordinate (vSE) at (6,0);
  \coordinate (vNE) at (6,6);
  \coordinate (vNW) at (0,6);

  % The two two-cells and their common diagonal edge.
  \draw[cell edge]
    (vSW) -- (vSE) -- (vNE) -- (vNW) -- cycle;
  \draw[cell edge] (vNW) -- (vSE);

  % One arrowhead marks the identified horizontal edges.
  \foreach \y in {0,6} {
    \draw[cell edge]
      (2.86,\y-0.14) -- (3.00,\y) -- (2.86,\y+0.14);
  }

  % Two arrowheads mark the diagonal edge, oriented northwest.
  \draw[cell edge] (2.84,3.16) -- (3.08,3.16);
  \draw[cell edge] (2.84,3.16) -- (2.84,2.92);
  \draw[cell edge] (3.00,3.00) -- (3.24,3.00);
  \draw[cell edge] (3.00,3.00) -- (3.00,2.76);

  % Three arrowheads mark the identified vertical edges.
  \foreach \x in {0,6} {
    \draw[cell edge]
      (\x-0.14,3.24) -- (\x,3.10) -- (\x+0.14,3.24);
    \draw[cell edge]
      (\x-0.14,3.06) -- (\x,2.92) -- (\x+0.14,3.06);
    \draw[cell edge]
      (\x-0.14,2.88) -- (\x,2.74) -- (\x+0.14,2.88);
  }

  % Vertices.
  \node[vertex] at (vSW) {};
  \node[vertex] at (vSE) {};
  \node[vertex] at (vNE) {};
  \node[vertex] at (vNW) {};

  % Labels of the two-cells.
  \node[cell label] at (1.72,1.68)
    {$x_1x_2x_4y_3$};
  \node[cell label] at (4.03,4.43)
    {$\dfrac{x_1x_2x_4^2y_3^2}{x_3y_4}$};

  % Labels of the one-cells.
  \node[edge label] at (3.00,-0.48)
    {$x_1x_4y_3$};
  \node[edge label] at (3.00,6.58)
    {$\dfrac{x_1x_2x_4^2y_3^2}{x_3y_2y_4}$};
  \node[edge label] at (-0.78,3.00)
    {$x_2x_4y_3$};
  \node[edge label] at (7.00,3.00)
    {$\dfrac{x_1x_2x_4^2y_3^2}{x_3y_1y_4}$};
  \node[edge label] at (1.78,3.38)
    {$\dfrac{x_1x_2x_4y_3}{x_3y_4}$};

  % Labels of the vertices in the chosen fundamental domain.
  \node[vertex label,anchor=east] at (-0.18,6.00)
    {$\dfrac{x_2x_4y_3}{x_3y_2y_4}$};
  \node[vertex label,anchor=west] at (6.28,6.00)
    {$\dfrac{x_1x_2x_4^2y_3^2}
    {x_3^2y_1y_2y_4^2}$};
  \node[vertex label,anchor=north east] at (-0.08,-0.15)
    {$1$};
  \node[vertex label,anchor=north west] at (6.28,-0.15)
    {$\dfrac{x_1x_4y_3}{x_3y_1y_4}$};
\end{tikzpicture}
\caption{$\faktor{\R L\cap\mathcal{H}_L}{L}$ for $Bl_p\PP^2$.}
\label{fig: BlpP2}
\end{figure}

The corresponding locally-free resolution in $D^b_{Coh}(Bl_p \PP^2 \times Bl_p \PP^2)$ is given by 

\[\begin{array}{ccccccc}
0 \rightarrow \makecell{\mathcal{O}(-1,-1)\times(-2,-1) \\ \oplus \\ \mathcal{O}(-2,-1)\times(-1,-1) } & \rightarrow & \makecell{ \mathcal{O}(-1,-1)\times(-1,-1) \\ \oplus \\ \mathcal{O}(-1,0)\times\mathcal{O}(-1,0) \\ \oplus \\ \mathcal{O}(-1,-1)\times (-1,-1) } & \rightarrow & \mathcal{O} & \rightarrow & 0 \\
& & & & \downarrow \\
& & & & \mathcal{O}_\Delta & \rightarrow & 0.\end{array} \]

Reading off the line bundles from either the left-hand or right-hand side yields \[ \mathcal{E} = \{ \mathcal{O}(-2,-1), \mathcal{O}(-1,-1), \mathcal{O}(-1,0), \mathcal{O} \} \]

which can be shown to be exceptional from computing $H^*(Bl_p\PP^2, \mathcal{O}(a_1, a_2))$ for $(a_1,a_2)$ appearing in Figure~\ref{fig:blphoms}.
That is, the fact that $\mathcal{E}$ is exceptional is given by $H^*(Bl_p\PP^2, \mathcal{O}(a_1,a_2)) = 0$ for $(a_1, a_2)$ in $\{(-2,-1), (-1,-1), (-1,0), (0,-1)\}$. These are exactly the arrows appearing in Figure~\ref{fig:blphoms}.

\begin{figure}[!htbp]
\centering
\begin{tikzpicture}[
  bundle/.style={
    font=\small,
    fill=white,
    inner xsep=5pt,
    inner ysep=2pt
  },
  vanishing arrow/.style={
    ->,
    >=latex,
    line width=0.65pt,
    shorten <=2pt,
    shorten >=2pt
  },
  arrow label/.style={
    font=\scriptsize,
    fill=white,
    inner sep=1.5pt
  }
]
  % The collection is ordered from E_0 at the bottom to E_3 at the top.
  \node[bundle] (E4) at (0,0)
    {$E_3=\mathcal{O}$};
  \node[bundle] (E3) at (0,-1.75)
    {$E_2=\mathcal{O}(-1,0)$};
  \node[bundle] (E2) at (0,-3.50)
    {$E_1=\mathcal{O}(-1,-1)$};
  \node[bundle] (E1) at (0,-5.25)
    {$E_0=\mathcal{O}(-2,-1)$};

  % Adjacent pairs.
  \draw[vanishing arrow] (E4) -- (E3);
  \draw[vanishing arrow] (E3) -- (E2);
  \draw[vanishing arrow] (E2) -- (E1);

  % Pairs separated by one intermediate object.
  \draw[vanishing arrow]
    (E4.west) to[bend right=43] (E2.west);
  \draw[vanishing arrow]
    (E3.east) to[bend left=43] (E1.east);

  % The outermost pair.
  \draw[vanishing arrow]
    (E4.east) to[bend left=68] (E1.east);

  % Draw labels after the arrows so that no curve crosses the text.
  \path (E4) -- node[arrow label,right=4pt]
    {$(-1,0)$} (E3);
  \path (E3) -- node[arrow label,right=4pt]
    {$(0,-1)$} (E2);
  \path (E2) -- node[arrow label,right=4pt]
    {$(-1,0)$} (E1);

  \path (E4.west) to[bend right=43]
    node[arrow label,left=3pt] {$(-1,-1)$}
    (E2.west);
  \path (E3.east) to[bend left=43]
    node[arrow label,right=3pt] {$(-1,-1)$}
    (E1.east);

  \path (E4.east) to[bend left=68]
    node[arrow label,right=4pt] {$(-2,-1)$}
    (E1.east);
\end{tikzpicture}
\caption{Backward Ext-groups for $\mathcal{E}$ on
$X=Bl_p\PP^2$. An arrow $E_i\to E_j$, with $i>j$, labeled
$a=(a_1,a_2)$ records
$\operatorname{Ext}^{\bullet}(E_i,E_j)\cong
H^{\bullet}(X,\mathcal{O}_X(a))=0$, proving that $\mathcal{E}$ is
exceptional.}
\label{fig:blphoms}
\end{figure}

$\mathcal{E}$ is strong, since $H^*(Bl_p \PP^2, \mathcal{O}(b_1,b_2))$ is concentrated in degree $0$ for $(b_1, b_2)$ appearing in Figure~\ref{fig:BlpStrongHoms}.

\begin{figure}[!htbp]
\centering
\begin{tikzpicture}[
  bundle/.style={
    font=\small,
    fill=white,
    inner xsep=5pt,
    inner ysep=2pt
  },
  cohomology arrow/.style={
    ->,
    >=Latex,
    line width=0.65pt,
    shorten <=2pt,
    shorten >=2pt
  },
  arrow label/.style={
    font=\scriptsize,
    fill=white,
    inner sep=1.5pt
  }
]
  % The collection is ordered from E_0 at the bottom to E_3 at the top.
  \node[bundle] (E4) at (0,0)
    {$E_3=\mathcal{O}$};
  \node[bundle] (E3) at (0,-1.75)
    {$E_2=\mathcal{O}(-1,0)$};
  \node[bundle] (E2) at (0,-3.50)
    {$E_1=\mathcal{O}(-1,-1)$};
  \node[bundle] (E1) at (0,-5.25)
    {$E_0=\mathcal{O}(-2,-1)$};

  % Adjacent pairs.
  \draw[cohomology arrow] (E3) -- (E4);
  \draw[cohomology arrow] (E2) -- (E3);
  \draw[cohomology arrow] (E1) -- (E2);

  % Pairs separated by one intermediate object.
  \draw[cohomology arrow]
    (E2.west) to[bend left=43] (E4.west);
  \draw[cohomology arrow]
    (E1.east) to[bend right=43] (E3.east);

  % The outermost pair.
  \draw[cohomology arrow]
    (E1.east) to[bend right=68] (E4.east);

  % Labels are placed after the arrows so that no curve crosses the text.
  \path (E3) -- node[arrow label,right=4pt]
    {$(1,0)$} (E4);
  \path (E2) -- node[arrow label,right=4pt]
    {$(0,1)$} (E3);
  \path (E1) -- node[arrow label,right=4pt]
    {$(1,0)$} (E2);

  \path (E2.west) to[bend left=43]
    node[arrow label,left=3pt] {$(1,1)$}
    (E4.west);
  \path (E1.east) to[bend right=43]
    node[arrow label,right=3pt] {$(1,1)$}
    (E3.east);

  \path (E1.east) to[bend right=68]
    node[arrow label,right=4pt] {$(2,1)$}
    (E4.east);
\end{tikzpicture}
\caption{Forward Ext-groups for $\mathcal{E}$ on
$X=Bl_p\PP^2$. An arrow $E_j\to E_i$, with $i>j$, labeled
$b=(b_1,b_2)$ records
$\operatorname{Ext}^{k}(E_j,E_i)\cong
H^{k}(X,\mathcal{O}_X(b))$ for $k\geq 0$. For every such pair,
$H^{k}(X,\mathcal{O}_X(b))=0$ for $k>0$; together with
exceptionality, this proves that $\mathcal{E}$ is strong.}
\label{fig:BlpStrongHoms}
\end{figure}

\subsection{$Bl_{p,q}\PP^2$} Let the blow-up of $\PP^2$ at 2 torus invariant points have primitive ray generators $\Sigma(1) = \{ \rho_1 = e_1, \rho_2 = e_2, \rho_3 = -e_1 - e_2, \rho_4 = e_1 + e_2, \rho_5 = -e_2 \}$. This gives $B = \left[\begin{matrix} 1 & 0 \\ 0 & 1 \\ -1 & -1 \\ 1 & 1 \\ 0 & -1 \end{matrix}\right]$. We have relations $D_1 + D_4 \sim D_3$ and $D_2 + D_4 \sim D_3 + D_5$ which implies $D_2 \sim D_1 + D_5$ in $\mathrm{Cl}(Bl_{p,q}\PP^2)$, so we choose presentation $\mathrm{Cl}(Bl_{p,q}\PP^2) \cong \left< D_1, D_4, D_5\right>$. This gives finite cellular complex $\faktor{\R L \cap \mathcal{H}_L}{L}$ in Figure~\ref{fig: BlpqP2}.

\begin{figure}[!htbp]
\centering
\begin{tikzpicture}[
  x=1cm,
  y=1cm,
  cell edge/.style={line width=0.8pt},
  vertex/.style={
    circle,
    fill=black,
    inner sep=0pt,
    minimum size=7pt
  },
  cell label/.style={
    font=\small,
    fill=white,
    inner sep=1.5pt
  },
  edge label/.style={
    font=\small,
    fill=white,
    inner sep=1.5pt
  },
  vertex label/.style={
    font=\small,
    fill=white,
    inner sep=1pt
  }
]
  \coordinate (vSW) at (0,0);
  \coordinate (vSE) at (6,0);
  \coordinate (vNE) at (6,6);
  \coordinate (vNW) at (0,6);

  % The two two-cells and their common diagonal edge.
  \draw[cell edge]
    (vSW) -- (vSE) -- (vNE) -- (vNW) -- cycle;
  \draw[cell edge] (vNW) -- (vSE);

  % One arrowhead marks the identified horizontal edges.
  \foreach \y in {0,6} {
    \draw[cell edge]
      (2.86,\y-0.14) -- (3.00,\y) -- (2.86,\y+0.14);
  }

  % Two arrowheads mark the diagonal edge, oriented northwest.
  \draw[cell edge] (2.84,3.16) -- (3.08,3.16);
  \draw[cell edge] (2.84,3.16) -- (2.84,2.92);
  \draw[cell edge] (3.00,3.00) -- (3.24,3.00);
  \draw[cell edge] (3.00,3.00) -- (3.00,2.76);

  % Three arrowheads mark the identified vertical edges.
  \foreach \x in {0,6} {
    \draw[cell edge]
      (\x-0.14,3.24) -- (\x,3.10) -- (\x+0.14,3.24);
    \draw[cell edge]
      (\x-0.14,3.06) -- (\x,2.92) -- (\x+0.14,3.06);
    \draw[cell edge]
      (\x-0.14,2.88) -- (\x,2.74) -- (\x+0.14,2.88);
  }

  % Vertices.
  \node[vertex] at (vSW) {};
  \node[vertex] at (vSE) {};
  \node[vertex] at (vNE) {};
  \node[vertex] at (vNW) {};

  % Labels of the two-cells.
  \node[cell label] at (1.72,1.68)
    {$x_1x_2x_4y_3y_5$};
  \node[cell label] at (4.08,4.48)
    {$\dfrac{x_1x_2x_4^2y_3^2y_5}{x_3y_4}$};

  % Labels of the one-cells.
  \node[edge label] at (3.00,-0.50)
    {$x_1x_4y_3$};
  \node[edge label] at (3.00,6.66)
    {$\dfrac{x_1x_2x_4^2y_3^2y_5}
    {x_3x_5y_2y_4}$};
  \node[edge label,anchor=east] at (-0.38,3.00)
    {$x_2x_4y_3y_5$};
  \node[edge label,anchor=west] at (6.38,3.00)
    {$\dfrac{x_1x_2x_4^2y_3^2y_5}
    {x_3x_5y_1y_4}$};
  \node[edge label,anchor=south east] at (2.55,3.58)
    {$\dfrac{x_1x_2x_4y_3y_5}{x_3y_4}$};

  % Labels of the vertices in the chosen fundamental domain.
  \node[vertex label,anchor=east] at (-0.28,6.00)
    {$\dfrac{x_2x_4y_3y_5}{x_3x_5y_2y_4}$};
  \node[vertex label,anchor=west] at (6.28,6.00)
    {$\dfrac{x_1x_2x_4^2y_3^2y_5}
    {x_3^2x_5y_1y_2y_4^2}$};
  \node[vertex label,anchor=north east] at (-0.08,-0.15)
    {$1$};
  \node[vertex label,anchor=north west] at (6.28,-0.15)
    {$\dfrac{x_1x_4y_3}{x_3y_1y_4}$};
\end{tikzpicture}
\caption{$\faktor{\R L\cap\mathcal{H}_L}{L}$ for
$Bl_{p,q}\PP^2$.}
\label{fig: BlpqP2}
\end{figure}

The corresponding locally-free resolution in $D^b_{Coh}(Bl_{p,q}\PP^2 \times Bl_{p,q}\PP^2) $ is

\[\begin{array}{ccccccc}
0 \rightarrow\makecell{\mathcal{O}(-1,-1, -1)\times(-2, -1, -1) \\ \oplus \\ \mathcal{O}(-2,-1, -1)\times(-1,-1, -1) } & \rightarrow & \makecell{ \mathcal{O}(-1,-1, 0)\times(-1,-1, 0) \\ \oplus \\ \mathcal{O}(-1,0, -1)\times\mathcal{O}(-1,0, 
 -1) \\ \oplus \\ \mathcal{O}(-1,-1, -1)\times (-1,-1, -1) } & \rightarrow & \mathcal{O} & \rightarrow & 0 \\
& & & & \downarrow \\
& & & & \mathcal{O}_\Delta & \rightarrow & 0.\end{array} \]

Reading off the line bundles appearing on the left-hand side yields

\[ \mathcal{E} = \{ \mathcal{O}(-2,-1,-1), \mathcal{O}(-1,-1,-1), \mathcal{O}(-1,-1,0), \mathcal{O}(-1,0,-1), \mathcal{O} \}. \] 

Here, when we blow-up two torus-invariant points, both spheres are orthogonal, so this is one of two possible orders for $\mathcal{E}$. To show that $\mathcal{E}$ is exceptional, we check vanishing of $H^*(Bl_{p,q}\mathcal{O}(a_1,a_2,a_3))$ in all cohomological degrees, for $(a_1, a_2,a_3)$ appearing in Figure~\ref{fig:BlpqExc}. We have $Hom_{D^b_{Coh}(Bl_{p,q}\PP^2)}(E_i, E_i) = H^*(Bl_{p,q}\PP^2, \mathcal{O}) = k$ in degree $0$ for all $E_i$ in $\mathcal{E}$. $\mathcal{E}$ is also strong, since $H^*(Bl_{p,q}\PP^2, \mathcal{O}(b_1,b_2,b_3))$ is concentrated in degree $0$ for $(b_1, b_2,b_3)$ appearing in Figure~\ref{fig:BlpqStr}.

\begin{figure}[!htbp]
\centering
\begin{tikzpicture}[
  bundle/.style={
    font=\small,
    fill=white,
    inner xsep=5pt,
    inner ysep=2pt
  },
  vanishing arrow/.style={
    ->,
    >=latex,
    line width=0.65pt,
    shorten <=2pt,
    shorten >=2pt
  },
  crossing arrow/.style={
    vanishing arrow,
    preaction={draw=white,line width=2.6pt,-}
  },
  arrow label/.style={
    font=\scriptsize,
    fill=white,
    inner sep=1.5pt
  }
]
  % The collection is ordered from E_0 at the bottom to E_4 at the top.
  \node[bundle] (E5) at (0,0)
    {$E_4=\mathcal{O}$};
  \node[bundle] (E4) at (0,-1.75)
    {$E_3=\mathcal{O}(-1,0,-1)$};
  \node[bundle] (E3) at (0,-3.50)
    {$E_2=\mathcal{O}(-1,-1,0)$};
  \node[bundle] (E2) at (0,-5.25)
    {$E_1=\mathcal{O}(-1,-1,-1)$};
  \node[bundle] (E1) at (0,-7.00)
    {$E_0=\mathcal{O}(-2,-1,-1)$};

  % Adjacent pairs.
  \draw[vanishing arrow] (E5) -- (E4);
  \draw[vanishing arrow] (E4) -- (E3);
  \draw[vanishing arrow] (E3) -- (E2);
  \draw[vanishing arrow] (E2) -- (E1);

  % Pairs separated by one intermediate object.
  \draw[vanishing arrow]
    (E5.west) to[bend right=42] (E3.west);
  \draw[vanishing arrow]
    (E4.east) to[bend left=42] (E2.east);
  \draw[vanishing arrow]
    (E3.west) to[bend right=42] (E1.west);

  % Pairs separated by two intermediate objects.
  \draw[vanishing arrow]
    (E5.east) to[bend left=53] (E2.east);
  \draw[crossing arrow]
    (E4.west) to[bend right=53] (E1.west);

  % The outermost pair.
  \draw[vanishing arrow]
    (E5.east) to[bend left=68] (E1.east);

  % Place labels after drawing all arrows.
  \path (E5) -- node[arrow label,right=4pt]
    {$(-1,0,-1)$} (E4);
  \path (E4) -- node[arrow label,right=4pt]
    {$(0,-1,1)$} (E3);
  \path (E3) -- node[arrow label,right=4pt]
    {$(0,0,-1)$} (E2);
  \path (E2) -- node[arrow label,right=4pt]
    {$(-1,0,0)$} (E1);

  \path (E5.west) to[bend right=42]
    node[arrow label,left=3pt] {$(-1,-1,0)$}
    (E3.west);
  \path (E4.east) to[bend left=42]
    node[arrow label,right=3pt] {$(0,-1,0)$}
    (E2.east);
  \path (E3.west) to[bend right=42]
    node[arrow label,left=3pt] {$(-1,0,-1)$}
    (E1.west);

  \path (E5.east) to[bend left=53]
    node[arrow label,pos=0.32,right=3pt] {$(-1,-1,-1)$}
    (E2.east);
  \path (E4.west) to[bend right=53]
    node[arrow label,left=3pt] {$(-1,-1,0)$}
    (E1.west);

  \path (E5.east) to[bend left=68]
    node[arrow label,pos=0.64,right=4pt] {$(-2,-1,-1)$}
    (E1.east);
\end{tikzpicture}
\caption{Backward Ext-groups for $\mathcal{E}$ on
$X=Bl_{p,q}\PP^2$. An arrow $E_i\to E_j$, with $i>j$, labeled
$(a_1,a_2,a_3)$ records
$\operatorname{Ext}^{\bullet}(E_i,E_j)\cong
H^{\bullet}(X,\mathcal{O}_X(a_1,a_2,a_3))=0$, proving that $\mathcal{E}$ is
exceptional.}
\label{fig:BlpqExc}
\end{figure}

\begin{figure}[!htbp]
\centering
\begin{tikzpicture}[
  bundle/.style={
    font=\small,
    fill=white,
    inner xsep=5pt,
    inner ysep=2pt
  },
  cohomology arrow/.style={
    ->,
    >=latex,
    line width=0.65pt,
    shorten <=2pt,
    shorten >=2pt
  },
  crossing arrow/.style={
    cohomology arrow,
    preaction={draw=white,line width=2.6pt,-}
  },
  arrow label/.style={
    font=\scriptsize,
    fill=white,
    inner sep=1.5pt
  }
]
  \node[bundle] (E5) at (0,0)
    {$E_4=\mathcal{O}$};
  \node[bundle] (E4) at (0,-1.75)
    {$E_3=\mathcal{O}(-1,0,-1)$};
  \node[bundle] (E3) at (0,-3.50)
    {$E_2=\mathcal{O}(-1,-1,0)$};
  \node[bundle] (E2) at (0,-5.25)
    {$E_1=\mathcal{O}(-1,-1,-1)$};
  \node[bundle] (E1) at (0,-7.00)
    {$E_0=\mathcal{O}(-2,-1,-1)$};

  % Adjacent pairs.
  \draw[cohomology arrow] (E4) -- (E5);
  \draw[cohomology arrow] (E3) -- (E4);
  \draw[cohomology arrow] (E2) -- (E3);
  \draw[cohomology arrow] (E1) -- (E2);

  % Pairs separated by one intermediate object.
  \draw[cohomology arrow]
    (E3.west) to[bend left=42] (E5.west);
  \draw[cohomology arrow]
    (E2.east) to[bend right=42] (E4.east);
  \draw[cohomology arrow]
    (E1.west) to[bend left=42] (E3.west);

  % Pairs separated by two intermediate objects.
  \draw[cohomology arrow]
    (E2.east) to[bend right=53] (E5.east);
  \draw[crossing arrow]
    (E1.west) to[bend left=53] (E4.west);

  % The outermost pair.
  \draw[cohomology arrow]
    (E1.east) to[bend right=68] (E5.east);

  % Place the labels after drawing all arrows.
  \path (E4) -- node[arrow label,right=4pt]
    {$(1,0,1)$} (E5);
  \path (E3) -- node[arrow label,right=4pt]
    {$(0,1,-1)$} (E4);
  \path (E2) -- node[arrow label,right=4pt]
    {$(0,0,1)$} (E3);
  \path (E1) -- node[arrow label,right=4pt]
    {$(1,0,0)$} (E2);

  \path (E3.west) to[bend left=42]
    node[arrow label,left=3pt] {$(1,1,0)$}
    (E5.west);
  \path (E2.east) to[bend right=42]
    node[arrow label,right=3pt] {$(0,1,0)$}
    (E4.east);
  \path (E1.west) to[bend left=42]
    node[arrow label,left=3pt] {$(1,0,1)$}
    (E3.west);

  \path (E2.east) to[bend right=53]
    node[arrow label,pos=0.68,right=3pt] {$(1,1,1)$}
    (E5.east);
  \path (E1.west) to[bend left=53]
    node[arrow label,left=3pt] {$(1,1,0)$}
    (E4.west);

  \path (E1.east) to[bend right=68]
    node[arrow label,pos=0.36,right=4pt] {$(2,1,1)$}
    (E5.east);
\end{tikzpicture}
\caption{Forward Ext-groups for $\mathcal{E}$ on
$X=Bl_{p,q}\PP^2$. An arrow $E_j\to E_i$, with $i>j$, labeled
$(b_1,b_2,b_3)$ records
$\operatorname{Ext}^{k}(E_j,E_i)\cong
H^{k}(X,\mathcal{O}_X(b_1,b_2,b_3))$ for $k\geq 0$. For every such pair,
$H^{k}(X,\mathcal{O}_X(b_1,b_2,b_3))=0$ for $k>0$; together with
exceptionality, this proves that $\mathcal{E}$ is strong.}
\label{fig:BlpqStr}
\end{figure}

\subsection{$Bl_{p,q,r}\PP^2$} Let the blow-up of $\PP^2$ at three torus-invariant points have primitive ray generators $\Sigma(1) = \{ \rho_1 = e_1, \rho_2 = e_2, \rho_3 = -e_1-e_2, \rho_4 = e_1 + e_2, \rho_5 = -e_2, \rho_6 = -e_1 \}$. This gives $B = \left[\begin{matrix} 1 & 0 \\ 0 & 1 \\ -1 & -1 \\ 1 & 1 \\ 0 & -1 \\ -1 & 0  \end{matrix} \right]$. We have $D_2 \sim D_1 + D_5 - D_6$ and $D_3 \sim D_1 + D_4- D_6$ in $\mathrm{Cl}(Bl_{p,q,r}\PP^2)$, so we choose the presentation $\mathrm{Cl}(Bl_{p,q,r}\PP^2) \cong \left< D_1, D_4, D_5, D_6\right>$.

This gives the finite cellular complex $\faktor{\R L \cap \mathcal{H}_L}{L}$ in Figure~\ref{fig: BlpqrP2}.

\begin{figure}[!htbp]
\centering
\begin{tikzpicture}[
  x=1cm,
  y=1cm,
  cell edge/.style={line width=0.8pt},
  vertex/.style={
    circle,
    fill=black,
    inner sep=0pt,
    minimum size=7pt
  },
  cell label/.style={
    font=\small,
    fill=white,
    inner sep=1.5pt
  },
  edge label/.style={
    font=\small,
    fill=white,
    inner sep=1.5pt
  },
  vertex label/.style={
    font=\small,
    fill=white,
    inner sep=1pt
  }
]
  \coordinate (vSW) at (0,0);
  \coordinate (vSE) at (6,0);
  \coordinate (vNE) at (6,6);
  \coordinate (vNW) at (0,6);

  % The two two-cells and their common diagonal edge.
  \draw[cell edge]
    (vSW) -- (vSE) -- (vNE) -- (vNW) -- cycle;
  \draw[cell edge] (vNW) -- (vSE);

  % One arrowhead marks the identified horizontal edges.
  \foreach \y in {0,6} {
    \draw[cell edge]
      (2.86,\y-0.14) -- (3.00,\y) -- (2.86,\y+0.14);
  }

  % Two arrowheads mark the diagonal edge, oriented northwest.
  \draw[cell edge] (2.84,3.16) -- (3.08,3.16);
  \draw[cell edge] (2.84,3.16) -- (2.84,2.92);
  \draw[cell edge] (3.00,3.00) -- (3.24,3.00);
  \draw[cell edge] (3.00,3.00) -- (3.00,2.76);

  % Three arrowheads mark the identified vertical edges.
  \foreach \x in {0,6} {
    \draw[cell edge]
      (\x-0.14,3.24) -- (\x,3.10) -- (\x+0.14,3.24);
    \draw[cell edge]
      (\x-0.14,3.06) -- (\x,2.92) -- (\x+0.14,3.06);
    \draw[cell edge]
      (\x-0.14,2.88) -- (\x,2.74) -- (\x+0.14,2.88);
  }

  % Vertices.
  \node[vertex] at (vSW) {};
  \node[vertex] at (vSE) {};
  \node[vertex] at (vNE) {};
  \node[vertex] at (vNW) {};

  % Labels of the two-cells.
  \node[cell label] at (1.75,1.68)
    {$x_1x_2x_4y_3y_5y_6$};
  \node[cell label] at (4.10,4.47)
    {$\dfrac{x_1x_2x_4^2y_3^2y_5y_6}{x_3y_4}$};

  % Labels of the one-cells.
  \node[edge label] at (3.00,-0.52)
    {$x_1x_4y_3y_6$};
  \node[edge label] at (3.00,6.68)
    {$\dfrac{x_1x_2x_4^2y_3^2y_5y_6}
    {x_3x_5y_2y_4}$};
  \node[edge label,anchor=east] at (-0.42,3.00)
    {$x_2x_4y_3y_5$};
  \node[edge label,anchor=west] at (6.42,3.00)
    {$\dfrac{x_1x_2x_4^2y_3^2y_5y_6}
    {x_3x_6y_1y_4}$};
  \node[edge label,anchor=south east] at (2.50,3.60)
    {$\dfrac{x_1x_2x_4y_3y_5y_6}{x_3y_4}$};

  % Labels of the vertices in the chosen fundamental domain.
  \node[vertex label,anchor=east] at (-0.28,6.00)
    {$\dfrac{x_2x_4y_3y_5}{x_3x_5y_2y_4}$};
  \node[vertex label,anchor=west] at (6.28,6.00)
    {$\dfrac{x_1x_2x_4^2y_3^2y_5y_6}
    {x_3^2x_5x_6y_1y_2y_4^2}$};
  \node[vertex label,anchor=north east] at (-0.08,-0.15)
    {$1$};
  \node[vertex label,anchor=north west] at (6.28,-0.15)
    {$\dfrac{x_1x_4y_3y_6}{x_3x_6y_1y_4}$};
\end{tikzpicture}
\caption{$\faktor{\R L\cap\mathcal{H}_L}{L}$ for
$Bl_{p,q,r}\PP^2$.}
\label{fig: BlpqrP2}
\end{figure}

The corresponding locally-free resolution in $D^b_{Coh}(Bl_{p,q,r}\PP^2 \times Bl_{p,q,r}\PP^2)$ is

\[\begin{array}{ccccccc}
0 \rightarrow \makecell{\mathcal{O}(-1,-1,-1,0)\times(-2,-1,-1,1) \\ \oplus \\ \mathcal{O}(-2,-1, -1, 1)\times(-1, -1, -1, 0) } & \rightarrow & \makecell{ \mathcal{O}(-1,-1,0,0)\times (-1,-1,0,0) \\ \oplus \\ \mathcal{O}(-1,0,-1,0)\times\mathcal{O}(-1,0,-1,0) \\ \oplus \\ \mathcal{O}(-1,-1,-1,1) \times (-1,-1,-1,1)} & \rightarrow & \mathcal{O} & \rightarrow & 0 \\
& & & & \downarrow \\
& & & & \mathcal{O}_\Delta & \rightarrow & 0.\end{array} \]

Reading off the line bundles on either the left- or right-hand side yields \[ \mathcal{E} = \{ \mathcal{O}(-2,-1,-1,1), \mathcal{O}(-1,-1,-1,0), \mathcal{O}(-1,0,-1,0), \mathcal{O}(-1,-1,0,0), \mathcal{O}(-1,-1,-1,1), \mathcal{O} \}. \]

To show that $\mathcal{E}$ is exceptional, we check vanishing of $H^*(Bl_{p,q,r}\PP^2, \mathcal{O}(a_1,a_2,a_3,a_4))=0$ for arrow label $(a_1,a_2,a_3,a_4)$ in Figure~\ref{fig: BlpqrExc2}. Strongness of $\mathcal{E}$ follows from $H^*(Bl_{p,q,r}\PP^2, \mathcal{O}(b_1,b_2,b_3,b_4))$ being concentrated in cohomological degree $0$ for arrow label $(b_1, b_2, b_3, b_4)$  in Figure \hyperref[fig: BlpqrStrong]{13}. 

\begin{figure}[!htbp]
\centering
\begin{tikzpicture}[
  bundle/.style={
    font=\small,
    fill=white,
    inner xsep=5pt,
    inner ysep=2pt
  },
  vanishing arrow/.style={
    ->,
    >=latex,
    line width=0.65pt,
    shorten <=2pt,
    shorten >=2pt
  },
  arrow label/.style={
    font=\scriptsize,
    fill=white,
    inner sep=1.5pt
  }
]
  % The collection is ordered from E_0 at the bottom to E_5 at the top.
  \node[bundle] (E6) at (0,0)
    {$E_5=\mathcal{O}$};
  \node[bundle] (E5) at (0,-1.75)
    {$E_4=\mathcal{O}(-1,-1,-1,1)$};
  \node[bundle] (E4) at (0,-3.50)
    {$E_3=\mathcal{O}(-1,-1,0,0)$};
  \node[bundle] (E3) at (0,-5.25)
    {$E_2=\mathcal{O}(-1,0,-1,0)$};
  \node[bundle] (E2) at (0,-7.00)
    {$E_1=\mathcal{O}(-1,-1,-1,0)$};
  \node[bundle] (E1) at (0,-8.75)
    {$E_0=\mathcal{O}(-2,-1,-1,1)$};

  % Adjacent pairs: the five arrows on the central spine.
  \draw[vanishing arrow] (E6) -- (E5);
  \draw[vanishing arrow] (E5) -- (E4);
  \draw[vanishing arrow] (E4) -- (E3);
  \draw[vanishing arrow] (E3) -- (E2);
  \draw[vanishing arrow] (E2) -- (E1);

  % Pairs separated by one intermediate object: short left-hand arcs.
  \draw[vanishing arrow]
    (E6.west) to[bend right=40] (E4.west);
  \draw[vanishing arrow]
    (E5.west) to[bend right=40] (E3.west);
  \draw[vanishing arrow]
    (E4.west) to[bend right=40] (E2.west);
  \draw[vanishing arrow]
    (E3.west) to[bend right=40] (E1.west);

  % Pairs separated by two intermediate objects: medium right-hand arcs.
  \draw[vanishing arrow]
    (E6.east) to[bend left=50] (E3.east);
  \draw[vanishing arrow]
    (E5.east) to[bend left=50] (E2.east);
  \draw[vanishing arrow]
    (E4.east) to[bend left=50] (E1.east);

  % Pairs separated by three intermediate objects: long left-hand arcs.
  \draw[vanishing arrow]
    (E6.west) to[bend right=62] (E2.west);
  \draw[vanishing arrow]
    (E5.west) to[bend right=62] (E1.west);

  % The outermost pair.
  \draw[vanishing arrow]
    (E6.east) to[bend left=69] (E1.east);

  % Place every label after drawing every arrow, so no curve crosses text.
  \path (E6) -- node[arrow label,right=4pt]
    {$(-1,-1,-1,1)$} (E5);
  \path (E5) -- node[arrow label,right=4pt]
    {$(0,0,1,-1)$} (E4);
  \path (E4) -- node[arrow label,right=4pt]
    {$(0,1,-1,0)$} (E3);
  \path (E3) -- node[arrow label,right=4pt]
    {$(0,-1,0,0)$} (E2);
  \path (E2) -- node[arrow label,right=4pt]
    {$(-1,0,0,1)$} (E1);

  \path (E6.west) to[bend right=40]
    node[arrow label,left=3pt] {$(-1,-1,0,0)$}
    (E4.west);
  \path (E5.west) to[bend right=40]
    node[arrow label,left=3pt] {$(0,1,0,-1)$}
    (E3.west);
  \path (E4.west) to[bend right=40]
    node[arrow label,left=3pt] {$(0,0,-1,0)$}
    (E2.west);
  \path (E3.west) to[bend right=40]
    node[arrow label,left=3pt] {$(-1,-1,0,1)$}
    (E1.west);

  \path (E6.east) to[bend left=50]
    node[arrow label,right=3pt] {$(-1,0,-1,0)$}
    (E3.east);
  \path (E5.east) to[bend left=50]
    node[arrow label,right=3pt] {$(0,0,0,-1)$}
    (E2.east);
  \path (E4.east) to[bend left=50]
    node[arrow label,right=3pt] {$(-1,0,-1,1)$}
    (E1.east);

  \path (E6.west) to[bend right=62]
    node[arrow label,pos=0.36,left=3pt] {$(-1,-1,-1,0)$}
    (E2.west);
  \path (E5.west) to[bend right=62]
    node[arrow label,pos=0.36,left=3pt] {$(-1,0,0,0)$}
    (E1.west);

  \path (E6.east) to[bend left=69]
    node[arrow label,pos=0.42,right=4pt] {$(-2,-1,-1,1)$}
    (E1.east);
\end{tikzpicture}
\caption{Backward Ext-groups for $\mathcal{E}$ on
$X=Bl_{p,q,r}\PP^2$. An arrow $E_i\to E_j$, with $i>j$, labeled
$a=(a_1,a_2,a_3,a_4)$ represents
$\operatorname{Ext}^{\bullet}(E_i,E_j)\cong
H^{\bullet}(X,\mathcal{O}_X(a))=0$.}
\label{fig: BlpqrExc2}
\end{figure}

\begin{figure}[!htbp]
\centering
\begin{tikzpicture}[
  bundle/.style={
    font=\small,
    fill=white,
    inner xsep=5pt,
    inner ysep=2pt
  },
  vanishing arrow/.style={
    ->,
    >=latex,
    line width=0.65pt,
    shorten <=2pt,
    shorten >=2pt
  },
  arrow label/.style={
    font=\scriptsize,
    fill=white,
    inner sep=1.5pt
  }
]
  % The collection is ordered from E_0 at the bottom to E_5 at the top.
  \node[bundle] (E6) at (0,0)
    {$E_5=\mathcal{O}$};
  \node[bundle] (E5) at (0,-1.75)
    {$E_4=\mathcal{O}(-1,-1,-1,1)$};
  \node[bundle] (E4) at (0,-3.50)
    {$E_3=\mathcal{O}(-1,-1,0,0)$};
  \node[bundle] (E3) at (0,-5.25)
    {$E_2=\mathcal{O}(-1,0,-1,0)$};
  \node[bundle] (E2) at (0,-7.00)
    {$E_1=\mathcal{O}(-1,-1,-1,0)$};
  \node[bundle] (E1) at (0,-8.75)
    {$E_0=\mathcal{O}(-2,-1,-1,1)$};

  % Adjacent pairs: the five arrows on the central spine.
  \draw[vanishing arrow] (E5) -- (E6);
  \draw[vanishing arrow] (E4) -- (E5);
  \draw[vanishing arrow] (E3) -- (E4);
  \draw[vanishing arrow] (E2) -- (E3);
  \draw[vanishing arrow] (E1) -- (E2);

  % Pairs separated by one intermediate object.
  \draw[vanishing arrow]
    (E4.west) to[bend left=40] (E6.west);
  \draw[vanishing arrow]
    (E3.west) to[bend left=40] (E5.west);
  \draw[vanishing arrow]
    (E2.west) to[bend left=40] (E4.west);
  \draw[vanishing arrow]
    (E1.west) to[bend left=40] (E3.west);

  % Pairs separated by two intermediate objects.
  \draw[vanishing arrow]
    (E3.east) to[bend right=50] (E6.east);
  \draw[vanishing arrow]
    (E2.east) to[bend right=50] (E5.east);
  \draw[vanishing arrow]
    (E1.east) to[bend right=50] (E4.east);

  % Pairs separated by three intermediate objects.
  \draw[vanishing arrow]
    (E2.west) to[bend left=62] (E6.west);
  \draw[vanishing arrow]
    (E1.west) to[bend left=62] (E5.west);

  % The outermost pair.
  \draw[vanishing arrow]
    (E1.east) to[bend right=69] (E6.east);

  % Draw the labels after the arrows, so no curve crosses any label.
  \path (E5) -- node[arrow label,right=4pt]
    {$(1,1,1,-1)$} (E6);
  \path (E4) -- node[arrow label,right=4pt]
    {$(0,0,-1,1)$} (E5);
  \path (E3) -- node[arrow label,right=4pt]
    {$(0,-1,1,0)$} (E4);
  \path (E2) -- node[arrow label,right=4pt]
    {$(0,1,0,0)$} (E3);
  \path (E1) -- node[arrow label,right=4pt]
    {$(1,0,0,-1)$} (E2);

  \path (E4.west) to[bend left=40]
    node[arrow label,left=3pt] {$(1,1,0,0)$}
    (E6.west);
  \path (E3.west) to[bend left=40]
    node[arrow label,left=3pt] {$(0,-1,0,1)$}
    (E5.west);
  \path (E2.west) to[bend left=40]
    node[arrow label,left=3pt] {$(0,0,1,0)$}
    (E4.west);
  \path (E1.west) to[bend left=40]
    node[arrow label,left=3pt] {$(1,1,0,-1)$}
    (E3.west);

  \path (E3.east) to[bend right=50]
    node[arrow label,right=3pt] {$(1,0,1,0)$}
    (E6.east);
  \path (E2.east) to[bend right=50]
    node[arrow label,right=3pt] {$(0,0,0,1)$}
    (E5.east);
  \path (E1.east) to[bend right=50]
    node[arrow label,right=3pt] {$(1,0,1,-1)$}
    (E4.east);

  \path (E2.west) to[bend left=62]
    node[arrow label,pos=0.64,left=3pt] {$(1,1,1,0)$}
    (E6.west);
  \path (E1.west) to[bend left=62]
    node[arrow label,pos=0.64,left=3pt] {$(1,0,0,0)$}
    (E5.west);

  \path (E1.east) to[bend right=69]
    node[arrow label,pos=0.58,right=4pt] {$(2,1,1,-1)$}
    (E6.east);
\end{tikzpicture}
\caption{Forward higher-Ext groups for $\mathcal{E}$ on
$X=Bl_{p,q,r}\PP^2$. An arrow $E_j\to E_i$, with $i>j$, labeled
$b=(b_1,b_2,b_3,b_4)$ represents
$\operatorname{Ext}^{k}(E_j,E_i)\cong
H^{k}(X,\mathcal{O}_X(b))=0$ for every $k>0$.}
\label{fig: BlpqrStrong}
\end{figure}

\section{Comparing resolutions}

One guiding principle made apparent by the cellular resolution used for
smooth projective toric Fano surfaces in Section~\ref{sec: 4}, as well as
by the resolution of the diagonal due to Hanlon--Hicks--Lazarev
\cite{hanlon2023resolutions}, is the following. A closely related
construction appears in Favero--Huang. For a homotopy path algebra $A$,
they define its path category $\mathcal{C}_A$ and its geometric
realization
\[
X_A:=B(\mathcal{C}_A)
\]
in \cite[\S3.1, Definition~3.7]{favero2022homotopy}. They prove that
$X_A$ supports a projective cellular resolution of the diagonal
$A$-bimodule in
\cite[\S6.1, Corollary~6.7]{favero2022homotopy}. Moreover, for the
Bondal--Ruan homotopy path algebra $A_\Phi$, defined in
\cite[\S5.2, Definition~5.10]{favero2022homotopy}, they prove that
\[
X_{A_\Phi}\simeq \mathbb{T}^m=M_{\mathbb{R}}/M
\]
in \cite[\S5.2, Theorem~5.12]{favero2022homotopy}. These constructions
motivate the following result.

\begin{theorem}\label{thm:5.1}
For a smooth projective toric variety $X_\Sigma$, there exists a cellular resolution $\mathcal{F}_{\mathcal{H}_L/L}^\bullet \rightarrow \mathcal{O}_\Delta$ supported on a real torus $\mathbb{T}^m$.    \end{theorem}

\begin{proof} Up to isomorphism of the fan $\Sigma \subseteq N_\R$, there exists a maximal cone $\sigma\preceq \Sigma$ which can be taken to be $\R_{>0}\left<e_1,\dots e_n\right>$. Thus, the quotient $\faktor{\mathcal{H}_L}{L}$ hyperplane arrangement may be taken inside of a fundemental domain which is a Minkowski sum $\sum_{i=1}^n [0,e_i] = \left[0,1\right]^n$ with opposite pairs of faces identified, creating the torus $\faktor{M_\R}{M}$.\end{proof}

A second proof comes from the construction of the resolution of the diagonal in \cite{hanlon2023resolutions} which is constructed on the same hyperplane arrangement, and our Theorem 5.1 above is proven explicitly for the resolution of the diagonal bimodule in \cite{favero2022homotopy}. \\

Since the resolution of the diagonal constructed in \cite{hanlon2023resolutions} is minimal, we use the Hanlon-Hicks-Lazarev resolution of the diagonal to investigate the following question:

\textbf{Question}: \textit{For which smooth toric Fano threefolds and fourfolds does the Hanlon-Hicks-Lazarev resolution of the diagonal yield a full strong exceptional collection of line bundles?}

\section{Smooth Projective Toric Fano Threefolds}

In this section, we apply the Hanlon-Hicks-Lazarev resolution of $\OO_{\Delta}$ as constructed in Section~\ref{subsec:hhl-resolution} to each of the 18 smooth projective toric Fano threefolds, and check whether the resulting collection of line bundles appearing on one side (say, the left-hand side) factor of the box product in the resolution yields a full strong exceptional collection of line bundles. The results of this section do not require Theorem~\ref{thm:5.1}.

Smooth projective toric Fano threefolds are classified by smooth reflexive
lattice polytopes \cite{C-L-S}, which are available through polymake
\cite{Assarf_Gawrilow_Herr_Joswig_Lorenz_Paffenholz_Rehn_2016,Gawrilow2000}.
Let \(P\subseteq\R^3\) be one of these polytopes, and let \(r\) denote the
number of its facets. Polymake records the half-space representation of
\(P\) as an \(r\times 4\) matrix
\[
H_P=(\mathbf{c}\mid D),
\]
where \(\mathbf{c}\in\R^r\) is the first column and
\(D\in\R^{r\times 3}\) consists of the remaining three columns. If
\(\mathbf{x}=(x_1,x_2,x_3)^{\mathsf T}\), then the row
\[
(c_i,d_{i1},d_{i2},d_{i3})
\]
represents the inequality
\[
c_i+d_{i1}x_1+d_{i2}x_2+d_{i3}x_3\geq 0.
\]
Equivalently,
\[
P=\left\{\mathbf{x}\in\R^3:
\mathbf{c}+D\mathbf{x}\geq 0\right\},
\]
where the inequality is interpreted componentwise. For the reflexive
polytopes considered below, \(\mathbf{c}=\mathbf{1}_r\).

Macaulay2 instead writes a half-space representation in the form
\[
A_{\mathrm{M2}}\mathbf{x}\leq\mathbf{b}_{\mathrm{M2}}.
\]
Consequently, the two conventions are related by
\[
A_{\mathrm{M2}}=-D,
\qquad
\mathbf{b}_{\mathrm{M2}}=\mathbf{c}=\mathbf{1}_r.
\]
Thus, in each polymake matrix displayed below, the first column is the
constant column \(\mathbf{1}_r\), while the final three columns constitute
the matrix \(D\).

Throughout this section, for line bundles \(L\) and \(M\) on \(X\), we write
\[
\operatorname{Ext}^{k}_{X}(L,M)
:=
\operatorname{Hom}_{D^b(X)}(L,M[k]),
\qquad k\in\mathbb{Z}.
\]
For each threefold \(X\) considered below, let \(\mathcal{L}_X\) denote
the set of distinct line bundles occurring on one fixed factor of the
box products in the Hanlon--Hicks--Lazarev resolution of
\(\mathcal{O}_{\Delta}\). After fixing the displayed basis
\(\operatorname{Cl}(X)\cong\mathbb{Z}^{\rho}\), a tuple
\(\mathbf{a}=(a_1,\ldots,a_{\rho})\) appearing in a collection denotes
the line bundle \(\mathcal{O}_X(\mathbf{a})\).

Whenever an ordering of \(\mathcal{L}_X\) is chosen, we denote the
resulting ordered collection by
\[
\mathcal{E}=(E_0,\ldots,E_{N-1}),
\qquad N=\lvert\mathcal{L}_X\rvert,
\]
where the indexing agrees with the displayed order.
The rows and columns of every Hom-dimension matrix below are indexed by
$0,\ldots,N-1$ in this same order.

\subsection{Polytope F.3D.0000}\label{sec:6.1} The polytope $P$ (denoted F.3D.0000 in polymake) has polymake half-space representation

\[ \left[ \begin{matrix} 1 &-1& 1& 2\\
1& -1 &0 &0\\
1 &-1& 0& 1 \\
1 &0 &-1& 0 \\
1 &0& 0& -1 \\
1& 1 &0 &-1 \end{matrix} \right] \] and vertices given by the columns of 

\[  \left( \begin{matrix} 1 & 1 & -4 & 0 & 0 & 1 & 0 & 1 \\ 0 & 1 & 1 & 1 & -3 & -2 & 1 & 1 \\ 0 & 0 & -3 & -1 & 1 & 1 & 1 & 1 \end{matrix}\right). \] 

Let $X$ denote the complete toric variety associated to $P$, with primitive ray generators 
\[ \{ (-1,0,0),(0,-1,0),(0,0,-1),(1,0,-1),(-1,0,1),(-1,1,2)\} \] 

for the fan $\Sigma$ and presentation of class group given by 

\[ \Z^6 \stackrel{ \left( \begin{matrix} 1 & 0 & -1 & 1 & 0 & 0 \\ -1 & 0 & 1 & 0 & 1 & 0 \\ -1 & 1 & 2 & 0 & 0 & 1 \end{matrix}\right) }{\longrightarrow} \Z^3. \]

The Hanlon-Hicks-Lazarev resolution of the diagonal gives the free ranks (written here as ungraded $S$-modules, for $S$ the homogeneous coordinate ring of $Y \cong X \times X$)

\[ 0 \rightarrow S^6 \rightarrow S^{15} \rightarrow S^{11} \rightarrow S^2 \rightarrow 0. \]

Again, this resolution is symmetric in $\pi_1$ and $\pi_2$ for $\pi_i: Y \cong X \times X \rightarrow X$ projection to each copy of $X$. Looking on the left (or right) side, we obtain the collection
\[
\begin{aligned}
\mathcal{E}=(E_0,\ldots,E_8):=\bigl(&
\mathcal{O}_X(-1,0,-2),\mathcal{O}_X(0,-1,-3),
\mathcal{O}_X(-1,0,-1),\\
&\mathcal{O}_X(0,-1,-2),\mathcal{O}_X(0,-1,-1),
\mathcal{O}_X(-1,0,0),\\
&\mathcal{O}_X(0,0,-2),\mathcal{O}_X(0,0,-1),
\mathcal{O}_X\bigr).
\end{aligned}
\]

For every $k\geq 0$, we have
\begin{align*}
\operatorname{Ext}^k_X(E_4,E_0)
&\cong H^k\bigl(X,E_4^{\vee}\otimes E_0\bigr)\\
&\cong H^k\bigl(X,\mathcal{O}_X(-1,1,-1)\bigr)\\
&\cong
\begin{cases}
\C, & k=2,\\
0, & k\neq 2,
\end{cases}
\end{align*}
whereas
\begin{align*}
\operatorname{Ext}^k_X(E_0,E_4)
&\cong H^k\bigl(X,E_0^{\vee}\otimes E_4\bigr)\\
&\cong H^k\bigl(X,\mathcal{O}_X(1,-1,1)\bigr)\\
&\cong
\begin{cases}
\C^6, & k=0,\\
0, & k\neq 0.
\end{cases}
\end{align*}
Thus, no ordering of the members of $\mathcal{E}$ is exceptional. Indeed, if
$E_4$ is placed after $E_0$, then the nonvanishing of
$\operatorname{Ext}^2_X(E_4,E_0)$ violates the semiorthogonality required of
an exceptional collection. If $E_0$ is placed after $E_4$, then the
nonvanishing of $\operatorname{Ext}^0_X(E_0,E_4)$ violates the same condition.
Nevertheless, $\mathcal{E}$ is full, since its members occur on one side of
the locally free resolution of $\mathcal{O}_{\Delta}$.

\subsection{Polytope F.3D.0001} The polytope $P$ (denoted F.3D.0001 in polymake) has polymake half-space representation 

\[ \left[ \begin{matrix} 1 & 1 & 1 & 2 \\ 1 & -1 & 0 & 0 \\ 1 & 0 & 0 & 1 \\ 1 & 0 & -1 & 0 \\ 1 & 0 & 0 & -1 \end{matrix}\right]   \]

and vertices given by the columns of \[ \left( \begin{matrix} 1 & 0 & 1 & 1 & -4 & 1 \\ 0 & 1 & 1 & -4 & 1 & 1 \\ -1 & -1 & -1 & 1 & 1 & 1 \end{matrix}\right). \]

Let $X$ denote the complete toric variety associated to $P$, with primitive ray generators 

\[ \{  (-1,0,0),(0,-1,0),(0,0,-1),(0,0,1),(1,1,2) \} \]

for the fan $\Sigma$ and presentation of the class group given by

\[ \Z^5 \stackrel{ \left( \begin{matrix} 0 & 0 & 1 & 1 & 0 \\ 1 & 1 & 2 & 0 & 1 \end{matrix}\right) }{\longrightarrow} \Z^2. \]

The Hanlon-Hicks-Lazarev resolution of the diagonal gives the free ranks (written as ungraded $S$-modules, for $S$ the homogeneous coordinate ring of $Y \cong X\times X$ ):

\[  0 \rightarrow S^4 \rightarrow S^{11} \rightarrow S^9 \rightarrow S^2 \rightarrow 0. \]

Again, the line bundles appearing on the left-hand side give the collection
\[
\begin{aligned}
\mathcal{E}=(E_0,\ldots,E_6):=\bigl(&
\mathcal{O}_X(-1,-4),\mathcal{O}_X(-1,-3),
\mathcal{O}_X(0,-2),\mathcal{O}_X(-1,-2),\\
&\mathcal{O}_X(0,-1),\mathcal{O}_X(-1,-1),
\mathcal{O}_X\bigr).
\end{aligned}
\]
In particular,
\[
E_0=\mathcal{O}_X(-1,-4)
\qquad\text{and}\qquad
E_5=\mathcal{O}_X(-1,-1).
\]
For every $k\geq 0$, we have
\begin{align*}
\operatorname{Ext}^k_X(E_5,E_0)
&\cong H^k\bigl(X,E_5^{\vee}\otimes E_0\bigr)\\
&\cong H^k\bigl(X,\mathcal{O}_X(0,-3)\bigr)\\
&\cong
\begin{cases}
\C, & k=1,\\
0, & k\neq 1,
\end{cases}
\end{align*}
whereas
\begin{align*}
\operatorname{Ext}^k_X(E_0,E_5)
&\cong H^k\bigl(X,E_0^{\vee}\otimes E_5\bigr)\\
&\cong H^k\bigl(X,\mathcal{O}_X(0,3)\bigr)\\
&\cong
\begin{cases}
\C^{10}, & k=0,\\
0, & k\neq 0.
\end{cases}
\end{align*}
Thus, no ordering of the members of $\mathcal{E}$ is exceptional.
Indeed, if $E_5$ is placed after $E_0$, then
$\operatorname{Ext}^1_X(E_5,E_0)\neq 0$; if $E_0$ is placed after
$E_5$, then $\operatorname{Ext}^0_X(E_0,E_5)\neq 0$.

\subsection{Polytope F.3D.0002} The polytope $P$ (denoted F.3D.0002 in polymake) has half-space representation

\[ \left[ \begin{matrix} 1 &-1& 1& 1\\
1 &-1& 0& 0 \\
1 &-1& 0& 1 \\
1 &0 &-1& 0 \\
1 &0 &0 &1 \\
1 &0 &0 &-1 \\
1 &1 &0 &0 \end{matrix}\right] \] 

and vertices given by the columns of 

\[ \left( \begin{matrix} 1 &1& 0& -1& 0 &-1& -1& 1& -1& 1\\
0 &1 &0 &-1& 1 &1 &-3 &-1 &1 &1\\
0 &0 &-1& -1& -1& -1 &1 &1 &1& 1 \end{matrix} \right). \]

Let $X$ denote the complete toric variety associated to $P$, with primitive ray generators 

\[ \{(-1,0,0),(1,0,0),(0,-1,0),(0,0,-1),(0,0,1),(-1,0,1),(-1,1,1)\}\]

and presentation of the class group given by 

\[ \Z^7 \stackrel{ \left( \begin{matrix} 1 &1 & 0 & 0 & 0 & 0 & 0 \\ 0 & 0 & 0 & 1 & 1 & 0 & 0 \\ -1 & 0 & 0 & 1 & 0 & 1 & 0 \\ -1 & 0 & 1 & 1 & 0 & 0 & 1 \end{matrix}\right)}{\longrightarrow} \Z^4. \]

The Hanlon-Hicks-Lazarev resolution of the diagonal yields the free ranks

\[
0 \rightarrow S^4 \rightarrow S^9 \rightarrow S^6 \rightarrow S^1 \rightarrow 0 \] 

and the full, strong exceptional collection $\{ \mathcal{O}_X(a_1, a_2, a_3, a_4) \}$ for $(a_1, a_2, a_3, a_4)$ in  

\begin{align*} \mathcal{E} &= \{ 
(-1,-1,-1,-2), (-1,-1,-1,-1),(-1,-1,0,-1),(0,-1,-1,-2),(-1,0,0,-1), \dots \\ &\dots (-1,-1,0,0),(0,-1,-1,-1),(-1,0,0,0),(0,0,0,-1),(0,0,0,0)\}. \end{align*} 

\begin{figure}[!htbp]
\includegraphics[width=.5\textwidth]{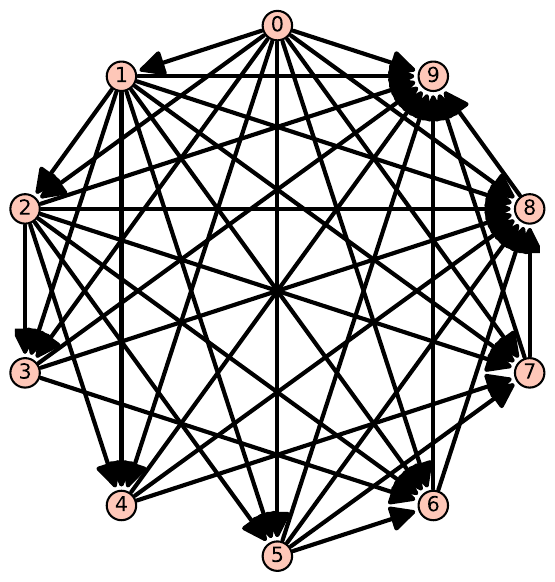}
\caption{Quiver showing nonzero $Hom^0(E_i, E_j)$ for $i\neq j$ in $\mathcal{E}$. Indexing on vertices from $0$ to $9=|\mathcal{E}|-1$.} 
\label{fig: 0002}
\end{figure}

The full quiver is given by the fact that via this ordering on $\mathcal{E}$, $Hom(E_i, E_j)$ is concentrated in degree $0$, with the rank of $Hom(E_j, E_i)$ given by the $(i,j)$ entry of the following matrix:

\[ F =  \left( \begin{matrix} 1 &0& 0& 0& 0& 0& 0& 0& 0& 0 \\ 
2 &1& 0& 0& 0& 0& 0& 0& 0& 0 \\
2 &1& 1& 0& 0& 0& 0& 0& 0& 0\\
3 &1& 1& 1& 0& 0& 0& 0& 0& 0 \\
3 &1& 1& 0& 1& 0& 0& 0& 0& 0 \\
3 &2& 2& 0& 0& 1& 0& 0& 0& 0 \\
5 &3 &2& 2& 0& 1& 1& 0& 0& 0 \\
5 &3& 2& 0& 2& 1& 0& 1& 0& 0 \\
8 &4 &4 &3 &3& 1& 1& 1& 1& 0 \\
12 &8& 7& 5& 5& 4& 3& 3& 2& 1 \end{matrix}\right) \]

i.e., the $(i,j)$ entry of this matrix records the rank of $Hom^0(E_j, E_i)$. The quiver $\mathcal{Q}$ showing nonzero $Hom^0(E_j, E_i)$ for $i\neq j$ is then given by the following directed graph in Figure~\ref{fig: 0002}.

Here, the property that $\mathcal{E}$ is exceptional is visible from the fact that $F$ is lower-triangular, and that starting from vertex $0$ in $\mathcal{Q}$ (which corresponds to $\mathcal{O}_X(-1,-1,-1,-2))$ and moving counter-clockwise, no arrows appear with vertex label of the head less than the vertex label of the tail. 

\subsection{Polytope F.3D.0003}: The polytope $P$ (denoted F.3D.0003 in polymake) has half-space representation

\[ \left[ \begin{matrix}1 &-1& 1& 1 \\
1& -1& 0& 0\\
1 &-1 &0 &1\\
1 &0& -1& 0\\
1 &0 &0 &1\\
1 &0& 0& -1\\
1 &1 &0 &0\\
1 &1& 0& -1 \end{matrix}\right]\]

and vertices given by the columns of

\begin{align*}
\left( \begin{matrix} 
1 & -1 & -1 & 1 & 0 & -1 & 0 & -1 & 0 & 1 & 0 & 1 \\
0 & -2 & 1 & 1 & 0 & -1 & 1 & 1 & -2 & -1 & 1 & 1 \\ 
0 & 0 & 0 & 0 & -1 & -1 & -1 & -1 & 1 & 1 & 1 & 1 \end{matrix}\right).
\end{align*}

Let $X$ denote the complete toric variety associated to $P$, with primitive ray generators 

\[ \{ (-1,0,0),(1,0,0),(0,-1,0),(0,0,-1),(1,0,-1),(0,0,1),(-1,0,1),(-1,1,1) \} \] 

and presentation of the class group given by 

\[ \Z^8 \stackrel{ \left( \begin{matrix} 
1 &1 &0 &0& 0& 0& 0& 0\\
1 &0 &0 &-1& 1& 0& 0& 0 \\
0& 0& 0& 1& 0& 1& 0& 0 \\
-1 &0 &0 &1 &0 &0 &1 &0 \\
-1 &0 &1 &1 &0 &0 &0 & 1 \end{matrix}\right)}{\longrightarrow} \Z^5.  \]

The Hanlon-Hicks-Lazarev resolution of the diagonal yields the free ranks

\[ 0 \rightarrow S^4 \rightarrow S^9 \rightarrow S^6 \rightarrow S^1 \rightarrow 0 \] 

and the full, strong exceptional collection $\{\mathcal{O}_X(a_1, a_2, a_3, a_4, a_5)\}$ of 12 line bundles for $(a_1,a_2,a_3,a_4,a_5)$ in 

\begin{align*}
    \mathcal{E} = &\{ (-1,0,-1,-1,-2), (-1,-1,-1,0,-1), (-1,-1,0,0,-1),(-1,0,-1,-1,-1), \dots\\
    &\dots (-1,0,-1,0,-1), (0,0,-1,-1,-2), (-1,-1,-1,0,0), (-1,0,-1,0,0), (0,0,-1,-1,-1), (-1,-1,0,0,0),\dots \\  &\dots (0,0,0,0,-1), (0,0,0,0,0)\}. 
\end{align*}

Here, $Hom_{D^b(X)}(E_i, E_j)$ is concentrated in degree $0$ for all $i$ and $j$, with the rank of $Hom^0(E_j, E_i)$ given by the $(i,j)$ entry of the following matrix:

\[ F = 
\left(
\begin{matrix} 
1 &0& 0& 0& 0& 0& 0& 0& 0& 0& 0& 0\\
2 &1 &0 &0 &0 &0 &0 &0 &0 &0 &0 &0 \\
3 &1& 1& 0& 0& 0& 0& 0& 0& 0& 0& 0 \\
3 &1 &0 &1 &0 &0 &0 &0 &0 &0 &0 &0 \\
2 &0& 0 &0 &1 &0 &0 &0 &0 &0 &0 &0 \\
4 &1 &1 &1 &1 &1 &0 &0 &0 &0 &0 &0 \\
5 &3 &0 &2 &0 &0 &1 &0 &0 &0 &0 &0 \\
4 &2 &0 &0 &2 &0 &0 &1 &0 &0 &0 &0 \\
5 &2 &1 &1 &2 &0 &0 &1 &1 &0 &0 &0 \\
5 &3 &2 &0 &0 &0 &0 &0 &0 &1 &0 &0 \\
7 &4 &2 &2 &2 &2 &1 &1 &0 &1 &1 &0\\
8 &5 &2 &2 &3 &0 &1 &2 &2 &1 &0 &1 \end{matrix} \right). \]

\begin{figure}[!htbp]
\includegraphics[width=.5\textwidth]{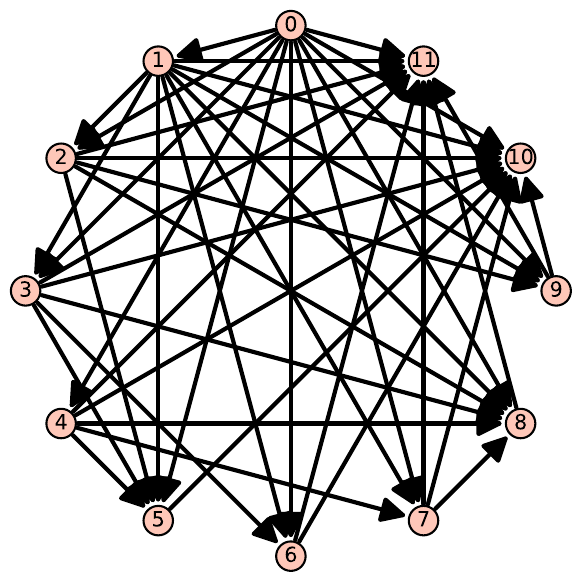}
\caption{Quiver showing nonzero $Hom^0(E_i, E_j)$ for $i\neq j$ in $\mathcal{E}$. Indexing on vertices from $0$ to $11=|\mathcal{E}|-1$.} 
\label{fig:15}
\end{figure}

The fact $F$ is lower-triangular shows that $\mathcal{E}$ is exceptional. The quiver $\mathcal{Q}$ showing nonzero $Hom^0(E_j, E_i)$ for $i\neq j$ is then given by the directed graph in Figure~\ref{fig:15}. Starting from vertex $0$ in $\mathcal{Q}$ (which corresponds to $\mathcal{O}_X(-1,0,-1,-1,-2))$ and moving counter-clockwise, no arrows appear with vertex label of the head less than the vertex label of the tail.

\subsection{Polytope F.3D.0004}: The polytope $P$ (denoted F.3D.0004 in polymake) has half-space representation 

\[ \left[  \begin{matrix} 1 &-1 & 1 & 1 \\
1 & -1&  0 & 0 \\ 1 & -1 &0  &1 \\ 1 &0 &-1 &0
 \\ 1 &0 &0 &1 \\ 1 &0 &0 &-1
 \\ 1 &1 &0 &-1\end{matrix} \right] \] 

 and vertices given by the columns of 

 \[ \left( \begin{matrix} 
 1 &1& 0& -2& 0& -2& 0& 1 &0 &1 \\
 0 &1 &0 &-2 &1 &1 &-2 &-1 &1& 1 \\
 0 &0& -1& -1 &-1& -1& 1 &1 &1 &1 \end{matrix} \right). \] 

Let $X$ denote the complete toric variety associated to $P$, with primitive ray generators 

\[ \{ (-1,0,0),(0,-1,0),(0,0,-1),(1,0,-1),(0,0,1), (-1,0,1),(-1,1,1) \} \] 

and presentation of the class group given by \[ \Z^7 \stackrel{ \left( \begin{matrix} 1 & 0 & -1 & 1 & 0 & 0 & 0 \\ 0 & 0 & 1 & 0 & 1 & 0 & 0 \\ -1 & 0 & 1 & 0 & 0 & 1 & 0 \\ -1 & 1 & 1 & 0 & 0 & 0 & 1 \end{matrix}\right)        }{\longrightarrow} \Z^4.  \]

The Hanlon-Hicks-Lazarev resolution of the diagonal yields the free ranks (written as ungraded S-modules,
for S the homogeneous coordinate ring of $Y \cong X \times X )$ \[ 
 0 \rightarrow S^4 \rightarrow S^9 \rightarrow S^6 \rightarrow S^1 \rightarrow 0 \] where here $S$ denotes the homogeneous coordinate ring for $Y = X\times X$. Here we have a full strong exceptional collection of 10 line bundles $\{ \mathcal{O}_X(a_1, a_2, \dots a_5) \}$ for $(a_1, \dots, a_5)$ in 

\begin{align*}
\mathcal{E} &= \{ (-1,-1,0,-1), (0,-1,-1,-2), (-1,0,0,-1), (-1,-1,0,0), (0,-1,-1,-1), (0,-1,0,-1), (0,-1,0,0), \dots \\
&\dots (-1,0,0,0), (0,0,0,-1), (0,0,0,0) \}. \end{align*} 

Here, $Hom_{D^b(X)}(E_i, E_j)$ is concentrated in degree $0$ for all $i$ and $j$, with the rank of $Hom^0(E_j, E_i)$ given by the $(i,j)$ entry of the following matrix:

\[ 
F =  \left( \begin{matrix} 1&0 &0& 0& 0& 0& 0& 0& 0& 0 \\
2 &1 &0 &0 &0 &0 &0 &0 &0 &0 \\
3& 1& 1& 0& 0& 0& 0& 0& 0& 0 \\
2& 0& 0& 1& 0& 0& 0& 0& 0& 0 \\
4 &2 &0 &2 &1 &0 &0 &0 &0 &0 \\
5 &2 &1 &2 &1 &1 &0 &0 &0 &0 \\
6 &2 &2 &3 &1 &1 &1 &0 &0 &0 \\
5 &3 &2 &0 &0 &0 &0 &1 &0 &0 \\
8& 5 &2 &3 &2 &2 &0 &1 &1 &0 \\
10 &6 &4 &5 &3 &2 &2 &2 &1 &1 \end{matrix} \right).   \]

\begin{figure}[!htbp]
\includegraphics[width=.5\textwidth]{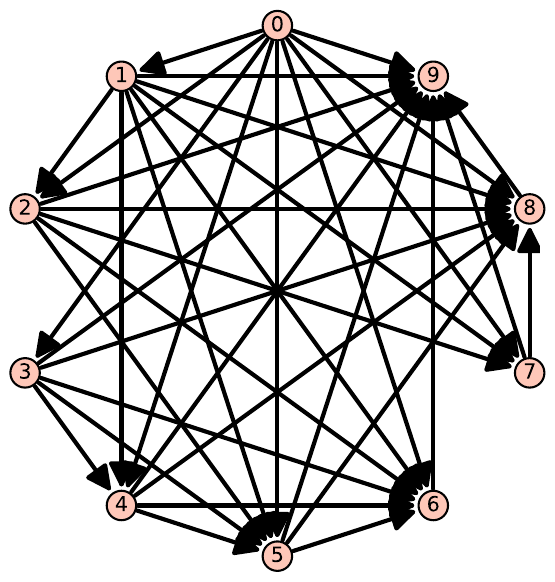}
\caption{Quiver showing nonzero $Hom^0(E_i, E_j)$ for $i\neq j$ in $\mathcal{E}$. Indexing on vertices from $0$ to $9=|\mathcal{E}|-1$.} 
\label{fig: 16}
\end{figure}

The fact $F$ is lower-triangular shows that $\mathcal{E}$ is exceptional. The quiver $\mathcal{Q}$ showing nonzero $Hom^0(E_j, E_i)$ for $i\neq j$ is then given by the directed graph in Figure~\ref{fig: 16}. Starting from vertex $0$ in $\mathcal{Q}$ (which corresponds to $\mathcal{O}_X( -1,-1,0,-1 ))$ and moving counter-clockwise, no arrows appear with vertex label of the head less than the vertex label of the tail.

\subsection{Polytope F.3D.0005} The polytope $P$ (denoted F.3D.0005 in polymake) has half-space representation

\[ \left[ \begin{matrix} 1 & 0&  1& 1\\
1 &-1& 0& 0\\
1 & 0& 0& 1\\
1 &0 &-1& 0\\
1& 1 &0 &1\\
1 &0 &0 &-1 \end{matrix} \right] \] 

and vertices given by the columns of 

\[ \left( \begin{matrix} 0 &1 &0 &1 &-2& 1& -2 &1 \\
0 &0 &1 &1 &-2& -2& 1& 1\\
-1 &-1& -1& -1& 1 &1 &1 &1 \end{matrix} \right). \]

Let $X$ denote the complete toric variety associated to $P$, with primitive ray generators 

\[ \{ (-1,0,0),(0,-1,0),(0,0,-1),(0,0,1),(1,0,1),(0,1,1) \} \] 

and presentation of the class group given by 

\[ \Z^6 \stackrel{ \left(\begin{matrix} 0 & 0 &1 &1 &0 &0\\
1 &0 &1 &0 &1& 0\\
0 &1 &1 &0 &0 &1 \end{matrix}\right)}{\longrightarrow} \Z^3. \]  

The Hanlon-Hicks-Lazarev resolution of the diagonal yields the free ranks (written as ungraded $S$-modules, for $S$ the homogeneous coordinate ring of $Y\cong X\times X$) 

\[ 0 \rightarrow S^4 \rightarrow S^9 \rightarrow S^6 \rightarrow S^1 \rightarrow 0 \] 

and the collection of 8 line bundles $\mathcal{E}$ which appear on, say, the left-hand side are $\mathcal{O}(a_1,a_2,a_3)$ for $(a_1,a_2,a_3)$ in

\begin{align*} 
\mathcal{E} = \{ (-1,-2,-2),(-1,-1,-2) ,(-1,-2,-1)  ,(-1,-1,-1), {0,-1,-1}, (0,-1,0) ,  (0,0,-1) ,  (0,0,0) \}. \end{align*}

$Hom_{D^b(X)}^\bullet (E_i, E_j)$ is concentrated in degree $0$ for all $i$ and $j$ here, with the rank of $Hom^0_{D^b(X)}(E_j, E_i)$ given by the $(i,j)$ entry of the following matrix:

\[ F = \left( \begin{matrix} 
1 &0& 0& 0& 0& 0& 0& 0\\
2 &1 &0 &0 &0 &0 &0 &0\\
2 &0 &1 &0 &0 &0 &0 &0\\
4 &2 &2 &1 &0 &0 &0 &0\\
5 &2 &2 &1 &1 &0 &0 &0\\
8 &3 &5 &2 &2 &1 &0 &0\\
8 &5 &3 &2 &2 &0 &1 &0\\
13 &8 &8 &5 &4 &2 &2 &1
\end{matrix}\right). \]

\begin{figure}[!htbp]
\includegraphics[width=.5\textwidth]{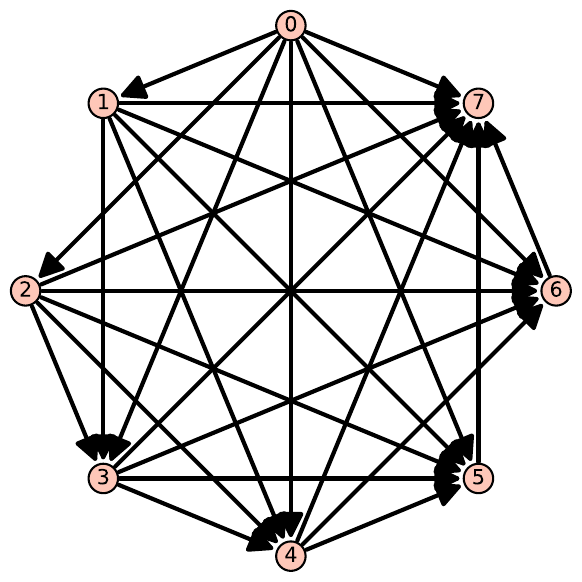}
\caption{Quiver showing nonzero $Hom^0(E_i, E_j)$ for $i\neq j$ in $\mathcal{E}$. Indexing on vertices from $0$ to $7=|\mathcal{E}|-1$.} 
\label{fig: 17}
\end{figure}

The fact that $F$ is lower-triangular shows that $\mathcal{E}$ is exceptional. The quiver $\mathcal{Q}$ showing nonzero $Hom^0(E_j, E_i)$ for $i\neq j$ is then given by the directed graph in Figure~\ref{fig: 17}. Starting from vertex $0$ in $\mathcal{Q}$ (which corresponds to $\mathcal{O}_X( -1,-1,-2)$ and moving counter-clockwise, no arrows appear with vertex label of the head less than the vertex label of the tail.

\subsection{Polytope F.3D.0006} The polytope $P$ (denoted F.3D.0006 in polymake) has half-space representation

\[ \left[ \begin{matrix} 1& 0 &-1& 1 \\ 1 & -1 &0& 0 \\ 1& 0& -1& 0\\ 1 &0 &0& -1 \\  1& -1& 1& 0 \\1 & 1 &0 &0 \end{matrix}\right] \] 

and vertices given by the columns of 

\[ \left( \begin{matrix} 
-1 &1& -1& 1 &1 &-1& -1& 1 \\
1 &1 &-2& 0& 0& -2& 1& 1 \\
0 &0 &-3& -1& 1& 1& 1& 1 \end{matrix} \right). \]

Let $X$ denote the complete toric variety associated to $P$, with primitive ray generators of $\Sigma$ \[ \{ (-1,0,0), (1,0,0), (0,-1,0), (-1,1,0), (0,0,-1), (0,-1,1) \}  \] and presentation of the class group given by

\[ \Z^6 \stackrel{ \left( \begin{matrix} 1 &1 &0 &0 &0 &0 \\ -1 &0& 1& 1& 0& 0 \\ -1 &0 &0 &1 &1 &1 \end{matrix} \right) }{\longrightarrow} \Z^3. \]

The Hanlon-Hicks-Lazarev resolution of the diagonal yields the free ranks (written as ungraded $S$-modules, for $S$ the homogeneous coordinate ring of $Y\cong X \times X$)

\[ 0 \rightarrow S^4 \rightarrow S^9 \rightarrow S^6 \rightarrow S^1 \rightarrow  0     \]

and the collection of $8$ line bundles $\mathcal{E}$ which appear on the left-hand side are $\mathcal{O}(a_1, a_2, a_3)$ for $(a_1,a_2,a_3)$ appearing in

\[ \mathcal{E} = \{ (-1,-1,-2),(-1,-1,-1),(-1,0,-1),(0,-1,-2),(-1,0,0), (0,-1,-1),(0,0,-1),(0,0,0) \}. \]

Here, $Hom_{D^b(X)}^\bullet(E_i, E_j)$ is concentrated in degree $0$ for all $i$ and $j$, with the rank of $Hom^0(E_j, E_i)$ given by the $(i,j)$ entry of the following matrix:

\[ F = \left( \begin{matrix} 
1 &0& 0& 0& 0& 0& 0& 0\\
2 &1 &0 &0 &0 &0 &0& 0\\
3 &1 &1 &0 &0 &0 &0 &0\\
4 &1 &1 &1 &0 &0 &0 &0\\
5 &3 &2 &0 &1 &0 &0 &0\\
7 &4 &2 &2 &1 &1 &0 &0\\
9 &4 &4 &3 &1 &1 &1 &0\\
14 &9 &7 &5 &4 &3 &2 &1\end{matrix} \right). \]

\begin{figure}[!htbp]
\includegraphics[width=.5\textwidth]{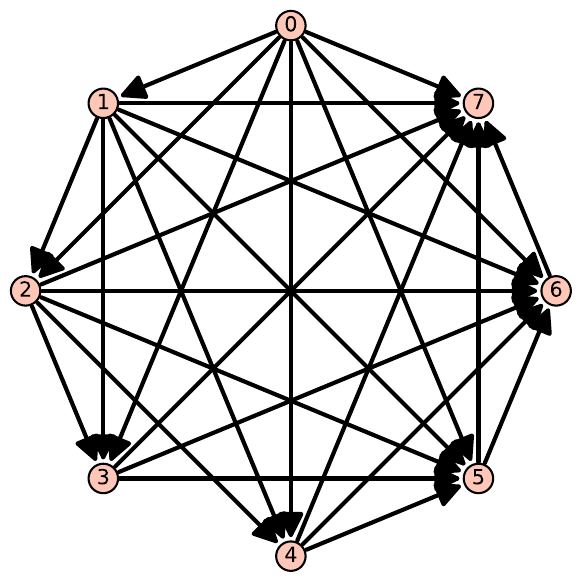}
\caption{Quiver showing nonzero $Hom^0(E_i, E_j)$ for $i\neq j$ in $\mathcal{E}$. Indexing on vertices from $0$ to $7=|\mathcal{E}|-1$.} 
\label{fig: 18}
\end{figure}

The fact that $F$ is lower-triangular shows that $\mathcal{E}$ is exceptional. The quiver $\mathcal{Q}$ showing nonzero $Hom^0(E_j, E_i)$ for $i\neq j$ is then given by the directed graph in Figure~\ref{fig: 18}. Starting from vertex $0$ in $\mathcal{Q}$ (which corresponds to $\mathcal{O}_X( -1,-1,-2)$ ) and moving counter-clockwise, no arrows appear with vertex label of the head less than the vertex label of the tail.

\subsection{Polytope F.3D.0007} The polytope $P$ (denoted F.3D.0007 in polymake) has half-space representation

\[ \left[ \begin{matrix} 1 & -1 & 1 & 1 \\ 1 & -1 & 0 & 0 \\ 1 & -1 & 0 & 1 \\ 1 & 0 & -1 & 0 \\ 1 & 0 & 0 & -1 \\ 1 & 1 & 0 & 0 \\ 1 & 1 & 0 & -1 \end{matrix}\right]  \]

and vertices given by the columns of 

\[ \left( \begin{matrix} 1 &-1& -1& 1& -1& -1& 0& 1& 0& 1 \\ 0 &-2& 1& 1& 0& 1& -2& -1& 1& 1 \\ 0& 0& 0& 0& -2& -2& 1 &1 &1 &1 \end{matrix} \right). \]

Let $X$ denote the complete toric variety associated to $P$, with primitive ray generators of $\Sigma$ 

\[ \{ (-1,0,0),(1,0,0),(0,-1,0),(0,0,-1),(1,0,-1),(-1,0,1),(-1,1,1) \}  \]

and presentation of the class group given by 

\[  \Z^7 \stackrel{ \left( \begin{matrix} 1& 1& 0& 0& 0& 0& 0\\ 1 &0 &0 &-1& 1& 0& 0 \\ -1& 0& 0& 1& 0& 1& 0 \\ -1 &0 &1 &1 &0 &0 &1 \end{matrix} \right)    }{\longrightarrow} \Z^4. \]

The Hanlon-Hicks-Lazarev resolution of the diagonal yields the free ranks (written as ungraded $S$-modules, for $S$ the homogeneous coordinate ring of $Y \cong X \times X$)

 \[0 \rightarrow S^4 \rightarrow S^9 \rightarrow S^6 \rightarrow S^1 \rightarrow 0 \] 

and the collection of 10 line bundles 

\begin{align*} \mathcal{E} &= \{  (-1,0,-1,-2), (-1,-1,0,-1),(-1,0,-1,-1),(-1,0,0,-1),(0,0,-1,-2),(-1,-1,0,0), \dots \\
&\dots (-1,0,0,0),(0,0,-1,-1),(0,0,0,-1),(0,0,0,0)         \} \end{align*} 

is full, strong, and exceptional. That is, $Hom_{D^b(X)}^\bullet(E_i, E_j)$ is concentrated in degree $0$ for all $i$ and $j$, with the rank of $Hom^0(E_j, E_i)$ given by the $(i,j)$ entry of the following matrix:

\[ F = \left( \begin{matrix}
1 &0 &0 &0 &0 &0 &0 &0 &0 &0\\
2 &1 &0 &0 &0 &0 &0 &0 &0 &0\\
3 &1 &1 &0 &0 &0 &0 &0 &0 &0\\
2 &0 &0 &1 &0 &0 &0 &0 &0 &0\\
4 &1 &1 &1 &1 &0 &0 &0 &0 &0\\
5 &3 &2 &0 &0 &1 &0 &0 &0 &0\\
4& 2& 0& 2& 0& 0& 1& 0& 0& 0\\
6 &2 &2 &3 &1 &0 &1 &1 &0 &0\\
7 &4 &2 &2 &2 &1 &1 &0 &1 &0 \\
10 &6 &4 &5 &2 &2 &3 &2 &1 &1\end{matrix} \right). \] 

\begin{figure}[!htbp]
\includegraphics[width=.5\textwidth]{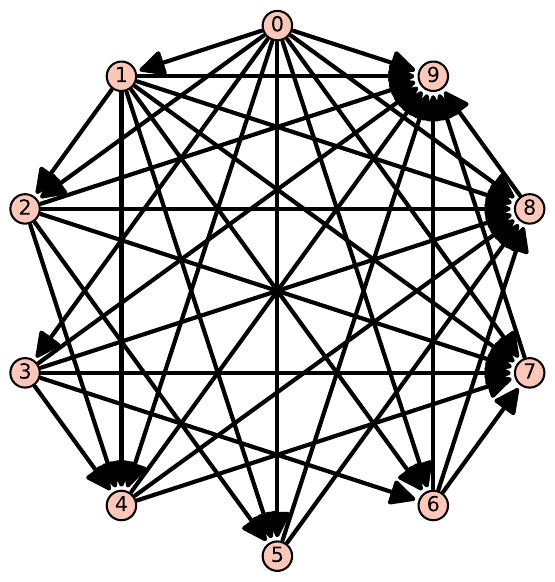}
\caption{Quiver showing nonzero $Hom^0(E_i, E_j)$ for $i\neq j$ in $\mathcal{E}$. Indexing on vertices from $0$ to $9=|\mathcal{E}|-1$.} 
\label{fig: 19}
\end{figure}

The fact that $F$ is lower-triangular shows that $\mathcal{E}$ is exceptional. The quiver $\mathcal{Q}$ showing nonzero $Hom^0(E_j, E_i)$ for $i\neq j$ is then given by the directed graph in Figure~\ref{fig: 19}. Starting from vertex $0$ in $\mathcal{Q}$ (which corresponds to $\mathcal{O}_X( -1,0,-1,-2)$ ) and moving counter-clockwise, no arrows appear with vertex label of the head less than the vertex label of the tail. 

\subsection{Polytope F.3D.0008} The polytope $P$ (denoted F.3D.0008 in polymake) has half-space representation  

\[ \left[ \begin{matrix} 1 &0& 0& 1\\
1 &-1& 0& 0\\\
1 &0& 1& 1\\
1 &0 &1& 0\\
1 &0 &-1& 0\\
1 &0 &0 &-1\\
1 &1 &0 &0 \end{matrix} \right] \] 

and vertices given by the columns of

\[ \left( \begin{matrix}   
-1 &1 &-1& 1& -1& 1& -1& 1& -1& 1\\
-1 &-1& 0& 0& 1& 1& -1& -1& 1 &1\\
0 &0 &-1& -1& -1& -1& 1 &1 &1 &1\end{matrix}\right). \] 

Let $X$ denote the complete toric variety associated to $P$, with primitive ray generators of $\Sigma$ 

\[ \{ (-1,0,0), (1,0,0),(0,-1,0),(0,1,0), (0,0,-1), (0,0,1), (0,1,1) \} \]

and presentation of the class group given by 

\[ \Z^7 \stackrel{ \left( \begin{matrix} 1 &1 &0 &0 &0 &0 &0 \\ 0 &0 &1 &1 &0 &0 &0 \\ 0& 0& 0& 0& 1& 1& 0 \\ 0 &0 &1 &0 &1 &0 &1 \end{matrix}\right)   }{\longrightarrow} \Z^4. \]

The Hanlon-Hicks-Lazarev resolution of the diagonal yields the free ranks (written as ungraded $S$-modules, for $S$ the homogeneous coordinate ring of $Y \cong X \times X)$ \[ 0 \rightarrow S^2 \rightarrow S^5 \rightarrow S^4 \rightarrow S^1 \rightarrow 0. \] 

The collection of 10 line bundles which appear on the left-hand side are $\mathcal{O}(a_1,a_2,a_3,a_4)$ for $(a_1, a_2, a_3, a_4)$ appearing in

\begin{align*}   
\mathcal{E} &= \{ (-1,-1,-1,-2), (-1,-1,-1,-1), (-1,0,-1,-1),(-1,-1,0,-1),(0,-1,-1,-2), \dots \\ &\dots (0,-1,-1,-1), (0,0,-1,-1), (0,-1,0,-1), (-1,0,0,0), (0,0,0,0) \} \end{align*}

is full, strong, and exceptional. That is, $Hom_{D^b(X)}^\bullet(E_i, E_j)$ is concentrated in degree $0$ for all $i$ and $j$, with the rank of $Hom^0(E_j, E_i)$ given by the $(i,j)$ entry of the following matrix:

\[ F = \left( \begin{matrix}
1 &0 & 0 & 0 & 0 & 0 & 0 & 0 & 0  &0\\
2 & 1 & 0 & 0 & 0 & 0 & 0 & 0 & 0 & 0\\
2 & 0 & 0 & 1 & 0 & 0 & 0 & 0 & 0 & 0\\
3  &0  &1  &1  &1  &0  &0  &0  &0  &0\\
4  &0  &2  &2  &1  &1  &0  &0  &0  &0\\
4 & 2 & 2 & 0 & 0 & 0 & 1 & 0 & 0  &0\\
4 & 2 & 0 & 2 & 0 & 0 & 0 & 1 & 0 & 0\\
6 & 3  &2  &2  &2  &0  &1  &1  &1 &0\\
8  &4  &4  &4  &2  &2  &2  &2  &1 & 1 \end{matrix}\right). \] 

\begin{figure}[!htbp]
\includegraphics[width=.5\textwidth]{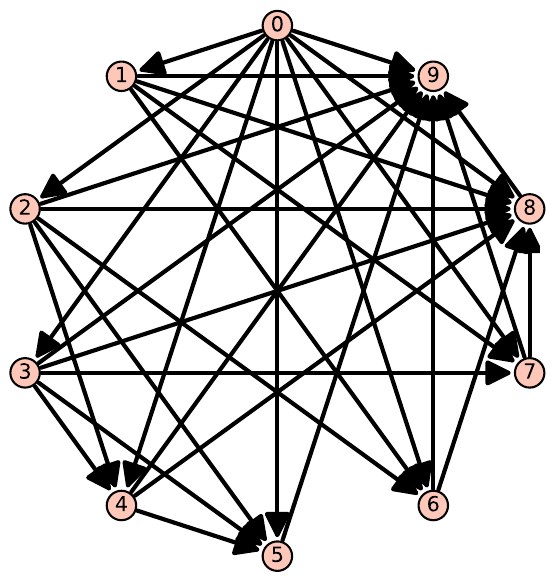}
\caption{Quiver showing nonzero $Hom^0(E_i, E_j)$ for $i\neq j$ in $\mathcal{E}$. Indexing on vertices from $0$ to $9=|\mathcal{E}|-1$.} 
\label{fig: 20}
\end{figure}

The fact that $F$ is lower-triangular shows that $\mathcal{E}$ is exceptional. The quiver $\mathcal{Q}$ showing nonzero $Hom^0(E_j, E_i)$ for $i\neq j$ is then given by the directed graph in Figure~\ref{fig: 20}. Starting from vertex $0$ in $\mathcal{Q}$ (which corresponds to $\mathcal{O}_X( -1,-1,-1,-2)$ ) and moving counter-clockwise, no arrows appear with vertex label of the head less than the vertex label of the tail. 

\subsection{Polytope F.3D.0009} The polytope $P$ (denoted F.3D.0009 in polymake) has half-space representation

\[ \left[ \begin{matrix} 1 &0 &-1& 1\\
1 &-1& 0& 0\\
1 &0 &0 &1 \\
1 &0 &1 &0 \\
1 &0 &-1& 0 \\
1 &0 &1 &-1 \\
1 &1 &0 &0 \\
1 &0 &0 &-1 \end{matrix} \right] \] and vertices given by the columns of 

\begin{align*}
    \left( \begin{matrix}   
    -1 &1 &-1 &1 &-1 &1 &-1 &1 &-1 &1 &-1 &1\\
-1&-1&1&1&0&0&-1&-1&0&0&1&1\\
0&0&0&0&-1&-1&-1&-1&1&1&1&1  \end{matrix} \right).
\end{align*}

Let $X$ denote the complete toric variety associated to $P$, with primitive ray generators of $\Sigma$ 

\[ \{(-1,0,0),(1,0,0), (0,-1,0), (0,1,0), (0,0,-1), (0,1,-1),(0,0,1), (0,-1,1) \} \] 

and presentation of the class group given by 

\[ \Z^8 \stackrel{ \left( \begin{matrix} 1 &1 &0 &0 &0 &0 &0 &0 \\
0 &0 &1 &1 &0 &0 &0 &0\\
0 &0 &1 &0 &-1 &1 &0 &0\\
0 &0 &0 &0 &1 &0 &1 &0\\
0 &0 &-1 &0 &1 &0 &0 &1 \end{matrix}\right)}{\longrightarrow} \Z^5. \]

The Hanlon-Hicks-Lazarev resolution of the diagonal yields the free ranks (written as ungraded $S$-modules, for $S$ the homogeneous coordinate ring of $Y \cong X \times X$):

\[  0 \rightarrow S^2 \rightarrow S^5 \rightarrow S^4 \rightarrow S^1 \rightarrow 0. \]

The collection of $12$ line bundles $\mathcal{E}$ which appear on the left-hand side are $\mathcal{O}(a_1, a_2, a_3,a_4, a_5)$ for $(a_1,a_2,a_3, a_4, a_5)$ appearing in

\begin{align*} \mathcal{E} &=  \{ (-1,-1,-1,-1,0), (-1,-1,0,-1,-1), (-1,-1,-1,0,0), (-1,0,0,-1,-1), (-1,-1,0,-1,0), \dots \\
& \dots (0,-1,-1,-1,0), (0,-1,0,-1,-1), (0,0,0,-1,-1), (0,-1,-1,0,0), (0,-1,0,-1,0), (-1,0,0,0,0), (0,0,0,0,0) \} \end{align*} 

is full, strong, and exceptional. That is, $Hom_{D^b(X)}^\bullet(E_i, E_j)$ is concentrated in degree $0$ for all $i$ and $j$, with the rank of $Hom^0(E_j, E_i)$ given by the $(i,j)$ entry of the following matrix:

\[ F = \left( \begin{matrix} 
1 &0 &0 &0 &0 &0 &0 &0 &0 &0 &0 &0 \\
2 &1 &0 &0 &0 &0 &0 &0 &0 &0 &0 &0\\
2 &0 &1 &0 &0 &0 &0 &0 &0 &0 &0 &0\\
2 &0 &0 &1 &0 &0 &0 &0 &0 &0 &0 &0\\
2 &0 &0 &0 &1 &0 &0 &0 &0 &0 &0 &0 \\
3 &0 &1 &1 &1 &1 &0 &0 &0 &0 &0 &0\\
3 &0 &1 &1 &1 &0 &1 &0 &0 &0 &0 &0\\
4 &2 &2 &0 &0 &0 &0 &1 &0&0 &0 &0 \\
4 &2 &0 &0 &2 &0 &0 &0 &1 &0 &0 &0\\
4 &2 &0 &2 &0 &0 &0 &0 &0 &1 &0 &0\\
6 &3 &2 &2 &2 &2 &0&1 &1 &1 &1 &0\\
6 &3 &2 &2 &2 &0 &2 &1 &1 &1 &0 &1
\end{matrix} \right). \]

\begin{figure}[!htbp]
\includegraphics[width=.5\textwidth]{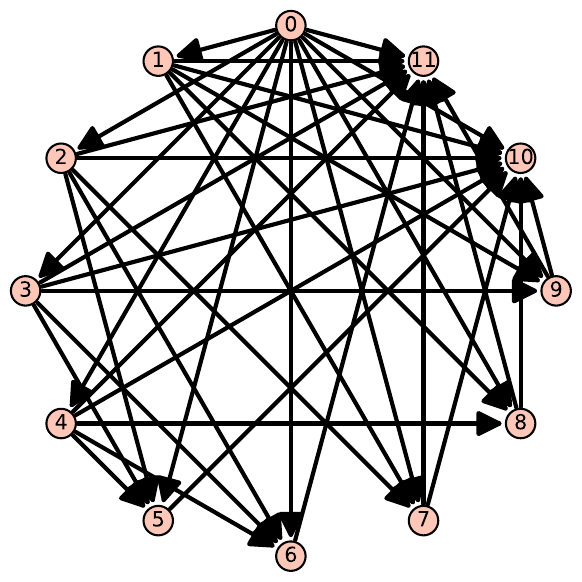}
\caption{Quiver showing nonzero $Hom^0(E_i, E_j)$ for $i\neq j$ in $\mathcal{E}$. Indexing on vertices from $0$ to $11=|\mathcal{E}|-1$.} 
\label{fig: 21}
\end{figure}

The fact that $F$ is lower-triangular shows that $\mathcal{E}$ is exceptional. The quiver $\mathcal{Q}$ showing nonzero $Hom^0(E_j, E_i)$ for $i\neq j$ is then given by the directed graph in Figure~\ref{fig: 21}. Starting from vertex $0$ in $\mathcal{Q}$ (which corresponds to $\mathcal{O}_X( -1,-1,-1,-1,0)$ and moving counter-clockwise, no arrows appear with vertex label of the head less than the vertex label of the tail.

\subsection{Polytope F.3D.0010} The polytope $P$ (denoted F.3D.0010 in polymake) has half-space representation

\[ \left[ \begin{matrix} 1 &1 &0& 1\\
1 &-1& 0& 0 \\
1 &1& 1& 1 \\
1 &0 &1 &0 \\
1 &0 &-1& 0 \\
1& 0& 0& -1 \end{matrix} \right] \] 

and vertices given by the columns of

\[ \left( \begin{matrix}
1 &1 &1 &-2 &-1 &1 &-2 &1\\
0 &1 &-1 &0 &-1 &-1 &1 &1\\
-2&-2 &-1 &1 &1 &1 &1 &1 \end{matrix}\right). \]  

Let $X$ denote the complete toric variety associated to $P$, with primitive ray generators of $\Sigma$ \[ \{  (-1,0,0),(0,-1,0), (0,1,0), (0,0,-1),(1,0,1), (1,1,1) \} \]

and presentation of the class group given by 

\[ \Z^6 \stackrel{\left( \begin{matrix} 0 &1 &1 &0 &0 &0 \\
1 &0 &0 &1 &1 &0\\
1&1 &0 &1 &0 &1 \end{matrix} \right)}{\longrightarrow} \Z^3. \] 

The Hanlon-Hicks-Lazarev resolution of the diagonal yields the free ranks (written as ungraded S-modules, for $S$ the homogeneous coordinate ring of $Y \cong X \times X$):\\

\[ 0 \rightarrow S^4 \rightarrow S^9 \rightarrow S^6 \rightarrow S^1 \rightarrow 0.\]

The collection of 8 line bundles $\mathcal{E}$ which appear on the left-hand side are $\mathcal{O}(a_1, a_2, a_3)$ for $(a_1, a_2, a_3)$ appearing in 

\[ \mathcal{E} = \{ (-1,-2,-3), (-1,-2,-2), (-1,-1,-2), (0,-2,-2), (-1,-1,-1),(-1,0,-1), (0,-1,-1), (0,0,0)\} \]

is full, strong, and exceptional. That is, $Hom_{D^b(X)}^\bullet(E_i, E_j)$ is concentrated in degree $0$ for all $i$ and $j$, with the rank of $Hom^0(E_j, E_i)$ given by the $(i,j)$ entry of the following matrix:

\[ F = \left( \begin{matrix} 1 &0 &0 &0 &0 &0 &0 &0\\    
3 &1 &0 &0 &0 &0 &0 &0\\
2 &0 &1 &0 &0 &0 &0 &0\\
4 & 1 &1 &1 &0 &0 &0 &0\\
6 &3 &0 &0 &1 &0 &0 &0\\
6 &2 &3 &1 &0 &1 &0 &0\\
9 &4 &3 &3 &1 &1 &1 &0 \\
12 &6 &6 &3 &2 &3 &1 &1 \end{matrix} \right). \]

\begin{figure}[!htbp]
\includegraphics[width=.5\textwidth]{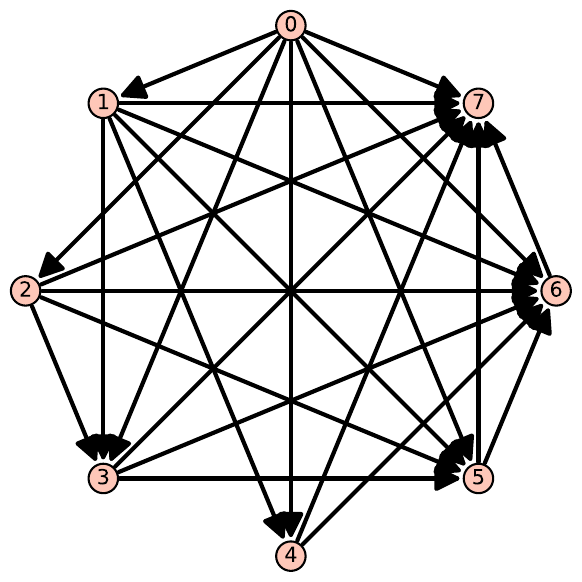}
 \caption{Quiver showing nonzero $Hom^0(E_i, E_j)$ for $i\neq j$ in $\mathcal{E}$. Indexing on vertices from $0$ to $7=|\mathcal{E}|-1$.} 
\label{fig: 22}
\end{figure}

The fact that $F$ is lower-triangular shows that $\mathcal{E}$ is exceptional. The quiver $\mathcal{Q}$ showing nonzero $Hom^0(E_j, E_i)$ for $i\neq j$ is then given by the directed graph in Figure~\ref{fig: 22}. Starting from vertex $0$ in $\mathcal{Q}$ (which corresponds to $\mathcal{O}_X( -1,-2,-3)$ and moving counter-clockwise, no arrows appear with vertex label of the head less than the vertex label of the tail.

\subsection{Polytope F.3D.0011} The polytope $P$ (denoted F.3D.0011 in polymake) has half-space representation

\[ \left[ \begin{matrix} 1 &0 &1 &1\\
1 &-1& 0& 0 \\
1 &0 &0 &1 \\
1 &0 &-1& 0 \\
1 &1 &0 &0 \\
1 &0 &0 &-1 \end{matrix} \right] \] 

and vertices given by the columns of 

\[ \left( \begin{matrix} -1 &1 &-1 &1 &-1 &1 &-1 &1 \\
0&0 &1 &1 &-2 &-2 &1 &1 \\
-1&-1& -1 &-1 &1 &1 &1 &1 \end{matrix} \right). \] 

Let $X$ denote the complete toric variety associated to $P$, with primitive ray generators of $\Sigma$ \[ \{  (-1,0,0), (1,0,0), (0,-1,0), (0,0,-1), (0,0,1 ), (0,1,1) \} \]

and presentation of the class group given by 

\[ \Z^6 \stackrel{ \left( \begin{matrix} 1 &1 &0 &0 &0 &0 \\ 
0 &0 &0 &1 &1 &0\\
0& 0 & 1 & 1& 0 &1 \end{matrix} \right) }{\longrightarrow} \Z^3. \] 

The Hanlon-Hicks-Lazarev resolution of the diagonal yields the free ranks (written as ungraded S-modules, for $S$ the homogeneous coordinate ring of $Y \cong X \times X$):\\

\[ 0 \rightarrow S^2 \rightarrow S^5 \rightarrow S^4 \rightarrow S^1 \rightarrow 0.\]

The collection of $8$ line bundles $\mathcal{E}$ which appear on the left-hand side are $\mathcal{O}(a_1,a_2,a_3)$ for $(a_1,a_2,a_3)$ in 

\[ \mathcal{E} = \{  (-1,-1,-2), (-1,-1,-1), (-1,0,-1), (0,-1,-2), (0,0,-1), (0,-1,-1), (-1,0,0), (0,0,0) \} \] 

gives a full, strong exceptional collection of line bundles. That is, $Hom_{D^b(X)}^\bullet(E_i, E_j)$ is concentrated in degree $0$ for all $i$ and $j$, with the rank of $Hom^0(E_j, E_i)$ given by the $(i,j)$ entry of the following matrix:

\[ F = \left( \begin{matrix} 
1 &0 &0 &0 &0 &0 &0 &0 \\
2 &1 &0 &0 &0 &0 &0 &0 \\
 3 &0 &1 &1 &0 &0 &0 &0 \\
 2& 0 &0 &1 &0 &0 &0 &0\\
 5&0 &2 &3 &1 &0 &0 &0\\
 4&2 &0 &2 &0 &1 &0 &0\\
 6&3 &2 &2 &0 &1 &1 &0\\
 10&5 &4 &6 &2 &3 &2 &1 \end{matrix} \right). \]

\begin{figure}[!htbp]
\includegraphics[width=.5\textwidth]{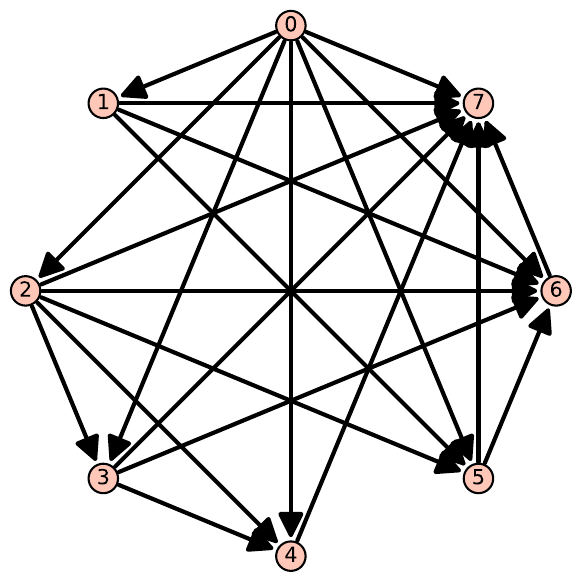}
 \caption{Quiver showing nonzero $Hom^0(E_i, E_j)$ for $i\neq j$ in $\mathcal{E}$. Indexing on vertices from $0$ to $7=|\mathcal{E}|-1$.} 
\label{fig: 23}
\end{figure}

The fact that $F$ is lower-triangular shows that $\mathcal{E}$ is exceptional. The quiver $\mathcal{Q}$ showing nonzero $Hom^0(E_j, E_i)$ for $i\neq j$ is then given by the directed graph in Figure~\ref{fig: 23}. Starting from vertex $0$ in $\mathcal{Q}$ (which corresponds to $\mathcal{O}_X( -1,-1,-2))$ and moving counter-clockwise, no arrows appear with vertex label of the head less than the vertex label of the tail.

\subsection{Polytope F.3D.0012} The polytope $P$ (denoted F.3D.0012 in polymake) has half-space representation 

\[ \left[ \begin{matrix}
1 &0 &1 &1 \\
1 &-1& 0& 0\\
1 &0 &0 &1\\
1 &0 &-1& 0\\
1 &0& 0 &-1\\
1& 1& 0& -1 \end{matrix} \right] \] 

with vertices given by the columns of 

\[ \left( \begin{matrix} -2 &1 &-2 &1 &0 &1 &0 &1 \\
0&0 &1 &1 &-2&-2&1&1\\
-1&-1 &-1 &-1 &1 &1 &1 &1 \end{matrix} \right). \] 

Let $X$ denote the complete toric variety associated to $P$, with primitive ray generators of $\Sigma$ \[ \{(-1,0,0),(0,-1,0),(0,0,-1),(1,0,-1),(0,0,1),(0,1,1)\} \] and presentation of the class group given by 

\[ \Z^6 \stackrel{ \left( \begin{matrix} -1 &0 &1 &-1 &0 &0 \\ 1 &0 &0 &1 &1 &0 \\ 1 &1 &0 &1 &0 &1 \end{matrix} \right) }{\longrightarrow } \Z^3.  \]  

The Hanlon-Hicks-Lazarev resolution of the diagonal yields the free ranks (written as ungraded $S$-modules, for $S$ the homogeneous coordinate ring of $Y \cong X \times X$):

\[ 0 \rightarrow S^4 \rightarrow S^9 \rightarrow S^6 \rightarrow S^1 \rightarrow 0. \] 

The collection of $8$ line bundles $\mathcal{E}$ which appear on the left-hand side are $\mathcal{O}(a_1, a_2, a_3)$ for $(a_1, a_2, a_3)$ appearing in 

\[ \mathcal{E} = \{ (1,-2,-3), (0,-1,-2), (1,-2,-2), (1,-1,-2), (0,-1,-1), (1,-1,-1), (0,0,-1), (0,0,0) \} \] 

is full, strong, and exceptional. That is, $Hom_{D^b(X)}^\bullet(E_i, E_j)$ is concentrated in degree $0$ for all $i$ and $j$, with the rank of $Hom^0(E_j, E_i)$ given by the $(i,j)$ entry of the following matrix:

\[  F  =  \left( \begin{matrix}  1&0 &0&0&0&0&0&0\\ 2&1&0&0&0&0&0&0 \\ 2&0&1&0&0&0&0&0 \\ 4&1&1&1&0&0&0&0 \\ 4&2&2&0&1&0&0&0 \\ 7&2&4&2&1&1&0&0 \\ 7&4&2&2&1&0&1&0 \\ 12&7&7&4&4&2&2&1 \end{matrix} \right). \] 

\begin{figure}[!htbp]
\includegraphics[width=.5\textwidth]{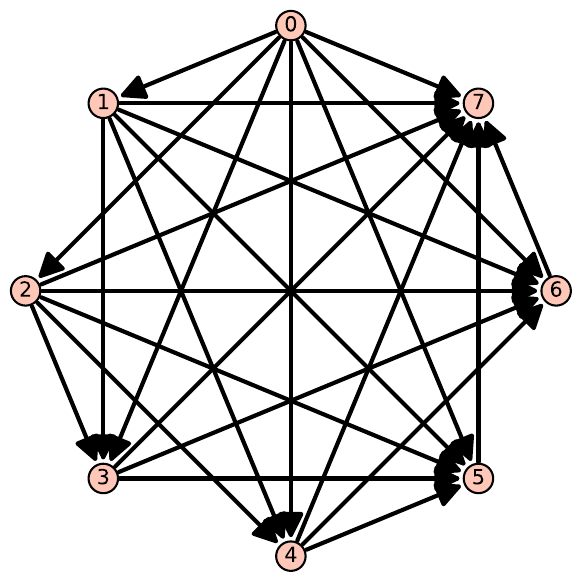}
 \caption{Quiver showing nonzero $Hom^0(E_i, E_j)$ for $i\neq j$ in $\mathcal{E}$. Indexing on vertices from $0$ to $7=|\mathcal{E}|-1$.} 
\label{fig: 24}
\end{figure}

The fact that $F$ is lower-triangular shows that $\mathcal{E}$ is exceptional. The quiver $\mathcal{Q}$ showing nonzero $Hom^0(E_j, E_i)$ for $i\neq j$ is then given by the directed graph in Figure~\ref{fig: 24}. Starting from vertex $0$ in $\mathcal{Q}$ (which corresponds to $\mathcal{O}_X( 1, -2, -3))$ and moving counter-clockwise, no arrows appear with vertex label of the head less than the vertex label of the tail.

\subsection{Polytope F.3D.0013} The polytope $P$ (denoted F.3D.0013 in polymake) has half-space representation

\[ \left[ \begin{matrix} 1& -1& 1& 0 \\
1 &1 &0 &1\\
1 &-1& 0& 0\\
1 &0 &-1& 0\\
1 &0 &0& -1   \end{matrix} \right] \] with vertices given by the columns of 

\[ \left( \begin{matrix} 
1&1&1&-2&-2&1\\
0&1&0&-3&1&1\\
-2&-2&1&1&1&1 \end{matrix} \right).   \] 

Let $X$ denote the complete toric variety associated to $P$, with primitive ray generators of $\Sigma$ \[    \]

and presentation of the class group given by 

\[ \Z^5 \stackrel{ \left( \begin{matrix} -1&1 &1 &0 &0 \\ 1 &0 &0 &1 &1 \end{matrix} \right) }{\longrightarrow} \Z^2. \] 

The Hanlon-Hicks-Lazarev resolution of the diagonal yields the free ranks (written as ungraded $S$-modules, for $S$ the homogeneous coordinate ring of $Y \cong X \times X$):

\[ 0 \rightarrow S^4 \rightarrow S^9 \rightarrow S^6 \rightarrow S^1 \rightarrow 0. \]

The collection of $6$ line bundles $\mathcal{E}$ which appear on the left-hand side are $\mathcal{O}(a_1, a_2)$ for $(a_1, a_2)$ appearing in 

\[ \mathcal{E} = \{ (-1,-2),(0,-2),(-1,-1),(0,-1),(-1,0),(0,0)\} \]

is full, strong, and exceptional. That is, $Hom_{D^b(X)}^\bullet(E_i, E_j)$ is concentrated in degree $0$ for all $i$ and $j$, with the rank of $Hom^0(E_j, E_i)$ given by the $(i,j)$ entry of the following matrix:

 \[ F = \left( \begin{matrix} 1 &0 &0 &0 &0 &0 \\
 2 &1 &0 &0 &0 &0\\
 4 &1 &1 &0 &0 &0\\
 7 &4 &2 &1 &0 &0\\
 10 &4 &4 &1 &1& 0\\
 16&10 &7 &4 &2 &1\end{matrix} \right). \] 

\begin{figure}[!htbp]
\includegraphics[width=.5\textwidth]{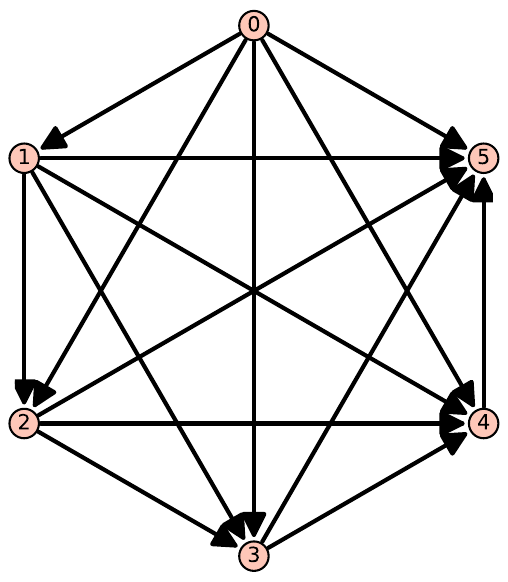}
 \caption{Quiver showing nonzero $Hom^0(E_i, E_j)$ for $i\neq j$ in $\mathcal{E}$. Indexing on vertices from $0$ to $5=|\mathcal{E}|-1$.} 
\label{fig: 25}
\end{figure}

The fact that $F$ is lower-triangular shows that $\mathcal{E}$ is exceptional. The quiver $\mathcal{Q}$ showing nonzero $Hom^0(E_j, E_i)$ for $i\neq j$ is then given by the directed graph in Figure~\ref{fig: 25}. Starting from vertex $0$ in $\mathcal{Q}$ (which corresponds to $\mathcal{O}_X( -1,-2))$ and moving counter-clockwise, no arrows appear with vertex label of the head less than the vertex label of the tail.

\subsection{Polytope F.3D.0014} The polytope $P$ (denoted F.3D.0014 in polymake) has half-space representation 

\[ \left[ \begin{matrix} 
1 &1 &1 &1 \\
1 &-1& 0& 0 \\
1 &0 &0 &1 \\
1& 0& -1& 0 \\
1 &0 &0 &-1 \end{matrix} \right] \] 

with vertices given by the columns of 

\[ \left( \begin{matrix} 1&-1&1&1&-3&1\\ -1&1&1&-3&1&1 \\-1&-1&-1&1&1&1 \end{matrix} \right).  \]

Let $X$ denote the complete toric variety associated to $P$, with primitive ray generators of $\Sigma$:

\[ \{  (-1,0,0), (0,-1,0), (0,0,-1), (0,0,1), (1,1,1)  \}      \]

and presentation of the class group given by 

\[ \Z^5 \stackrel{  \left( \begin{matrix} 0&0&1&1&0 \\ 1&1&1&0&1 \end{matrix}\right)}{ \longrightarrow} \Z^2. \] 

The Hanlon-Hicks-Lazarev resolution of the diagonal yields the free ranks (written as ungraded $S$-modules, for $S$ the homogeneous coordinate ring of $Y \cong X \times X$):

\[ 0 \rightarrow S^3 \rightarrow S^8 \rightarrow S^6 \rightarrow S^1 \rightarrow  0.\]

The collection of 6 line bundles $\mathcal{E}$ which appear on the left-hand side are $\mathcal{O}(a_1, a_2)$ for $(a_1, a_2)$ appearing in 

\[ \mathcal{E} = \{ (-1,-3), (-1,-2), (0, -2), (-1,-1), (0,-1), (0,0) \} \] 

is full, strong, and exceptional. That is, $Hom_{D^b(X)}^\bullet(E_i, E_j)$ is concentrated in degree $0$ for all $i$ and $j$, with the rank of $Hom^0(E_j, E_i)$ given by the $(i,j)$ entry of the following matrix:

\[  F  =  \left( \begin{matrix}  1 &0 &0 &0 &0 &0 \\ 3 &1 &0 &0 &0 &0 \\4 &1 &1 &0 &0 &0 \\ 6 &3 &0 &1 &0 &0 \\ 9 &4 &3 &1 &1 &0 \\ 16 &9 &6 &4 &3 &1 \end{matrix} \right). \] 

\begin{figure}[!htbp]
\includegraphics[width=.5\textwidth]{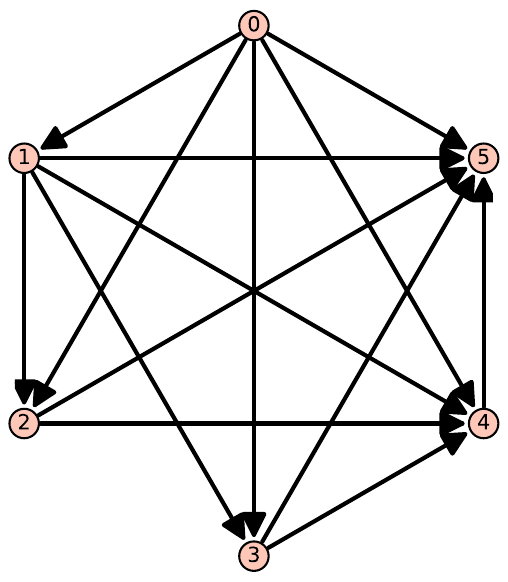}
 \caption{Quiver showing nonzero $Hom^0(E_i, E_j)$ for $i\neq j$ in $\mathcal{E}$. Indexing on vertices from $0$ to $5=|\mathcal{E}|-1$.} 
\label{fig: 26}
\end{figure}

The fact that $F$ is lower-triangular shows that $\mathcal{E}$ is exceptional. The quiver $\mathcal{Q}$ showing nonzero $Hom^0(E_j, E_i)$ for $i\neq j$ is then given by the directed graph in Figure~\ref{fig: 26}. Starting from vertex $0$ in $\mathcal{Q}$ (which corresponds to $\mathcal{O}_X(-1,-3))$ and moving counter-clockwise, no arrows appear with vertex label of the head less than the vertex label of the tail.

\subsection{Polytope F.3D.0015} The polytope $P$ (denoted F.3D.0015 in polymake) has half-space representation

\[ \left[ \begin{matrix} 
1 &0 &1& 0 \\
1 &-1& 0& 0 \\
1 &0 &0 &1 \\
1 &0 &-1& 0 \\
1 &1 &0 &0 \\
1 &0 &0& -1 
\end{matrix} \right] \] 

and vertices given by the columns of 

\[ \left( \begin{matrix} -1&1&-1&1&-1&1&-1&1 \\
-1&-1&1&1&-1&-1&1&1 \\
-1&-1&-1&-1&1&1&1&1 \end{matrix} \right).  \] 

Let $X$ denote the complete toric variety associated to $P$, with primitive ray generators of $\Sigma$ \[ \{ (-1,0,0), (1,0,0), (0,-1,0), (0,1,0), (0,0,-1),(0,0,1)  \} \]

and presentation of the class group given by 

\[ \Z^6 \stackrel{ \left(\begin{matrix} 1&1&0&0&0&0 \\ 0&0&1&1&0&0 \\0&0&0&0&1&1 \end{matrix} \right) } {\longrightarrow} \Z^3. \] 

The Hanlon-Hicks-Lazarev resolution of the diagonal yields the free ranks (written as ungraded $S$-modules, for $S$ the homogeneous coordinate ring of $Y\cong X \times X)$

\[ 0 \rightarrow S^1 \rightarrow S^3 \rightarrow S^3 \rightarrow S^1 \rightarrow  0. \]

The collection of 8 line bundles $\mathcal{E}$ which appear on the left-hand side are $\mathcal{O}(a_1, a_2, a_3)$ for $(a_1, a_2, a_3)$ appearing in \[ \mathcal{E} = \{ (-1,-1,-1), (-1,-1,0), (0,-1,-1), (-1,0,-1), (0,0,-1), (-1,0,0), (0,-1,0), (0,0,0) \} \]

is full, strong, and exceptional. That is, $Hom_{D^b(X)}^\bullet(E_i, E_j)$ is concentrated in degree $0$ for all $i$ and $j$, with the rank of $Hom^0(E_j, E_i)$ given by the $(i,j)$ entry of the following matrix:

\[ F = \left( \begin{matrix}      1&0&0&0&0&0&0&0\\ 2&1&0&0&0&0&0&0\\2&0&1&0&0&0&0&0\\ 2&0&0&1&0&0&0&0 \\4&0&2&2&1&0&0&0\\ 4&2&0&2&0&1&0&0 \\4&2&2&0&0&0&1&0 \\ 8&4&4&4&2&2&2&1                   \end{matrix} \right). \] 

\begin{figure}[!htbp]
\includegraphics[width=.5\textwidth]{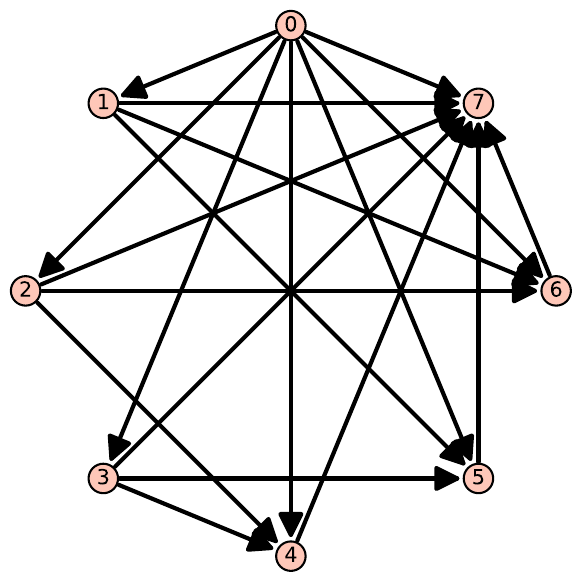}
 \caption{Quiver showing nonzero $Hom^0(E_i, E_j)$ for $i\neq j$ in $\mathcal{E}$. Indexing on vertices from $0$ to $7=|\mathcal{E}|-1$.} 
\label{fig: 27}
\end{figure}

The fact that $F$ is lower-triangular shows that $\mathcal{E}$ is exceptional. The quiver $\mathcal{Q}$ showing nonzero $Hom^0(E_j, E_i)$ for $i\neq j$ is then given by the directed graph in Figure~\ref{fig: 27}. Starting from vertex $0$ in $\mathcal{Q}$ (which corresponds to $\mathcal{O}_X( -1,-1,-1))$ and moving counter-clockwise, no arrows appear with vertex label of the head less than the vertex label of the tail.

\subsection{Polytope F.3D.0016} \label{sec:0016} The polytope $P$ (denoted F.3D.0016 in polymake) has half-space representation

\[ \left[ \begin{matrix} 1 &0 &1 &0\\
1& 1& 0& 1\\
1 &-1& 0& 0\\
1 &0 &-1& 0\\
1& 0& 0& -1 \end{matrix} \right] \] 

and vertices given by the columns of

\[ \left( \begin{matrix} 1&1 & -2 &1 &-2 &1 \\ -1 &1 &-1&-1&1&1\\-2&-2&1&1&1&1 \end{matrix} \right). \]

Let $X$ denote the complete toric variety associated to $P$, with primitive ray generators of $\Sigma$ \[ \{ (-1,0,0), (0,-1,0), (0,1,0), (0,0,-1), (1,0,1)   \} \]

and presentation of the class group given by 

\[ \Z^5 \stackrel{ \left( \begin{matrix} 0 &1 &1 &0 &0 \\ 1 &0 &0 &1 &1 \end{matrix} \right) }{\longrightarrow} \Z^2. \]

The Hanlon-Hicks-Lazarev resolution of the diagonal yields the free ranks (written as ungraded $S$-modules, for $S$ the homogeneous coordinate ring of $Y\cong X \times X)$

\[ 0 \rightarrow S^2 \rightarrow S^5 \rightarrow S^4 \rightarrow S^1 \rightarrow 0. \]

The collection of 6 line bundles $\mathcal{E}$ which appear on the left-hand side are $\mathcal{O}(a_1, a_2)$ for $(a_1, a_2)$ appearing in \[ \mathcal{E} = \{ (-1,-2), (-1,-1), (0,-2),(-1,0), (0,-1), (0,0) \} \] 

is full, strong, and exceptional. That is, $Hom_{D^b(X)}^\bullet(E_i, E_j)$ is concentrated in degree $0$ for all $i$ and $j$, with the rank of $Hom^0(E_j, E_i)$ given by the $(i,j)$ entry of the following matrix:

\[ F = \left( \begin{matrix}  1&0&0&0&0&0 \\ 3&1&0&0&0&0 \\ 2&0&1&0&0&0\\ 6&3&0&1&0&0 \\ 6&2&3&0&1&0\\ 12&6&6&2&3&1                    \end{matrix} \right). \] 

\begin{figure}[!htbp]
\includegraphics[width=.5\textwidth]{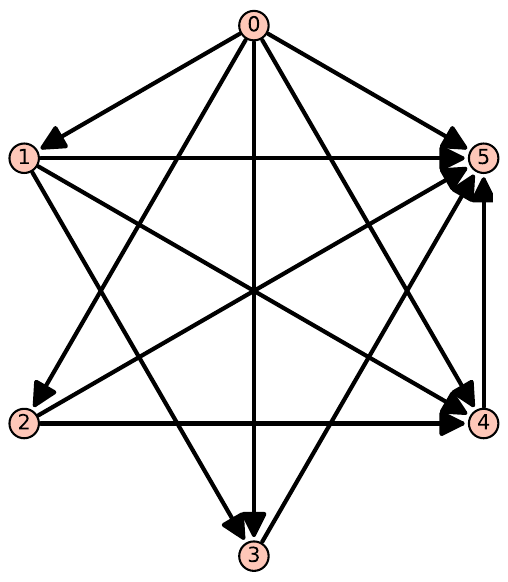}
 \caption{Quiver showing nonzero $Hom^0(E_i, E_j)$ for $i\neq j$ in $\mathcal{E}$. Indexing on vertices from $0$ to $5=|\mathcal{E}|-1$.} 
\label{fig: 28}
\end{figure}

The fact that $F$ is lower-triangular shows that $\mathcal{E}$ is exceptional. The quiver $\mathcal{Q}$ showing nonzero $Hom^0(E_j, E_i)$ for $i\neq j$ is then given by the directed graph in Figure~\ref{fig: 28}. Starting from vertex $0$ in $\mathcal{Q}$ (which corresponds to $\mathcal{O}_X( -1,-2))$ and moving counter-clockwise, no arrows appear with vertex label of the head less than the vertex label of the tail.

\subsection{Polytope F.3D.0017} \label{sec:0017}The polytope $P$ (denoted F.3D.0017 in polymake) has half-space representation

\[ \left[ \begin{matrix}  1& -1& 0& 0\\
1 &0 &-1& 0\\
1 &0 &0 &-1\\
1 &1& 1& 1     \end{matrix} \right] \]

with vertices given by the columns of

\[ \left( \begin{matrix}  1&1&-3&1\\1&-3&1&1 \\-3&1&1&1 \end{matrix}  \right). \] 

Let $X$ denote $\PP^3$, the complete toric variety associated to $P$, with primitive ray generators of $\Sigma$ \[ \{ (-1,0,0), (0,-1,0), (0,0,-1),(1,1,1)  \} \]

and presentation of the class group given by 

\[ \Z^4 \stackrel{ \left( \begin{matrix} 1 &1& 1& 1 \end{matrix} \right) }{\longrightarrow} \Z^1.  \]

The Hanlon-Hicks-Lazarev resolution of the diagonal yields the free ranks (written as ungraded $S$-modules, for $S$ the homogeneous coordinate ring of $Y\cong X \times X)$

\[ 0 \rightarrow S^3 \rightarrow S^8 \rightarrow S^6 \rightarrow S^1 \rightarrow 0. \] 

The collection of 4 line bundles $\mathcal{E}$ which appear on the left-hand side is the Beilinson collection

\[ \mathcal{E} = \{ \mathcal{O}(-3), \mathcal{O}(-2), \mathcal{O}(-1), \mathcal{O} \} \] 

which is full, strong, and exceptional. That is, $Hom_{D^b(X)}^\bullet(E_i, E_j)$ is concentrated in degree $0$ for all $i$ and $j$, with the rank of $Hom^0(E_j, E_i)$ given by the $(i,j)$ entry of the following matrix:

\[ F = \left( \begin{matrix} 1&0 &0 &0 \\ 4&1 &0 &0 \\ 10& 4&1 &0 \\20&10&4 &1  \end{matrix} \right). \] 

\begin{figure}[!htbp]
\includegraphics[width=.5\textwidth]{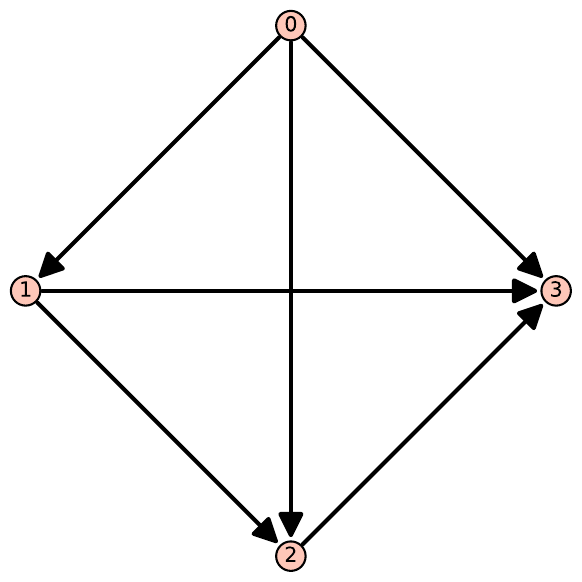}
 \caption{Quiver showing nonzero $Hom^0(E_i, E_j)$ for $i\neq j$ in $\mathcal{E}$. Indexing on vertices from $0$ to $3=|\mathcal{E}|-1$.} 
\label{fig: 29}
\end{figure}

The fact that $F$ is lower-triangular shows that $\mathcal{E}$ is exceptional. The quiver $\mathcal{Q}$ showing nonzero $Hom^0(E_j, E_i)$ for $i\neq j$ is then given by the directed graph in Figure~\ref{fig: 29}. Starting from vertex $0$ in $\mathcal{Q}$ (which corresponds to $\mathcal{O}_X( -3))$ and moving counter-clockwise, no arrows appear with vertex label of the head less than the vertex label of the tail.

\subsection{Discussion}
Theorem~\ref{thm: MainTheorem2} is proved in Sections~\ref{sec:6.1}--~\ref{sec:0017}, which does not require the use of Theorem~\ref{thm:5.1}. The polytopes F.3D.0001 and F.3D.0000 correspond to II(d) and III(k) in Watanabe--Watanabe's classification \cite{watanabe-watanabe}, respectively, which are $Z = \PP(\mathcal{O}_{\PP^2} \oplus \mathcal{O}_{\PP^2}(2) )$ and $W$, the blowup of $Z$ along a line, respectively. That a numerical condition for $Z$ and $W$ fails for the presence of a full, strong exceptional collection of line bundles is given in \cite{bondalOberwolfach}. We reiterate that this shows only that Bondal's construction does not yield full strong exceptional collections for these two threefolds. Nevertheless, every smooth toric Fano threefold over $\C$ has such a collection: Bernardi--Tirabassi constructed the strong collections \cite{bernardi2010derivedcategoriestoricfano}, and Uehara proved fullness \cite{uehara2015exceptionalcollectionstoricfano}.

\end{document}